\documentclass{article}
\usepackage[utf8]{inputenc}
\usepackage{amsmath,amsbsy,multirow,multicol,tabularx,booktabs}
\usepackage{graphicx,placeins}
\usepackage{color}
\usepackage[shortlabels]{enumitem}
\usepackage{a4wide}
\usepackage{amsthm}
\usepackage[psamsfonts]{amssymb}
\usepackage{latexsym}
\usepackage{oldlfont}
\usepackage{lmodern}

\newtheorem{conj}{Conjecture}
\newtheorem{thm}[conj]{\bf Theorem}

\newtheorem{prop}[conj]{\bf Proposition}
\newtheorem{lemma}[conj]{\bf Lemma}
\newtheorem{rem}[conj]{\bf Remark}

\newtheorem{example}[conj]{\bf Example}

\newcommand*{\bigdot}[1]{%
  \overset{\mbox{\large\bfseries .}}{#1}}
\newcommand*{\bigdott}[1]{%
  \overset{\mbox{\large\bfseries ..}}{#1}}

\makeatletter
\def\saveenum{\xdef\@savedenum{\the\c@enumi\relax}}
\def\resetenum{\global\c@enumi\@savedenum}
\makeatother

\usepackage[top=3cm,bottom=3cm,right=3.4cm,left=3.4cm]{geometry}

\def\bar{\overline}
\def\to{\rightarrow}

\def\Cc{\mbox{$\mathcal C$}}
\def\Ac{\mbox{$\mathcal A$}}
\def\Dc{\mbox{$\mathcal D$}}
\def\Pc{\mbox{$\mathcal P$}}
\def\Hc{\mbox{$\mathcal H$}}
\def\Fc{\mbox{$\mathcal G$}}
\def\Mc{\mbox{$\mathcal M$}}
\def\Xc{\mbox{$\mathcal X$}}

\def\Fc{\mbox{$\mathcal F$}}
\def\Nc{\mbox{$\mathcal N$}}
\def\Ic{\mbox{$\mathcal I$}}

\def\Tc{\mbox{$\mathcal T$}}
\def\Uc{\mbox{$\mathcal U$}}

\def\Xc{\mbox{$\mathcal X$}}
\def\Wc{\mbox{$\mathcal W$}}
\def\Sc{\mbox{$\mathcal S$}}
\def\AA{{\mathbb A}}

\def\BB{{\mathbb B}}

\def\GG{{\mathbb G}}
\def\RR{{\mathbb R}}

\def\WW{{\mathbb W}}

\def\EE{ {\rm I} \kern-.15em {\rm E} }
\def\PP{ {\rm I} \kern-.15em {\rm P} }
\def\E{ {\rm I} \kern-.15em {\rm E} }

\def\contiguous{\triangleleft\kern-.20em\triangleright}

\def\ii{{\mathbf 1}( }
\def\1{{\mathbf 1} }

\def\x{ {\bf x}}

\def\y{ {\bf y}}
\def\z{ {\bf z}}

\def\u{ {\bf u}}

\def\X{ {\bf X}}

\def\Z{ {\bf Z}}

\def\mds{\medskip}

\begin{document}

\title{About tests of the ``simplifying'' assumption for conditional copulas}

\author{Alexis Derumigny\thanks{ENSAE, 3 avenue Pierre-Larousse, 92245 Malakoff cedex, France. alexis.derumigny@ensae.fr}, Jean-David Fermanian\thanks{ENSAE, J120, 3 avenue Pierre-Larousse, 92245 Malakoff cedex, France. jean-david.fermanian@ensae.fr. This research has been supported by the Labex Ecodec.}}

\date{\today}

\maketitle

\abstract{We discuss the so-called ``simplifying assumption'' of conditional copulas in a general framework. We introduce several tests of the latter assumption for non- and semiparametric copula models. Some related test procedures based on conditioning subsets instead of point-wise events are proposed. The limiting distribution of such test statistics under the null are approximated by several bootstrap schemes, most of them being new. We prove the validity of a particular semiparametric bootstrap scheme.
Some simulations illustrate the relevance of our results.
}

\mds


\mds

{\bf MCS:} 62G05, 62G08, 62G09.

\section{Introduction}
In statistical modelling and applied science more generally, it is very common to distinguish two subsets of variables: a random vector of interest (also called explained/exogenous variables) and a vector of covariates (explanatory/endogenous variables). The objective is to predict the law of the former vector given the latter vector belongs to some subset, possibly a singleton. This basic idea constitutes the first step towards forecasting some important statistical sub-products as conditional means, quantiles, volatilities, etc.
Formally, consider a $d$-dimensional random vector $\X$. We are faced with two random sub-vectors $\X_I$ and $\X_J$, s.t. $\X=(\X_I,\X_J)$, $I\cup J =\{1,\ldots,d\}$, $I\cap J=\emptyset$, and our models of interest specify the conditional law of $\X_I$ knowing $\X_J=\x_J$ or knowing $\X_J\in A_J$ for some subset $A_J\subset \RR^{|J|}$.
We use the standard notations for vectors: for any set of indices $I$,
$\x_I$ means the $|I|$-dimensional vector whose arguments are the $x_k$, $k\in I$. For convenience and without a loss of generality, we will set $I=\{1,\ldots,p\}$ and $ J=\{p+1,\ldots,d\}$.

\mds

Besides, the problem of dependence among the components of $d$-dimensional random vectors has been extensively studied in the academic literature and among practitioners in a lot of different fields. The raise of copulas for more than twenty years illustrates the need of flexible and realistic multivariate models and tools.
When covariates are present and with our notations, the challenge is to study the dependence among the components of $\X_I$ given $\X_J$.
Logically, the concept of conditional copulas has emerged. 
{\color{black} First introduced for pointwise (atomic) conditioning events by Patton (2006a, 2006b), the definition has been generalized in Fermanian and Wegkamp (2012) for arbitrary measurable conditioning subsets.
In this paper, we rely on the following definition:} for any borel subset $A_J \subset \RR^{d-p}$, a conditional copula of $\X_I$ given $(\X_J \in A_J)$ is denoted by $C_{I|J}^{A_J}(\cdot | \X_J \in A_J)$.
This is the cdf of the random vector $(F_{1|J}(X_{1}|\X_J\in A_J),\ldots, F_{p|J}(X_{p}|\X_J\in A_J))$ given $(\X_J \in A_J)$.
Here, $F_{k|J}(\cdot | \X_J\in A_J)$ denotes the conditional law of $X_k$ knowing $\X_J\in A_J$, $k=1,\ldots,p$.
The latter conditional distributions will be assumed continuous in this paper, implying
the existence and uniqueness of $C_{I|J}^{A_J}$ (Sklar's theorem).
In other words, for any $\x_I\in \RR^p$,
\begin{equation*}
    \PP \left(\X_I \leq \x_I | \X_J \in A_J \right) = C_{I|J}^{A_J}\Big( F_{1|J}(x_1 | \X_J\in A_J),\ldots,  F_{p|J}(x_p | \X_J \in A_J) \, \Big| \,  \X_J \in A_J \Big).
    \label{cond_cop_def_0}
\end{equation*}
{\color{black} Note that the influence of $A_J$ on $C_{I|J}^{A_J}$ is twofold: when $A_J$ changes, the conditioning event $(\X_J \in A_J)$ changes, but the conditioned random vector $(F_{1|J}(X_{1}|\X_J\in A_J),\ldots, F_{p|J}(X_{p}|\X_J\in A_J))$ changes too.}

\mds

In particular, when the conditioning events are reduced to singletons, we get that the conditional copula of $\X_I$ knowing $\X_J=\x_J$ is a cdf $C_{I|J}(\cdot | \X_J=\x_J)$ on $[0,1]^{p}$ s.t., for every $\x_I\in \RR^{p}$,
\begin{equation*}
    \PP \left(\X_I \leq \x_I | \X_J=\x_J \right)
    = C_{I|J} \left( F_{1|J}(x_1 | \X_J=\x_J), \ldots, F_{p|J}(x_p | \X_J=\x_J) \, | \,  \X_J=\x_J \right).
    \label{cond_cop_def}
\end{equation*}

\mds

With generalized inverse functions, an equivalent definition of a conditional copula is as follows:
\begin{equation*}
    \label{cond_cop_def_geninverse}
    C_{I|J}\left( \u_I \, | \,  \X_J=\x_J \right) =
    F_{I|J} \big( F^{-}_{1|J}(u_1 | \X_J=\x_J),\ldots,F^{-}_{p|J}(u_p | \X_J=\x_J) | \X_J = \x_J \big),
\end{equation*}
for every $\u_I$ and $\x_J$, setting $F_{I|J}(\x_I | \X_J=\x_J):=\PP \left(\X_I \leq \x_I | \X_J = \x_J \right)$.

\mds

Most often, the dependence of $C_{I|J}(\cdot | \X_J=\x_J)$ w.r.t. to $\x_J$ is a source of significant complexities, in terms of model specification and inference.
Therefore, most authors assume that the following ``simplifying assumption'' is fulfilled.

\mds

{\it Simplifying assumption} $(\Hc_0)$: the conditional copula $C_{I|J}(\cdot | \X_J=\x_J)$ does not depend on $\x_J$, i.e., for every $\u_I\in [0,1]^p$, the function
$ \x_J\in \RR^{d-p}\mapsto C_{I|J}(\u_I | \X_J=\x_J)$ is a constant function (that depends on $\u_I$).

\mds

Under the simplifying assumption, we will set
$C_{I|J}(\u_I | \X_J=\x_J) =: C_{s,I|J}(\u_I)$.
The latter identity means that the dependence on $\X_J$ across the components of $\X_I$ is passing only through their conditional margins. Note that $C_{s,I|J}$ is different from the usual copula of $\X_I$:
{\color{black} $C_I(\cdot)$ is always the cdf of the vector $(F_1(X_1), \dots, F_p(X_p))$ whereas, under $\Hc_0$, $C_{s,I|J}$ is the cdf of the vector $\Z_{I|J}:=(F_{1|J}(X_1|X_J), \dots, F_{p|J}(X_p|X_J))$ (see Proposition~\ref{prop_indep_H0} below).
Note that the latter copula is identical to the partial copula introduced by Bergsma (2011), and recently studied by Gijbels et al. (2015b), Spanhel and Kurz (2015) in particular.
Such a partial copula $C_{I|J}^P$ can always be defined (whether $\Hc_0$ is satisfied or not) as the cdf of $\Z_{I|J}$, and it satisfies an interesting ``averaging'' property: $C_{I|J}^P(\u_I) :=
\int_{\RR^{d-p}} C_{I|J}(\u_I|\X_J=\x_J) dP_J(\x_J)$.}

\mds

\begin{rem}
    \label{ex_simple}
    The simplifying assumption $\Hc_0$ {\it does not imply} that
    $C_{s,I|J}(\cdot)$ is $C_I(\cdot)$, the usual copula of $\X_{I}$.
    This can be checked with a simple example: let $\X=(X_1,X_2,X_3)$ be a trivariate random vector s.t., given $X_3$, $X_1\sim \Nc(X_3,1)$ and $X_2\sim \Nc(X_3,1)$. Moreover, $X_1$ and
    $X_2$ are independent given $X_3$. The latter variable may be $\Nc(0,1)$, to fix the ideas. Obviously, with our notations, $I=\{1,2\}$, $J=\{3\}$, $d=3$ and $p=2$. Therefore, for any couple $(u_1,u_2)\in [0,1]^2$ and any real number $x_3$, $C_{1,2|3}(u_1,u_2 | x_3)=u_1 u_2$ and does not depend on $x_3$. Assumption $\Hc_0$ is then satisfied. But the copula of $(X_1,X_2)$ is not the independence copula, simply because $X_1$ and $X_2$ are not independent.
\end{rem}

\mds

Basically, it is far from obvious to specify and estimate relevant conditional copula models in practice, especially when the conditioning and/or conditioned variables are numerous.
The simplifying assumption is particularly relevant with vine models (Aas et al. 2009, among others).
Indeed, to build vines from a $d$-dimensional random vector $\X$, it is necessary to consider sequences of conditional bivariate copulas $C_{I|J}$, where $I=\{i_1,i_2\}$ is a couple of indices in $\{1,\ldots,d\}$, $J\subset \{1,\ldots,d\}$, $I \cap J = \emptyset$, and $(i_1,i_2|J)$ is a node of the vine.
In other words, a bivariate conditional copula is needed at every node of any vine, and the sizes of the conditioning subsets of variables are increasing along the vine.
Without additional assumptions, the modelling task becomes rapidly very cumbersome (inference and estimation by maximum likelihood).
Therefore, most authors adopt the simplifying assumption $\Hc_0$ at every node of the vine.
Note that the curse of dimensionality still apparently remains because conditional marginal cdfs $F_{k|J}(\cdot |\X_J)$ are invoked with different subsets $J$ of increasing sizes. But this curse can be avoided by calling recursively the non-parametric copulas that have been estimated before (see Nagler and Czado, 2015).

\mds

Nonetheless, the simplifying assumption has appeared to be rather restrictive, even if it may be
seen as acceptable for practical reasons and in particular situations. The debate between pro and cons of the simplifying assumption is still largely open, particularly
when it is called in some vine models. On one side, Hob\ae k-Haff et al. (2010) affirm that this simplifying assumption is not only required for
fast, flexible, and robust inference, but that it provides ‘‘a rather good approximation, even when the simplifying assumption
is far from being fulfilled by the actual model’’. On the other side, Acar et al. (2012) maintain that ``this view is too optimistic''.
They propose a visual test of $\Hc_0$ when $d=3$ and in a parametric framework.
Their technique was based on local linear approximations and sequential likelihood maximizations. They illustrate the limitations of $\Hc_0$ by simulation and through real datasets. They note that ``an uncritical use of the simplifying assumption may be misleading''. Nonetheless, they do not provide formal test procedures.
Beside, Acar et al. (2013) have proposed a formal likelihood test of the simplifying assumption but when the conditional marginal distributions are known, a rather restrictive situation.
Some authors have exhibited classes of parametric distributions for which $\Hc_0$ is satisfied: see Hob\ae k-Haff et al. (2010), significantly extended by St\"{o}ber et al. (2013). Nonetheless, such families are rather strongly constrained. Therefore, these two papers propose to approximate some conditional copula models by others for which the simplifying assumption is true.
This idea has been developed in Spanhel and Kurz (2015) in a vine framework, because
they recognize that ``it is very unlikely that the unknown data generating process satisfies the simplifying assumption in a strict mathematical sense.''

\mds

Therefore, there is a need for formal universal tests of the simplifying assumption. It is likely that the latter assumption is acceptable in some circumstances, whereas it is too rough in others.
This means, for given subsets of indices $I$ and $J$,
we would like to test
\begin{equation*}
    \Hc_0: C_{I|J}(\cdot| \X_J=\x_J) \;\text{does not depend on } \x_J,
    \label{atester}
\end{equation*}
against that opposite assumption.
Hereafter, we will propose several test statistics of $\Hc_0$, possibly assuming that the conditional copula belongs to some parametric family.

\mds

Note that several papers have already proposed estimators of conditional copula. Veraverbeke et al. (2011), Gijbels et al. (2011) and
Fermanian and Wegkamp (2012) have studied some nonparametric kernel based estimators.
Craiu and Sabeti (2012), Sabeti, Wei and Craiu (2014) studied bayesian additive models of conditional copulas.
Recently, Schellhase and Spanhel (2016) invoke B-splines to manage vectors of conditioning variables.
In a semiparametric framework, i.e. assuming an underlying parametric family of conditional copulas, numerous models and estimators have been proposed, notably
Acar et al. (2011), Abegaz et al. (2012), Fermanian and Lopez (2015) (single-index type models), Vatter and Chavez-Demoulin (2015) (additive models), among others.
But only a few of these papers have a focus on testing the simplifying assumption $\Hc_0$ specifically, although convergence of the proposed estimators
are necessary to lead such a task in theory. Actually, some tests of $\Hc_0$ is invoked ``in passing'' in these papers as potential applications,
but without a general approach and/or without some guidelines to evaluate p-values in practice. As an exception, in a very recent paper, Gijbels et al. (2016) have tackled the simplifying assumption directly through comparisons between conditional and unconditional Kendall's tau.

\mds

\begin{example}
    \label{example:Cases_of_SA}
    To illustrate the problem, let us consider a simple example of $\Hc_0$ in dimension $3$.
    Assume that $p=2$ and $d=3$. For simplicity, let us assume that $(X_1, X_2)$ follows a Gaussian distribution conditionally on $X_3$, that is :
    \begin{equation}
        \left( \begin{matrix} X_1 \\ X_2 \end{matrix} \right)
        \Big| X_3 = x_3 \sim \Nc
        \left( \left( \begin{matrix}
        \mu_1 (x_3) \\ \mu_2 (x_3) \end{matrix} \right) \, , \,
        \left( \begin{matrix}
            \sigma_1^2 (x_3) &
            \rho (x_3) \sigma_1(x_3) \sigma_2 (x_3) \\
            \rho (x_3) \sigma_1(x_3) \sigma_2 (x_3)
            & \sigma_2^2  (x_3)
        \end{matrix} \right)
        \right).
        \label{model_Gaussian:Cases_of_SA}
    \end{equation}
    Obviously, $\alpha(\cdot) := (\mu_1 , \mu_2 , \sigma_1 , \sigma_2)(\cdot)$ is a parameter that only affects the conditional margins.
    Moreover, the conditional copula of $(X_1,X_2)$ given $X_3=x_3$ is gaussian with the parameter $\rho(x_3)$.
    Six possible cases can then be distinguished:
    \begin{enumerate}[a.]
        \item All variables are mutually independent.
        \item $(X_1,X_2)$ is independent of $X_3$, but $X_1$ and $X_2$ are not independent.
        \item $X_1$ and $X_2$ are both marginally independent of $X_3$, but the conditional copula of $X_1$ and $X_2$ depends on $X_3$.

        \item $X_1$ (or $X_2$) and $X_3$ are not independent but $X_1$ and $X_2$ are independent conditionally given $X_3$.
        \item $X_1$ (or $X_2$) and $X_3$ are not independent but the conditional copula of $X_1$ and $X_2$ is independent of $X_3$.
        \item $X_1$ (or $X_2$) and $X_3$ are not independent and the conditional copula of $X_1$ and $X_2$ is dependent of $X_3$.
    \end{enumerate}
    %
    These six cases are summarized in the following table:

    \begin{center}
    \begin{tabular}{|c|c|c|c|}
        \hline
        & $\rho(\cdot) = 0$ & $\rho(\cdot) = \rho_0$ & $\rho(\cdot)$ is not constant \\
        \hline
        $\alpha(\cdot) = \alpha_0$
        & a & b & c \\
        \hline
        $\alpha(\cdot)$ is not constant
        & d & e & f \\
        \hline
    \end{tabular}
    \end{center}
    {\color{black} In the conditional Gaussian model (\ref{model_Gaussian:Cases_of_SA}), the simplifying assumption $\Hc_0$ consists in assuming that we live in one of the cases $\{a,b,d,e\}$, whereas the alternative cases are $c$ and $f$. In this model, the conditional copula is entirely determined by the conditional correlation. Note that, in other models, the conditional correlation can vary only because of the conditional margins, while the conditioning copula stay constant: see Property 8 of Spanhel and Kurz (2015).}
\end{example}

Note that, in general, there is no reason why the conditional margins would be constant in the conditioning variable (and in most applications, they are not).
Nevertheless, if we knew the marginal cdfs' were constant with respect to the conditioning variable, then the test of $\Hc_0$ (i.e. b against c) would become a classical test of independence between $\X_I$ and $\X_J$.

\mds

Testing $\Hc_0$ is closely linked to the $m$-sample copula problem, for which
we have $m$ different and independent samples of a $p$-dimensional variable $\X_I = (X_1, \dots, X_p)$. In each sample $k$, the observations are i.i.d., with
their own marginal laws and their own copula $C_{I,k}$.
The $m$-sample copula problem consists on testing whether the $m$ latter copulas $C_{I,k}$ are equal.
Note that we could merge all samples into a single one, and create discrete variables $Y_i$ that are equal to $k$ when $i$ lies in the sample $k$.
Therefore, the $m$-sample copula problem is formally equivalent to testing $\Hc_0$ with the conditioning variable $\X_J:=Y$.

\mds

Conversely, assume we have defined a partition $\{A_{1,J} , \dots , A_{m,J}\}$ of $\RR^{d-p}$ composed of borelian subsets such that $\PP(\X_{J} \in A_{k,J} ) > 0$ for all $k=1, \dots, m$,
and we want to test
$$\bar \Hc_0: k \in \{1, \dots, m\} \mapsto
C_{I|J}^{A_{k,J}} (\, \cdot \, | \X_J \in A_{k,J})
\, \text{does not depend on} \; k.$$
Then, divide the sample in $m$ different sub-samples, where any sub-sample $k$ contains the observations for which the conditioning variable belongs to $A_{k,J}$.
Then, $\bar \Hc_0$ is equivalent to a $m$-sample copula problem.
Note that $\bar\Hc_0$ looks like a ``consequence'' of $\Hc_0$ when it is not the case in general (see Section~\ref{Link_SA}), for continuous $\X_J$ variables.

\mds

Nonetheless, $\bar\Hc_0$ conveys the same intuition as $\Hc_0$. Since it can be led more easily in practice (no smoothing is required), some researchers could prefer the former assumption than the latter. That is why it will be discussed hereafter.
Note that the 2-sample copula problem has already been addressed by Rémillard and Scaillet (2009), and the $m$-sample by Bouzebda et al. (2011). However, both paper are designed only in a nonparametric framework, and these authors have not noticed the connection with the simplifying assumption.

\mds

The goal of the paper is threefold: first, to write a ``state-of-the art'' of the simplifying assumption problem; second to propose some ``reasonable''
test statistics of the simplifying assumption in different contexts; third, to introduce a new approach of the latter problem, through ``box-related'' zero assumptions
and some associated test statistics. Since it is impossible to state the theoretical properties of all these test statistics, we will rely on ``ad-hoc arguments'' to
convince the reader they are relevant, without trying
to establish specific results.
Globally, this paper can be considered also as a work program around the simplifying assumption $\Hc_0$ for the next years.

\mds

In Section \ref{Tests_SA}, we introduce different ways of testing $\Hc_0$.  We propose different test statistics under a fully nonparametric perspective,
i.e. when $C_{I|J}$ is not supposed to belong into a particular parametric copula family, through some comparisons between empirical cdfs' in Subsection~\ref{Bruteforce_SA},
or by invoking a particular independence property in Subsection~\ref{IndepProp_SA}. In Subsection~\ref{ParApproach_SA}, new tools are needed if we assume underlying parametric copulas.
To evaluate the limiting distributions of such tests, we propose several bootstrap techniques (Subsection \ref{Boot_SA}).
Section~\ref{Boxes} is related to testing $\bar \Hc_0$.
In Subsection \ref{Link_SA}, we detail the relations between $\Hc_0$ and $\bar \Hc_0$.
Then, we provide tests statistics of $\bar \Hc_0$ for both the nonparametric (Subsection \ref{NPApproach_m_SA}) and the parametric framework (Subsection \ref{ParApproach_m_SA}), as well as bootstrap methods (Subsection \ref{Boot_m_SA}). In particular, we prove the validity of the so-called ``parametric independent'' bootstrap when testing $\bar\Hc_0$.
The performances of the latter tests are assessed and compared by simulation in Section \ref{NumericalApplications}.
{\color{black}{{A table of notations is available in Appendix \ref{Section_notations} and} some of the proofs are collected in Appendix \ref{Section_Proofs}.}}

\mds

\section{Tests of the simplifying assumption}
\label{Tests_SA}

\subsection{``Brute-force'' tests of the simplifying assumption}
\label{Bruteforce_SA}

A first natural idea is to build a test of $\Hc_0$ based on a comparison between some estimates of the conditional copula $C_{I|J}$ with and without the simplifying assumption, for different conditioning events.
Such estimates will be called $\hat C_{I|J}$ and $\hat C_{s,I|J}$ respectively.
Then, introducing some distance $\Dc$ between conditional distributions, a test can be based on the statistics
$\Dc(\hat C_{I|J} , \hat C_{s,I|J})$.
Following most authors, we immediately think of Kolmogorov-Smirnov-type statistics
\begin{equation}
    \label{Tc0KS}
    \Tc^0_{KS,n}:=\|\hat C_{I|J} -\hat C_{s,I|J}\|_{\infty} =
    \sup_{\u_I \in [0,1]^p} \sup_{\x_J \in \RR^{d-p}} |\hat C_{I|J}(\u_I | \x_J) -\hat C_{s,I|J}(\u_I) |,
\end{equation}
or Cramer von-Mises-type test statistics
\begin{equation}
    \label{TcOCvM}
    \Tc^0_{CvM,n}:=\int \left( \hat C_{I|J}(\u_I | \x_J) -\hat C_{s,I|J}(\u_I) \right)^2\, w(d\u_I,d\x_J),
\end{equation}
for some weight function of bounded variation $w$, that could be chosen as random (see below).

\mds

To evaluate $\hat C_{I|J}$, we propose to invoke the nonparametric estimator of conditional copulas proposed by Fermanian and Wegkamp (2012).
Alternative kernel-based estimators of conditional copulas can be found in Gijbels et al. (2011), for instance.

\medskip

Let us start with an iid $d$-dimensional sample $(\X_i)_{i=1,\ldots,n}$.
Let $\hat F_k$ be the marginal empirical distribution function of $X_k$,
based on the sample $(X_{1,k},\ldots,X_{n,k})$, for any $k=1,\ldots,d$.
Our estimator of $C_{I|J}$ will be defined as
\begin{equation*}
    \hat C_{I|J} (\u_I| \X_J=\x_J) := \hat F_{I|J}
    \left( \hat F_{1|J}^{-}(u_1 | \X_J=\x_J)
    , \dots, \hat F_{p|J}^{-} (u_p | \X_J=\x_J) | \X_J=\x_J \right),
    \label{hatC_IJ}
\end{equation*}
\begin{equation}
    \hat F_{I|J} (\x_I | \X_J=\x_J) := \frac{1}{n}
    \sum_{i=1}^n K_n(\X_{i,J}, \x_J)
    \ii  \X_{i,I} \leq \x_I),
\label{hatF_I_J}
\end{equation}
where
$$ { \color{black} K_n(\X_{i,J}, \x_J) }
:= K_h \left(
\hat F_{p+1}(X_{i,p+1}) - \hat F_{p+1}(x_{p+1}) ,
\dots, \hat F_{d}(X_{i,d}) - \hat F_{d}(x_{d}) \right),$$
$$ K_h(\x_J) :=
h^{-(d-p)} K\left(x_{p+1}/h,\ldots,x_{d}/h\right),$$
and $K$ is a $(d-p)$-dimensional kernel.
Obviously, for $k\in I$, we have introduced some estimates of the marginal conditional cdfs' similarly:
\begin{equation}
    \hat F_{k|J} (x | \X_J=\x_J) :=
    {\color{black}
    \frac{\sum_{i=1}^n K_n(\X_{i,J}, \x_J)
    \ii  \X_{i,I}\leq \x_I)}{\sum_{j=1}^n K_n(\X_{j,J}, \x_J)}
    }\cdot
    \label{hatF_k_J}
\end{equation}

\medskip

Obviously, $h=h(n)$ is the term of a usual bandwidth sequence, where $h(n)\rightarrow 0$ when $n$ tends to the infinity.
Since $\hat F_{I|J}$ is a nearest-neighbors estimator, it does not necessitate a fine-tuning of local bandwidths (except for those values $\x_J$ s.t. $F_J(\x_J)$ is close to one or zero), contrary to more usual Nadaraya-Watson techniques.
In other terms, a single convenient choice of $h$ would provide ``satisfying'' estimates of $\hat C_{I|J}(\x_I | \X_J=\x_J)$ for most values of $\x$.
{\color{black} For practical reasons, it is important that $\hat F_{k|J}(x_k|\x_{J})$ belongs to $[0,1]$ and that $\hat F_{k|J}(\cdot|\x_{J})$ is a true distribution. This is the reason why we use a normalized version for the estimator of the conditional marginal cdfs.}

\medskip

To calculate the latter statistics~(\ref{Tc0KS}) and~(\ref{TcOCvM}), it is necessary to provide an estimate of the underlying conditional copula under $\Hc_0$.
This could be done naively by particularizing a point $\x_J^*\in \RR^{d-p}$ and by setting $\hat C_{s,I|J}^{(1)} (\cdot) := \hat C_{I|J}(\cdot | \X_J=\x_J^*)$.
Since the choice of $\x_J^*$ is too arbitrary, an alternative could be to set
\begin{equation*}
    \hat C_{s,I|J}^{(2)} (\cdot) :=
    \int \hat C_{I|J}(\cdot | \X_J=\x_J)\, w(d\x_J),
\end{equation*}
for some function $w$ that is of bounded variation, and $\int w(d\x_J)=1$. Unfortunately, the latter choice induce $(d-p)$-dimensional integration procedures, that becomes
a numerical problem rapidly when $d-p$ is larger than three.

\mds

Therefore, let us randomize the ``weight'' functions $w$, to avoid multiple integrations. For instance, choose the empirical distribution of $\X_J$ as $w$, providing
\begin{equation}
    \hat C_{s,I|J}^{(3)} (\cdot):=\int \hat C_{I|J}(\cdot | \X_J=\x_J)\, \hat F_J(d\x_J) = \frac{1}{n}\sum_{i=1}^n \hat C_{I|J}(\cdot | \X_J=\X_{i,J}).
    \label{estimator_meanCI_J}
\end{equation}

\mds

An even simpler estimate of $C_{s,I|J}$, the conditional copula of $\X_I$ given $\X_J$ under the simplifying assumption, can be obtained by noting that, under $\Hc_0$, $C_{s,I|J}$ is the joint law of $\Z_{I|J}:=(F_1(X_1|\X_J),\ldots,F_p(X_p|\X_J))$ (see Property~\ref{prop_indep_H0} below). Therefore, it is tempting to estimate $C_{s,I|J}(\u_I)$ by
\begin{equation}
    \hat C_{s,I|J}^{(4)} (\u_I) :=
    \frac{1}{n} \sum_{i=1}^n
    \1 \left( \hat F_{1|J}(X_{i,1}|\X_{i,J}) \leq u_1 , \dots, \hat F_{p|J}(X_{i,p}|\X_{i,J}) \leq u_p \right),
    \label{estimator_cond_cdf}
\end{equation}
when $\u_I\in [0,1]^p$,
for some consistent estimates $\hat F_{k|J}(x_k|\x_{J})$ of $F_{k|J}(x_k |\x_J)$.
A similar estimator has been promoted and studied in Gijbels et al. (2015a) or in Portier and Segers (2015), but they have considered the empirical copula associated to the
pseudo sample $((\hat F_1(X_{i1}|\X_{iJ}),\ldots,\hat F_p(X_{ip}|\X_{iJ})))_{i=1,\ldots,n}$ instead of its empirical cdf. It will be called $\hat C_{s,I|J}^{(5)}$.
Hereafter, we will denote $\hat C_{s,I|J}$ one of the ``averaged'' estimators $\hat C_{s,I|J}^{(k)}$, $k>1$ and we can forget the naive pointwise estimator $\hat C_{s,I|J}^{(1)}$.
Therefore, under some conditions of regularity, we guess that our estimators $\hat C_{s,I|J}(\u_I)$ of the conditional copula under $\Hc_0$ will be $\sqrt{n}$-consistent and asymptotically normal. It has been proved for $C^{(5)}_{s,I|J}$ in Gijbels et al. (2015a) or in Portier and Segers (2015), as a byproduct of the weak convergence of the associated process.

\mds

Under $\Hc_0$, we would like that the previous test statistics $\Tc^0_{KS,n}$ or $\Tc^0_{CvM,n}$ are convergent.
Typically, such a property is given as a sub-product by the weak convergence of a relevant empirical process, here
$(\u_I,\x_J)\in [0,1]^{p} \times \RR^{d-p} \mapsto \sqrt{nh_n^{d-p}}(\hat C_{I|J} - C_{I|J})(\u_I | \x_J)$.
Unfortunately, this will not be the case in general seing the previous process as a function indexed by $\x_J$, at least for wide ranges of bandwidths.
Due to the difficulty of checking the tightness of the process indexed by $\x_J$, some alternative techniques may be required
as Gaussian approximations (see Chernozhukov et al. 2014, e.g.). Nonetheless, they would lead us far beyond the scope of this paper.
Therefore, we simply propose to slightly modify the latter test statistics, to manage only a {\it fixed} set of arguments $\x_J$. For instance, in the case of
the Kolmogorov-Smirnov-type test, consider a simple grid $\chi_J:=\{\x_{1,J},\ldots, \x_{m,J}\}$, and the modified test statistics
\begin{equation*}
    \label{Tc0KSm}
    \Tc^{0,m}_{KS,n}:=
    \sup_{\u_I \in [0,1]^p} \sup_{\x_J \in \chi_J} |\hat C_{I|J}(\u_I | \x_J) -\hat C_{s,I|J}(\u_I) |.
\end{equation*}
In the case of the Cramer von-Mises-type test, we can approximate any integral by finite sums, possibly after a change of variable to manage a compactly supported integrand.
Actually, this is how they are calculated in practice!
For instance, invoking Gaussian quadratures, the modified statistics would be
\begin{equation}
    \label{TcOCvMm}
    \Tc^{0,m}_{CvM,n}:=\sum_{j=1}^m \omega_j  \left( \hat C_{I|J}(\u_{j,I} | \x_{j,J}) -\hat C_{s,I|J}(\u_{j,I}) \right)^2,
\end{equation}
for some conveniently chosen constants $\omega_j$, $j=1,\ldots,m$. Note that the numerical evaluation of $\hat C_{I|J}$ is relatively costly. Since
quadrature techniques require a lot less points $m$ than ``brute-force'' equally spaced grids (in dimension $d$, here), they have to be preferred most often.

\mds

Therefore, at least for such modified test statistics, we can insure the tests are convergent.
Indeed, under some conditions of regularity, it can be proved that $\hat C_{I|J}(\u_I| \X_J=\x_J) $ is consistent and asymptotically normal,
for every choice of $\u_I$ and $\x_J$ (see Fermanian and Wegkamp, 2012).
And a relatively straightforward extension of their Corollary 1 would provide that,
under ${\mathcal H}_0$ and for all $\Uc := (\u_{I,1},\ldots,\u_{I,q})\in [0,1]^{p(q+r)}$
and $\Xc := (\x_{J,1},\ldots,\x_{J,q})\in \RR^{(d-p)q}$,
\begin{eqnarray*}
    \lefteqn{
     \left\{
    \sqrt{nh_n^{d-p}}(\hat C_{I|J} - C_ {s,I|J}) (\u_{I,1}| \X_J=\x_{J,1}),\ldots,\sqrt{nh_n^{d-p}}(\hat C_{I|J} - C_{s,I|J}) (\u_{I,q}| \X_J=\x_{J,q}), \right. \nonumber }\\
    && \left. \sqrt{n}(\hat C_{s,I|J} -   C_ {s,I|J}) (\u_{I,q+1}),\ldots,\sqrt{n}(\hat C_{s,I|J} -  C_ {s,I|J}) (\u_{I,q+r}) \right\},
    \hspace{5cm}
\end{eqnarray*}
converges in law towards a Gaussian random vector.
As a consequence, $\sqrt{nh_n^{d-p}}\Tc^{0,m}_{KS,n}$ and $nh_n^{d-p}\Tc^{0,m}_{CvM,n}$ tends to a complex but not degenerate law under the $\Hc_0$.

\mds

\begin{rem}
    Other test statistics of $\Hc_0$ can be obtained by comparing directly the functions $\hat C_{I|J}(\cdot | \X_J=\x_J)$, for different values of $\x_J$.
    For instance, let us define
    \begin{eqnarray}
        \label{Tc0KS_bis}
        \lefteqn{\tilde\Tc^0_{KS,n}:=\sup_{\x_J,\x_J'\in \RR^{d-p}}\|\hat C_{I|J}(\cdot | \x_J) -\hat C_{I|J}(\cdot | \x'_J)\|_{\infty} \nonumber }\\
        &=&
        \sup_{\x_J,\x_J'\in \RR^{d-p}} \sup_{\u_I \in [0,1]^p}|\hat C_{I|J}(\u_I | \x_J) -\hat C_{I|J}(\u_I | \x'_J) |,
    \end{eqnarray}
    or
    \begin{equation}
        \label{TcOCvM_bis}
        \tilde\Tc^0_{CvM,n}:=\int \left( \hat C_{I|J}(\u_I | \x_J) -\hat C_{I|J}(\u_I | \x'_J)  \right)^2\, w(d\u_I,d\x_J,d\x'_J),
    \end{equation}
    for some function of bounded variation $w$. As above, modified versions of these statistics can be obtained considering fixed $\x_J$-grids.
    Since these statistics involve higher dimensional integrals/sums than previously, they will not be studied more in depth. 
\end{rem}


The $L^2$-type statistics $\Tc^0_{CvM,n}$ and $\tilde\Tc^0_{CvM,n}$ involve at least $d$ summations or integrals, which can become numerically expensive when the dimension of $\X$ is ``large''.
Nonetheless, we are free to set convenient weight functions. To reduce the computational cost, several versions of $\Tc^0_{CvM,n}$ are particularly well-suited,
by choosing conveniently the functions $w$. For instance, consider
\begin{equation*}
    \Tc^{(1)}_{CvM,n} := \int \left(
    \hat C_{I|J}(\u_I | \x_J) - \hat{C}_{s,I|J}(\u_I) \right)^2
    \, \hat C_I (d\u_I) \, \hat F_J(d\x_J),
    \label{CvM_Alexis}
\end{equation*}
where $\hat F_J$ and $\hat C_I$ denote the empirical cdf of $(\X_{i,J})$ and the empirical copula of $(\X_{i,I})$ respectively. Therefore, $\Tc^{(1)}_{CvM,n}$ simply becomes
\begin{equation}
    \Tc^{(1)}_{CvM,n} =
    \frac{1}{n^2} \sum_{j=1}^n \sum_{i=1}^n \left(
    \hat C_{I|J} (\hat U_{i,I} | \X_J=\X_{j,J}) -
    \hat C_{s,I|J} (\hat U_{i,I}) \right)^2, \label{T_CvM_1}
\end{equation}
where $\hat U_{i,I}=(\hat F_1 (X_{i,1}), \dots, \hat F_p (X_{i,p}))$, $i=1,\ldots,n$.
Similarly, we can choose
\begin{eqnarray*}
\lefteqn{    \tilde\Tc^{(1)}_{CvM,n} := \int \left(
    \hat C_{I|J}(\u_I | \x_J) - \hat C_{I|J}(\u_I | \x'_J)
    \right)^2 \, \hat C_I (d\u_I) \, \hat F_J(d\x_J) \, \hat F_J(d\x'_J) }\\
&=&
\frac{1}{n^3}\sum_{j=1}^n \sum_{j'=1}^n \sum_{i=1}^n \left( \hat C_{I|J}(\hat U_{i,I} | \X_J=\X_{j,J})
    - \hat C_{I|J}(\hat U_{i,I} |\X_J=\X_{j',J}) \right)^2.
\end{eqnarray*}

To deal with a single summations only, it is even possible to propose to set
\begin{align*}
    \Tc^{(2)}_{CvM,n} := \int &\left(
    \hat C_{I|J} (\hat F_{1|J} (x_1 | \x_J), \dots,
    \hat F_{p|J} (x_p|\x_J) |\x_J) \right. \\
    & \qquad - \left.
    \hat{C}_{s,I|J}( \hat F_{1|J} (x_1 | \x_J), \dots,
    \hat F_{p|J}(x_p|\x_J)) \right)^2 \, \hat F(d\x_I, d\x_J) ,
\end{align*}
where $\hat F$ denotes the empirical cdf of $\X$.
This means
\begin{align*}
    \Tc^{(2)}_{CvM,n} = \frac{1}{n} \sum_{i=1}^n
    \Big(& \hat C_{I|J} \big(
    \hat F_{1|J}(X_{i,1} | \X_{i,J}), \dots,
    \hat F_{p|J}(X_{i,p} | \X_{i,J}) | \X_J=\X_{i,J} \big)   \nonumber\\
    &- \hat{C}_{s,I|J} \big(
    \hat F_{1|J}(X_{i,1} | \X_{i,J}), \dots,
    \hat F_{p|J}(X_{i,p} | \X_{i,J}) \big) \Big)^2.
    \label{CvM_JDF}
\end{align*}

\mds

We have introduced some tests based on comparisons between empirical cdfs'. Obviously, the same idea could be applied to associated densities, as in Fermanian (2005) for instance, or even
to other functions of the underlying distributions.

\mds

Since the previous test statistics are complicated functionals of
some ``semi-smoothed'' empirical process, it is very challenging to evaluate their asymptotic laws under $\Hc_0$ analytically.
In every case, these limiting laws will not be distribution free, and their calculation would be very tedious.
Therefore, as usual with copulas, it is necessary to evaluate the limiting distributions of such tests
statistics by a convenient bootstrap procedure (parametric or nonparametric).
{\color{black} These bootstrap techniques will be presented in Section \ref{Boot_SA}}.

\subsection{Tests based on the independence property}
\label{IndepProp_SA}

Actually, testing $\Hc_0$ is equivalent to a test of the independence between the random vectors $\X_J$ and
$\Z_{I|J}:=(F_{1}(X_{1}|\X_J),\ldots,F_{p}(X_{p}|\X_J))$ strictly speaking, as proved in the following proposition.

\begin{prop}
    \label{prop_indep_H0}
    The vectors $\Z_{I|J}$ and $\X_J$ are independent iff
    $C_{I|J}(\u_I|\X_J=\x_J)$ does not depend on $\x_J$ for every
    vectors $\u_I$ and $\x_J$. In this case, the cdf of $\Z_{I|J}$ is $C_{s,I|J}$.
\end{prop}

\mds

{\it Proof:} For any vectors $\u_I \in [0,1]^p$ and any subset
$A_J\subset \RR^{d-p}$,
\begin{eqnarray*}
    \lefteqn{ \PP(\Z_{I|J} \leq \u_I , \X_J\in A_J) = \EE\left[ \1 (\X_J\in A_J ) \PP(\Z_{I|J} \leq \u_I |\X_J) \right] }\\
    &=& \int \1 (\x_J\in A_J ) \PP(\Z_{I|J} \leq \u_I | \X_J=\x_J) \, d\PP_{\X_J}(\x_J) \\
    &=& \int_{A_J} \PP(F_{k}(X_{k} | \X_J=\x_J)\leq u_{k},\forall k\in I | \X_J=\x_J) \, d\PP_{\X_J}(\x_J) \\
    &=& \int_{A_J} C_{I|J}( \u_I | \X_J=\x_J) \, d\PP_{\X_J}(\x_J).
\end{eqnarray*}
If $\Z_{I|J}$ and $\X_J$ are independent, then
$$\PP(\Z_{I|J} \leq \u_I) \PP(\X_J\in A_J) = \int \1 (\x_J\in A_J ) C_{I|J}( \u_I | \X_J=\x_J) \, d\PP_{\X_J}(\x_J),$$
for every $\u_I$ and $A_J$. This implies $\PP(\Z_{I|J} \leq \u_I) =
C_{I|J}( \u_I | \X_J=\x_J)$ for every $\u_I\in [0,1]^{p}$ and every
$\x_J$ in the support of $\X_J$. This means that
$C_{I|J}(\u_I|\X_J=\x_J)$ does not depend on $\x_J$, because
$\Z_{I|J}$ does not depend on any $\x_J$ by definition.

\medskip

Reciprocally, under $\Hc_0$, $C_{s,I|J}$ is the cdf of $\Z_{I|J}$.
Indeed,
\begin{eqnarray*}
    \lefteqn{\PP( \Z_{I|J} \leq \u_I ) =
    \PP \left(  F_{k}(X_{k}|\X_J) \leq u_{k},\forall k\in I    \right) }\\
    &=& \int \PP \left(  F_{k}(X_{k}|\X_J=\x_J) \leq u_{k},\forall k \in I |\, \X_J=\x_J \right)\, d\PP_{\X_J}(\x_J)   \\
    &=& \int C_{I|J}(\u_I | \X_J=\x_J)\, d\PP_{\X_J}(\x_J)= \int C_{s,I|J}(\u_I)\, d\PP_{\X_J}(\x_J) =C_{s,I|J}(\u_I).
\end{eqnarray*}
Moreover, due to Sklar's Theorem, we have
\begin{eqnarray*}
\lefteqn{    \PP(\Z_{I|J} \leq \u_I , \X_J\in A_J) = \int \1 (\x_J\in A_J ) C_{I|J}( \u_I | \X_J=\x_J) \, d\PP_{\X_J}(\x_J) }\\
    &= \int \1 (\x_J\in A_J ) C_{s,I|J}( \u_I ) \, d\PP_{\X_J}(\x_J) =\PP(\Z_{I|J} \leq \u_I ) \PP( \X_J\in A_J),
\end{eqnarray*}
implying the independence between $\Z_{I|J}$ and $\X_J$. $\Box$

\medskip

Then, testing $\Hc_0$ is formally equivalent to testing
$$\Hc^*_0:
\Z_{I|J}=(F_{1}(X_{1}|\X_J),\ldots,F_{p}(X_{p}|\X_J)) \; \text{and}\; \X_J\; \text{are independent}.$$

\mds

Since the conditional marginal cdfs' are not observable, keep in mind that we have to work with
pseudo-observations in practice, i.e. vectors of observations that are not
independent. In other words, our tests of independence should be based on pseudo-samples
\begin{equation}
    \left(  \hat F_{1|J}(X_{i,1}|\X_{i,J}), \dots,
    \hat F_{p|J}(X_{i,p}|\X_{i,J}) \right)_{i=1,\dots,n}
    :=(\hat \Z_{i,I|J})_{i=1,\ldots,n},
    \label{pseudoObs}
\end{equation}
for some consistent estimate $\hat F_{k|J}(\cdot|\X_J)$, $k\in I$ of the conditional cdfs',
for example as defined in Equation~(\ref{hatF_k_J}).
The chance of getting distribution-free asymptotic statistics will be very tiny, and we will have to rely on some bootstrap techniques again.
To summarize, we should be able to apply some usual tests of independence, but replacing iid observations with (dependent) pseudo-observations.

\mds

Most of the tests of $\Hc^*_0$ rely on the joint law of $(\Z_{I|J},\X_J)$, that may be evaluated empirically as
\begin{eqnarray*}
    \lefteqn{G_{I,J}(\x_I,\x_J):=\PP(\Z_{I|J}\leq \x_I,\X_J\leq \x_J)  }\\
    &\simeq &\hat G_{I,J}(\x):= n^{-1} \sum_{i=1}^n \1 (\hat \Z_{i,I|J} \leq \x_I,\X_{i,J} \leq \x_J).
\end{eqnarray*}

Now, let us propose some classical strategies to build independence tests.
\begin{itemize}
\item Chi-square-type tests of independence:
Let $B_1,\ldots,B_N$ (resp. $A_1,\ldots,A_{m}$) some disjoint subsets in $\RR^p$ (resp. $\RR^{d-p}$).
\begin{equation}
    \label{IcChi}
    \Ic_{\chi,n} = n\sum_{k=1}^{N} \sum_{l=1}^{m} \frac{\left( \hat G_{I,J}(B_k\times A_l) - \hat G_{I,J}(B_k\times \RR^{d-p}) \hat G_{I,J}(\RR^{p} \times A_l)  \right)^2}{\hat G_{I,J}(B_k\times \RR^{d-p}) \hat G_{I,J}(\RR^{p} \times A_l)} \cdot
\end{equation}


\item Distance between distributions:
\begin{equation}
    \label{IcKS}
    \Ic_{KS,n} = \sup_{\x \in \RR^d} | \hat G_{I,J}(\x) - \hat G_{I,J}(\x_I,\infty^{d-p}) \hat G_{I,J}(\infty^{p}, \x_J)|,\; \text{or}
\end{equation}
\begin{equation}
    \label{Ic2n}
    \Ic_{2,n} = \int \left( \hat G_{I,J}(\x) - \hat G_{I,J}(\x_I, \infty^{d-p}) \hat G_{I,J}(\infty^{p} , \x_J) \right)^2 \omega(\x)\, d\x,
\end{equation}
for some (possibly random) weight function $\omega$. Particularly, we can propose the single sum
\begin{eqnarray}
    \label{IcCvM}
    \lefteqn{ \Ic_{CvM,n} = \int \left( \hat G_{I,J}(\x) - \hat G_{I,J}(\x_I, \infty^{d-p}) \hat G_{I,J}(\infty^{p} , \x_J) \right)^2 \, \hat G_{I,J}(d\x)  \nonumber }\\
    &=& \frac{1}{n} \sum_{i=1}^n
    \left( \hat G_{I,J}(\hat \Z_{i,I|J},\X_{i,J} ) - \hat G_{I,J}(\hat \Z_{i,I|J}, \infty^{d-p}) \hat G_{I,J}(\infty^{p} ,\X_{i,J}) \right)^2 .
\end{eqnarray}

\item Tests of independence based on comparisons of copulas:
let $\breve C_{I,J}$ and $\hat C_{J}$ be the empirical copulas based on the pseudo-sample
$(\hat \Z_{i,I|J}, \X_{i,J})_{i=1,\ldots,n}$, and
$( \X_{i,J})_{i=1,\ldots,n}$ respectively. Set
\begin{equation*}
    \breve \Ic_{KS,n} = \sup_{\u \in [0,1]^d} | \breve C_{I,J}(\u) - \hat C_{s,I|J}^{(k)} (\u_I) \hat C_J (\u_J) | , k=1,\ldots,5,\, \text{or}
\end{equation*}
\begin{equation*}
    \breve \Ic_{2,n} = \int_{\u \in [0,1]^d}
    \left( \breve C_{I,J}(\u) - \hat C_{s,I|J}^{(k)}(\u_I)
    \hat C_J (\u_J) \right)^2 \omega(\u)\, d\u,
\end{equation*}
and in particular
\begin{equation*}
    \breve \Ic_{CvM,n} = \int_{\u \in [0,1]^d}
    \left( \breve C_{I,J}(\u) - \hat C_{s,I|J}^{(k)}(\u_I)
    \hat C_J (\u_J)\right)^2 \, \breve C_{I,J}(d\u).
\end{equation*}
The underlying ideas of the test statistics $\breve \Ic_{KS,n}$ and $\breve\Ic_{CvM,n}$ are similar to those that have been proposed by Deheuvels (1979,1981) in the case of unconditional copulas.
Nonetheless, in our case, we have to calculate pseudo-samples of the pseudo-observations
$(\hat \Z_{i,I|J})$ and $( \X_{i,J})$, instead of a usual pseudo-sample of $(\X_i)$.
\end{itemize}

\mds

Note that the latter techniques require the evaluation of some
conditional distributions, for instance by kernel smoothing.
Therefore, the level of numerical complexity of these test statistics of
$\Hc_0^*$ is comparable with those we have proposed before to test
$\Hc_0$ directly.

\subsection{Parametric tests of the simplifying assumption}
\label{ParApproach_SA}

In practice, modelers often assume a priori that the underlying
copulas belong to some specified parametric family $\Cc:=\{ C_\theta, \theta \in \Theta\subset \RR^m\}$.
Let us adapt our tests under this parametric assumption.
Apparently, we would like to test
$$ \check\Hc_0: C_{I|J}(\cdot| \X_J)= C_{\theta}(\cdot),\; \text{for some }
\theta\in \Theta \; \text{and almost every } \X_J.$$
Actually, $\check\Hc_0$ requires two different things: the fact that the conditional copula is a constant copula w.r.t. its conditioning events (test of $\Hc_0$) and, additionally, that the right copula belongs to $\Cc$ (classical composite Goodness-of-Fit test).
Under this point of view, we would have to adapt ``omnibus'' specification tests to manage conditional copulas and pseudo observations.
For instance, and among of alternatives, we could consider an amended version of Andrews (1997)'s specification test
$$ CK_n := \frac{1}{\sqrt{n}}\max_{j\leq n} |\sum_{i=1}^n \left[  \1(\hat\Z_{i,I|J}\leq \hat\Z_{j,I|J}) - C_{\hat\theta_0}(\hat \Z_{j,I|J})  \right]\1( \X_{i,J}\leq \X_{j,J})|,$$
recalling the notations in~(\ref{pseudoObs}).
For other ideas of the same type, see Zheng (2000) and the references therein.

\mds

The latter global approach is probably too demanding.
Here, we prefer to isolate the initial problem that was related to the simplifying assumption only.
Therefore, let us assume that, for every $\x_J$, there
exists a parameter $\theta(\x_J)$ such that $C_{I|J}(\cdot|
\x_J)=C_{\theta(\x_J)}(\cdot)$. To simplify, we assume the function $\theta(\cdot)$ is continuous.
Our problem is then reduced to testing the
constancy of $\theta$, i.e.
$$\Hc^c_0: \text{the function } \x_J \mapsto  \theta(\x_J) \text{ is
a constant, called} \;\theta_0.$$

\mds

For every $\x_J$, assume we estimate $\theta(\x_J)$ consistently. For instance, this can be done by modifying the standard semiparametric Canonical Maximum Likelihood methodology (Genest et al., 1995, Tsukahara, 2005): set
$$ \hat\theta(\x_J):= \arg\max_{\theta\in\Theta} \sum_{i=1}^n \log c_\theta \left(
\hat F_{1|J}(X_{i,1}|\X_J=\X_{i,J}), \dots,
\hat F_{p|J}(X_{i,p}|\X_J=\X_{i,J}) \right) \cdot K_n(\X_{i,J}, \x_J),$$
through usual kernel smoothing in $\RR^{d-p}$ {, \color{black} where
$c_\theta(\u) := \frac{\partial^p C_\theta(\u)}
{\partial u_1 \cdots \partial u_p}$
for $\theta\in\Theta$ and $\u\in [0,1]^p$}.
Alternatively, we could consider
$$\tilde\theta(\x_J) := \arg\max_{\theta\in\Theta} \sum_{i=1}^n \log c_\theta \left(
\hat F_{1|J}(X_{i,1}|\X_J=\x_J), \dots,
\hat F_{p|J}(X_{i,p}|\X_J=\x_J) \right) \cdot K_n(\X_{i,J}, \x_J),$$ instead of $\hat\theta(\x_J).$
See Abegaz et al. (2012) concerning the theoretical properties of $\tilde\theta(\x_J)$ and some choice of conditional cdfs'.
Those of $ \hat\theta(\x_J)$ remain to be stated precisely, to the best of our knowledge.
But there is no doubt both methodologies provide consistent estimators, even jointly, under some conditions of regularity.

\mds

Under $\Hc^c_0$, the natural ``unconditional'' copula parameter $\theta_0$ of the copula of the $\Z_{I|J}$ will be estimated by
\begin{equation}
    \hat\theta_0:= \arg\max_{\theta\in\Theta} \sum_{i=1}^n \log c_\theta \left(
    \hat F_{1|J}(X_{i,1}|\X_{i,J}), \dots,
    \hat F_{p|J}(X_{i,p}|\X_{i,J}) \right).
    \label{CML}
\end{equation}
Surprisingly, the theoretical properties of the latter estimator do not seem to have been established in the literature explicitly.
Nonetheless, the latter M-estimator is a particular case of those considered in Fermanian and Lopez (2015)
in the framework of single-index models when the link function is a known function (that does not depend on the index).
Therefore, by adapting their assumption in the current framework, we easily obtain that $\hat \theta_0$ is consistent and
asymptotically normal if $c_\theta$ is sufficiently regular, for convenient choices of bandwidths and kernels.

\mds

Now, there are some challengers to test $\Hc_0^c$:
\begin{itemize}
    \item Tests based on the comparison between $\hat\theta(\cdot)$ and $\hat\theta_0$:
        \begin{equation}
            \Tc_{\infty}^c := \sup_{\x_J\in \RR^{d-p}} \| \hat\theta(\x_J) - \hat\theta_0 \| ,\; \text{or}\;
            \Tc_{2}^c := \int \|\hat\theta(\x_J) - \hat\theta_0\|^2 \omega(\x_J)\, d\x_J,
            \label{Tc_infty_C}
        \end{equation}
    for some weight function $\omega$.

    \item Tests based on the comparison between $C_{\hat\theta(\cdot)}$ and $C_{\hat\theta_0}$:
    \begin{equation}
        \Tc_{dist}^c := \int dist\left(C_{\hat\theta(\x_J)},C_{\hat\theta_0} \right) \omega(\x_J)\,
        d\x_J,
        \label{Tcdfcomp}
    \end{equation}
    for some distance $dist(\cdot,\cdot)$ between cdfs'.

    \item Tests based on the comparison between copula densities (when they exist):
    \begin{equation}
        \Tc_{dens}^c := \int \left(c_{\hat\theta(\x_J)}(\u_I) -c_{\hat\theta_0}(\u_I)\right)^2 \omega(\u_I,\x_J)\,
        d\u_I\,d\x_J.
        \label{Tdensitycomp}
    \end{equation}
\end{itemize}

\begin{rem}
It might be difficult to compute some of these integrals numerically, because of unbounded supports.
One solution is to to make change of variables. For example,
\begin{equation*}
    \Tc_{2}^c = \int \| \hat \theta (F_J^{-} (\u_J))
    - \hat \theta_0 \|^2 \omega(F_J^{-} (\u_J)) \,
    \frac{d\u_J}{f_J(F_J^{-} (\u_J))} \cdot
\end{equation*}
Therefore, the choice $\omega = f_J$ allows us to simplify the latter statistics to
$\int \| \hat \theta (F_J^{-} (\u_J))
- \hat \theta_0 \|^2 d\u_J$, which is rather easy to evaluate. We used this trick in the numerical section below.
\end{rem}

\subsection{Bootstrap techniques for tests of $\Hc_0$}
\label{Boot_SA}

It is necessary to evaluate the limiting laws of the latter test statistics under the null.
As a matter of fact, we generally cannot exhibit explicit - and distribution-free a fortiori - expressions for these limiting laws.
The common technique is provided by bootstrap resampling schemes.

\mds

More precisely, let us consider a general statistics $\Tc$, built from the initial sample $\Sc:=(\X_1,\ldots,\X_n)$. The main idea of the bootstrap is to construct $N$ new samples $\Sc^*:=(\X_1^*,\ldots,\X_n^*)$ following a given resampling scheme given $\Sc$. Then, for each bootstrap sample $\Sc^*$, we will evaluate a bootstrapped test statistics $\Tc^*$, and the empirical law of all these $N$ statistics is used as an approximation of the limiting law of the initial statistic $\Tc$.

\subsubsection{Some resampling schemes}
\label{resampling_sch}

The first natural idea is to invoke Efron's usual ``nonparametric bootstrap'', where we draw independently with replacement $\X_i^*$ for $i=1,\ldots,n$ among the initial sample $\Sc=(\X_1,\ldots,\X_n)$.
This provides a bootstrap sample $\Sc^*:=(\X_1^*,\ldots,\X_n^*)$.

\mds

The nonparametric bootstrap is an ``omnibus'' procedure whose theoretical properties are well-known
but that may not be particularly adapted to the problem at hand.
Therefore, we will propose alternative sampling schemes that should be of interest, even if we do not state their validity on the theoretical basis.
Such a task is left for further researches.

\mds

An natural idea would be to use some properties of $\X$ under $\Hc_0$, in particular the characterization given in Proposition~\ref{prop_indep_H0}: under $\Hc_0$, we known that $\Z_{i,I|J}$ and $\X_{i,J}$ are independent.
This will be only relevant for the tests of Subsection~\ref{IndepProp_SA}, and for a few tests of Subsection~\ref{Bruteforce_SA},
where such statistics are based on the pseudo-sample $(\hat \Z_{i,I|J}, \X_{i,J})_{i=1,\ldots,n}$.
Therefore, we propose the following so-called ``pseudo-independent bootstrap'' scheme:

\mds

\noindent
Repeat, for $i=1$ to $n$,

\begin{enumerate}
    \item draw $\X^*_{i,J}$ among $(\X_{j,J})_{j=1,\ldots,n}$;
    \item draw $\hat \Z^*_{i,I|J}$ independently, among the observations $\hat \Z_{j,I|J}$, $j=1,\ldots,n$.
\end{enumerate}

\noindent This provides a bootstrap sample
$\Sc^*:=\big( (\hat \Z^*_{1,I|J} , \X^*_{1,J}), \ldots, (\hat \Z^*_{n,I|J} , \X^*_{n,J}) \big)$.

\mds

    Note that we could invoke the same idea, but with a usual nonparametric bootstrap perspective: draw with replacement a $n$-sample among the pseudo-observations
    $(\hat \Z_{i,I|J}, \X_{i,J})_{i=1,\ldots,n}$ for each bootstrap sample.
    This can be called a ``pseudo-nonparametric bootstrap'' scheme.

\mds

Moreover, note that we cannot draw independently $\X_{i,J}^*$ among $(\X_{j,J})_{j=1,\ldots,n}$,
and beside $\X_{i,I} ^*$ among $(\X_{j,I})_{j=1,\ldots,n}$ independently.
Indeed, $\Hc_0$ does not imply the independence between $\X_I$ and $\X_J$.
At the opposite, it makes sense to build a ``conditional bootstrap'' as follows:

\mds

\noindent
Repeat, for $i=1$ to $n$,

\begin{enumerate}
    \item draw $\X^*_{i,J}$ among $(\X_{j,J})_{j=1,\ldots,n}$;
    \item draw $\hat\X^*_{i,I}$ independently, along the estimated conditional law of $\X_I$ given $\X_J=\X^*_{i,J}$.
    This can be down by drawing a realization along the law $\hat F_{I|J} (\cdot |\X_J=\X^*_{i,J})$, for instance (see~(\ref{hatF_I_J})).
    This can be done easily because the latter law is purely discrete, with unequal weights that depend on $\X^*_{i,J}$ and $\Sc$.
\end{enumerate}

\noindent This provides a bootstrap sample
$\Sc^*:=\big( (\hat \X^*_{1,I} , \X^*_{1,J}), \ldots, (\hat \X^*_{n,I} , \X^*_{n,J}) \big)$.

\mds

\begin{rem}
Note that the latter way of resampling is not far from the usual nonparametric bootstrap. Indeed, when the bandwidths tend to zero, once $\x_J^*=\X_{i,J}$ is drawn, the procedure above will select the other components of $\X_i$ (or close values), i.e. the probability that $\x_I^* =\X_{i,I}$ is ``high''.
\end{rem}

\bigskip

In the parametric framework, we might also want to use an appropriate resampling scheme.
As a matter of fact, all the previous resampling schemes can be used, as in the nonparametric framework,
 but we would not take advantage of the parametric hypothesis, i.e. the fact that all conditional copulas belong to a known family.
We have also to keep in mind that even if the conditional copula has a parametric form, the global model is not fully parametric, because we have not provided a parametric model neither for the conditional marginal cdfs $F_{k|J}$, $k=1,\ldots,p$, nor for the cdf of $\X_J$.

\mds

Therefore, we can invoke the null hypothesis $\Hc^c_0$ and approximate the real copula $C_{\theta_0}$ of $\Z_{I|J}$ by $C_{\hat \theta_0}$.
This leads us to     define the following ``parametric independent bootstrap'':

\mds

\noindent
Repeat, for $i=1$ to $n$,

\begin{enumerate}
    \item draw $\X_{i,J}^*$ among $(\X_{j,J})_{j=1,\dots,n}$;
    \item sample $\Z_{i,I|J,\hat \theta_0}^*$ from the copula with parameter $\hat \theta_0$ independently.
\end{enumerate}

\noindent This provides a bootstrap sample
$\Sc^*:=\big( ( \Z^*_{1,I|J, \hat\theta_0} , \X^*_{1,J}), \ldots, ( \Z^*_{n,I|J, \hat\theta_0} , \X^*_{n,J}) \big)$.

\begin{rem}
    At first sight, this might seem like a strange mixing of parametric and nonparametric bootstrap.
    If $|J|=1$, we can nonetheless do a ``full parametric bootstrap'', by observing that all estimators of our previous test statistics do not depend on $\X_J$, but on realizations of $\hat F_{J}(\X_J)$ (see Equations (\ref{hatF_I_J}) and (\ref{hatF_k_J})).
    Since the law of latter variable is close to a uniform distribution, it is tempting to sample $V_{i,J}^* \sim {\mathcal U}_{[0,1]}$ at the first stage, $i=1,\ldots,n$, and then to replace $\hat F_{J}(\X_{i,J})$ with $V_{i,J}^*$ to get an alternative bootstrap sample.
\end{rem}

Without using $\Hc^c_0$, we could define the ``parametric conditional bootstrap'' as:

\mds

\noindent
Repeat, for $i=1$ to $n$,
\begin{itemize}
    \item draw $\X_{i,J}^*$ among $(\X_{j,J})_{j=1,\dots,n}$;
    \item sample $\Z_{i,I|J,\theta^*_i}^*$ from the copula with parameter $\hat \theta(\X^*_{i,J})$.
\end{itemize}
\noindent This provides a bootstrap sample
$\Sc^*:=\big( ( \Z^*_{1,I|J,\theta^*_i} , \X^*_{1,J}), \ldots, ( \Z^*_{n,I|J,\theta^*_i} , \X^*_{n,J}) \big)$.

\mds

Note that, in several resampling schemes, we should be able to keep the same $\X_J$ as in the original sample, and simulate only $\Z^*_{i,I|J}$ in step 2,
as in Andrews(1997), pages 10-11.
Such an idea has been proposed by Omelka et al. (2013), in a slightly different framework and univariate conditioning variables.
They proved that such a bootstrap scheme ``works'', after a fine-tuning of different smoothing parameters: see their Theorem 1.

\subsubsection{Bootstrapped test statistics}

The problem is now to evaluate the law of a given test statistic, say $\Tc$, under $\Hc_0$ by the some bootstrap techniques.
We recall the main technique in the case of the classical nonparametric bootstrap.
We conjecture that the idea is still theoretically sound under the other resampling schemes that have been proposed in Subsection~\ref{resampling_sch}.

\mds

The principle for the nonparametric bootstrap is based on the weak convergence of the underlying empirical process. Formally, if $\Sc:=\{\X_1,\ldots,\X_n\}$ in an iid sample in $\RR^d$, $\X\sim F$ and if $F_n$ denotes its empirical distribution, it is well-known that $\sqrt{n}\left(F_n - F\right)$ tends weakly in $\ell^{\infty}$ towards a $d$-dimensional Brownian bridge $\BB_F$.
And the nonparametric bootstrap works in the sense that
$\sqrt{n}\left(F_n^* - F_n\right)$ converges weakly towards a process $\BB_F'$, an independent version of $\BB_F$, given the initial sample $\Sc$.

\mds

Due to the Delta Method, for every Hadamard-differentiable functional $\chi$ from $\ell^\infty(\RR^d)$ to $\RR$, there exists a random variable
$H_\chi$ s.t. $ \sqrt{n}\left( \chi(F_n) - \chi(F)\right) \Rightarrow H_\chi.$
Assume a test statistics $\Tc_n$ of $\Hc_0$ can be written as a sufficiently regular functional of the underlying empirical process as
\begin{equation*}
    \label{eq:decomp_Tc_1}
    \Tc_n := \psi \left( \sqrt{n}\left( \chi_s(F_n) - \chi(F_n)\right) \right),
\end{equation*}
where $\chi_s(F)=\chi(F)$ under the null assumption.
Then, under $\Hc_0$, we can rewrite this expression as
\begin{equation}
    \label{eq:decomp_Tc_2}
    \Tc_n := \psi \left( \sqrt{n}\left( \chi_s(F_n) - \chi_s(F) + \chi(F) - \chi(F_n)\right) \right).
\end{equation}
Given any bootstrap sample $\Sc^*$ and the associated empirical distribution $F_n^*$,
the usual bootstrap equivalent of $\Tc_n$ is
\begin{equation*}
    \label{eq:decomp_Tc_st}
    \Tc_n^* := \psi \left( \sqrt{n}\left( \chi_s(F_n^*) - \chi_s(F_n) + \chi(F_n) - \chi(F_n^*)\right) \right),
\end{equation*}
from Equation (\ref{eq:decomp_Tc_2}).
See van der Vaart and Wellner (1996), Section 3.9, for details and mathematically sound statements.

\mds

Applying these ideas, we can guess the bootstrapped statistics corresponding to the
tests statistics of $\Hc_0$, at least when the usual nonparametric bootstrap is invoked.

Let us illustrate the idea with $\Tc_{KM,n}^0$.
Note that $\hat C_{I|J}( \cdot |\X_J=\cdot) = \chi_{KM}(F_n)(\cdot)$ and $\hat C_{s,I|J} = \chi_{s,KM}(F_n)$ for some
smoothed functional $\chi_{KM}$ and $\chi_{s,KM}$.
Under $\Hc_0$, $\chi_{KM}=\chi_{s,KM}$ and $\Tc^0_{KS,n}:=\|\chi_{KM}(F_n) - \chi_{KM}(F) -\chi_{s,KM}(F_n) + \chi_{s,KM}(F) \|_{\infty}$.
Therefore, its bootstrapped version is
\begin{eqnarray*}
    \label{Tc0KS_B}
    \lefteqn{ \Tc^{0,*}_{KS,n}
    :=\|\chi_{KM}(F_n^*) - \chi_{KM}(F_n) -\chi_{s,KM}(F_n^*) + \chi_{s,KM}(F_n) \|_{\infty}   }\\
    &=&\|\hat C_{I|J}^* - \hat C_{I|J} -\hat C^*_{s,I|J} + \hat C_{s,I|J}\|_{\infty}.  \hspace{5cm}\nonumber
\end{eqnarray*}
Obviously, the functions $\hat C_{I|J}^*$ and $\hat C_{s,I|J}^*$ have been calculated as $\hat C_{I|J}$ and $\hat C_{s,I|J}$ respectively, but replacing $\Sc$ by $\Sc^*$.
Similarly, the bootstrapped versions of some Cramer von-Mises-type test statistics are
\begin{equation*}
    \label{TcOCvM_B}
    \Tc^{0,*}_{CvM,n}:=\int \left( \hat C^*_{I|J}(\u_I | \x_J) - \hat C_{I|J}(\u_I | \x_J) -\hat C^*_{s,I|J}(\u_I) + \hat C_{s,I|J}(\u_I) \right)^2\, w(d\u_I,d\x_J).
\end{equation*}

When playing with the weight functions $w$, it is possible to keep the same weights for the bootstrapped versions, or to replace them with some functionals of $F_n^*$.
For instance, asymptotically, it is equivalent to consider
\begin{equation*}
    \Tc^{(1),*}_{CvM,n}:=\int \left( \hat C_{I|J}^*(\u_I | \x_J) - \hat C_{I|J}(\u_I | \x_J) -\hat{C}^*_{s,I|J}(\u_I) + \hat{C}_{s,I|J}(\u_I) \right)^2\, \hat{C}_{n}(d\u_I) \, \hat F_J(d\x_J),\; \text{or}
\end{equation*}
\begin{equation*}
    \Tc^{(1),*}_{CvM,n}:=\int \left( \hat C_{I|J}^*(\u_I | \x_J) - \hat C_{I|J}(\u_I | \x_J) -\hat{C}^*_{s,I|J}(\u_I) + \hat{C}_{s,I|J}(\u_I) \right)^2\, \hat{C}^*_{n}(d\u_I) \, \hat F^*_J(d\x_J).
\end{equation*}
Similarly, the limiting law of
\begin{eqnarray*}
    \lefteqn{ \Tc^{(2),*}_{CvM,n}:=\int \left( \hat C^*_{I|J}(\hat F^*_{n,1}(x_1 | \x_J),\ldots, \hat F^*_{n,p}(x_p|\x_J) |\x_J) \right. }\\
    &-& \left. \hat C_{I|J}(\hat F^*_{n,1}(x_1 | \x_J),\ldots, \hat F^*_{n,p}(x_p|\x_J) |\x_J) -
    \hat{C}_{s,I|J}^*( \hat F^*_{n,1}(x_1 | \x_J),\ldots, \hat F^*_{n,p}(x_p|\x_J))
    \right. \\
    &+& \left.
    \hat{C}_{s,I|J}( \hat F_{n,1}^*(x_1 | \x_J),\ldots, \hat F^*_{n,p}(x_p|\x_J))\right)^2 \, H_n(d\x_I, d\x_J) ,
\end{eqnarray*}
given $F_n$ is unchanged replacing $H_n$ by $H_n^*$.

\mds

The same ideas apply concerning the tests of Subsection~\ref{IndepProp_SA}, but they require some modifications.
Let $H$ be some cdf on $\RR^d$. Denote by $H_I$ and $H_J$ the associated cdf on the first $p$ and $d-p$ components respectively.
Denote by $\hat H$, $\hat H_I$ and $\hat H_J$ their empirical counterparts.
Under $\Hc_0$, and for any measurable subsets $B_I$ and $A_J$, $ H(B_I\times A_J) = H(B_I) H(A_J)$. Our tests will be based on the difference
\begin{eqnarray*}
    \lefteqn{\hat H(B_I\times A_J) - \hat H_I(B_I) \hat H_J(A_J)= (\hat H-H)(B_I\times A_J) }\\
    & -&   (\hat H_I - H_I)(B_I) \hat H_J(A_J) -
    (\hat H_J - H_J)(A_J) H_I(B_I) .
\end{eqnarray*}
Therefore, a bootstrapped approximation of the latter quantity will be
\begin{equation*}
    (\hat H^*-\hat H)(B_I\times A_J) -   (\hat H^*_I - \hat H_I)(B_I) \hat H^*_J(A_J) -
    (\hat H^*_J - \hat H_J)(A_J) \hat H_I(B_I).
\end{equation*}
To be specific, the bootstrapped versions of our tests are specified as below.
\begin{itemize}
    \item Chi-square-type test of independence:
    \begin{eqnarray*}
        \lefteqn{ \Ic^*_{\chi,n} := n\sum_{k=1}^{N} \sum_{l=1}^{m} \frac{1}{\hat G^*_{I,J}(B_k\times \RR^{d-p}) \hat G^*_{I,J}(\RR^{p} \times A_l)}
        \left( (\hat G^*_{I,J} - \hat G_{I,J})(B_k\times A_l)  \right. }\\
        &-& \left. \hat G^*_{I,J} (B_k\times \RR^{d-p}) \hat G^*_{I,J}(\RR^{p} \times A_l)+
        \hat G_{I,J}(B_k\times \RR^{d-p})  \hat G_{I,J}(\RR^{p} \times A_l)
         \right)^2.
    \end{eqnarray*}

    \item Distance between distributions:
    \begin{align*}
        \Ic^*_{KS,n} = \sup_{\x \in \RR^d} | (\hat G^*_{I,J}- \hat G_{I,J})(\x) - \hat G^*_{I,J}(\x_I,\infty^{d-p}) \hat G^*_{I,J}(\infty^{p}, \x_J)
        + \hat G_{I,J}(\x_I,\infty^{d-p})  \hat G_{I,J}(\infty^{p}, \x_J)|
    \end{align*}
    $$    \Ic^*_{2,n} = \int \bigg( (\hat G^*_{I,J}- \hat G_{I,J})(\x)
        - \hat G^*_{I,J}(\x_I, \infty^{d-p}) \hat G^*_{I,J}(\infty^{p} , \x_J)
        + \hat G_{I,J}(\x_I, \infty^{d-p}) \hat G_{I,J}(\infty^{p} , \x_J) \bigg)^2 \omega(\x)\, d\x,$$
    and $\Ic^*_{CvM,n}$ is obtained replacing $\omega(\x)\,d\x$ by $\hat G^*_{I,J}(d\x)$ (or even $\hat G_{I,J}(d\x)$).

    \item A test of independence based on the independence copula: Let $\breve C^*_{I,J}$, $\breve C^*_{I|J}$ and $\hat C^*_{J}$ be the empirical copulas based on a bootstrapped version of the pseudo-sample
    $(\hat \Z_{i,I|J}, \X_{i,J})_{i=1,\ldots,n}$, $(\hat \Z_{i,I|J})_{i=1,\ldots,n}$ and $( \X_{i,J})_{i=1,\ldots,n}$ respectively. This version can be obtained by nonparametric bootstrap, as usual, providing new vectors $\hat \Z^*_{i,I|J}$ at every draw.
    The associated bootstrapped statistics are
    \begin{align*}
        &\breve \Ic^*_{KS,n} = \sup_{\u \in [0,1]^d} | (\breve C^*_{I,J}- \breve C_{I,J})(\u) - \breve C^*_{I|J}(\u_I)\hat C^*_{J}(\u_J)+ \breve C_{I|J}(\u_I)\hat C_{J}(\u_J) | , \\
        &\breve \Ic^*_{2,n}  = \int_{\u \in [0,1]^d} \left( (\breve C^*_{I,J}- \breve C_{I,J}) (\u) - \breve C^*_{I|J}(\u_I)\hat C^*_{J}(\u_J)+  \breve C_{I|J}(\u_I)\hat C_{J}(\u_J)\right)^2 \omega(\u)\, d\u, \\
        &\breve \Ic^*_{CvM,n} = \int_{\u \in [0,1]^d} \left( (\breve C^*_{I,J}- \breve C_{I,J})(\u) - \breve C^*_{I|J}(\u_I)\hat C^*_{J}(\u_J) +
        \breve C_{I|J}(\u_I)\hat C_{J}(\u_J) \right)^2 \, \breve C^*_{I,J}(d\u).
    \end{align*}
\end{itemize}

\mds

In the case of the parametric statistics, the situation is pretty much the same, as long as we invoke the nonparametric bootstrap. For instance, the bootstrapped versions of some previous test statistics are
\begin{equation*}
    \left(\Tc_{2}^c \right)^* := \int \|
    \hat\theta^*(\x_J) - \hat\theta (\x_J) - \hat\theta^*_0 + \hat\theta_0\|^2 \omega(\x_J)\,d\x_J, \; \text{or}
    \label{T2compBootNP}
\end{equation*}
\begin{equation*}
    \left(\Tc_{dens}^c \right)^* := \int \left(
    c_{\hat\theta^*(\x_J)}(\u_I) - c_{\hat\theta (\x_J)}(\u_I)
    - c_{\hat\theta_0^*}(\u_I) + c_{\hat\theta_0}(\u_I)
    \right)^2 \omega(\u_I,\x_J)
    d\u_I\,d\x_J.
    \label{TdensitycompBootNP}
\end{equation*}
in the case of the nonparametric bootstrap.
We conjecture that the previous techniques can be applied with the other resampling schemes that have been proposed in Subsection~\ref{resampling_sch}.
Nonetheless, a complete theoretical study of all these alternative schemes and the statement of the validity of their associated bootstrapped statistics is beyond the scope of this paper.

\begin{rem}
    \label{tricky_schemes}
    For the ``parametric independent'' bootstrap scheme, we have observed
    that the test powers are a lot better by considering
    \begin{equation*}
        \left(\Tc_{2}^c \right)^{**} := \int \|
        \hat\theta^*(\x_J) - \hat\theta^*_0 \|^2 \omega(\x_J)\,d\x_J,\; \text{or}
        \label{T2compBootP}
    \end{equation*}
    \begin{equation*}
        \left(\Tc_{dens}^c \right)^{**} := \int \left(
        c_{\hat\theta^*(\x_J)}(\u_I)
        - c_{\hat\theta_0^*}(\u_I)
        \right)^2 \omega(\u_I,\x_J)
        d\u_I\,d\x_J
        \label{TdensitycompBootP},
    \end{equation*}
    instead. The relevance of such statistics may be theoretically justified in the slightly different context of ``box-type'' tests in the next Section (see Theorem~\ref{thm:ParIndepBoot}).
    Since our present case is close to the situation of ``many small boxes'', it is not surprising that we observe similar features. Note that, contrary to the nonparametric bootstrap
    or the ``parametric conditional'' bootstrap, the ``parametric independent'' bootstrap scheme uses $\Hc_0$.
    {\color{black} More generally, and following the same idea, we found that using the statistic $\Tc^{**} := \psi \left( \sqrt{n}\left( \chi_s(F_n^*) - \chi(F_n^*)\right) \right)$ for the pseudo-independent bootstrap yields much better performance than $\Tc^{*}$. In our simulations, we will therefore use $\Tc^{**}$ as the bootstrap test statistic (see Figures \ref{fig:I_chi_tau_max} and \ref{fig:I2n_tau_max}).
    }
\end{rem}

\mds

{\color{black}
\begin{rem}
For testing $\Hc_0$ at a node of a vine model, the realizations of the corresponding explanatory variables $\X_{i,J}$ are not observed in general. In practice, they have to be replaced with pseudo-observations in our previous test statistics. Their calculation involves the bivariate conditional copulas that are associated with the previous nodes in a recursive way. The theoretical analysis of the associated bootstrap schemes is challenging and falls beyond the scope of the current work. 
\end{rem}
}

\mds

\section{Tests with ``boxes''}
\label{Boxes}

\subsection{The link with the simplifying assumption}
\label{Link_SA}

As we have seen in Remark \ref{ex_simple}, we do not have $C_{s,I|J} = C_{I}$ in general, when $C_I(\u_I)=C_{I|J}(\u_I| \X_J \in \RR^{d-p})$ for every $\u_I$.
This is the hint there are some subtle relations between conditional copulas when the conditioning event is pointwise or when it is a measurable subset.
Actually, to test $\Hc_0$ in Section~\ref{Tests_SA}, we have relied on kernel estimates and smoothing parameters, at least to evaluate conditional marginal distributions empirically.
To avoid the curse of dimension (when $d-p$ is ``large'' i.e. larger than three in practice),
it is tempting to replace the pointwise conditioning events $\X_J=\x_J$ with $\X_J\in A_J$ for some borelian subsets $A_J\subset \RR^{d-p}$, $\PP (\X_J \in A_J) > 0$. 
As a shorthand notation, we shall write ${\mathcal A}_J$ the set of all such $A_J$. We call them ``boxes'' because choosing $d-p$-dimensional rectangles
(i.e. intersections of half-spaces separated by orthogonal hyperplans) is natural,
but our definitions are still valid for arbitrary borelian subsets in $\RR^{d-p}$.
Technically speaking, we will assume that the functions $\x_J \mapsto \1(\x_J\in A_J)$ are Donsker, to apply uniform CLTs' without any hurdle.
Actually, working with $\X_J$-``boxes'' instead of pointwise will simplify a lot the picture.
{\color{black} Indeed, the evaluation of conditional cdfs' given $\X_J\in A_J$ does not require kernel smoothing, bandwidth choices, or 
other techniques of curve estimation that deteriorate the optimal rates of convergence.}

\mds

Note that, by definition of the conditional copula of $\X_I$ given $(\X_J\in A_J)$, we have
\begin{eqnarray*}
\lefteqn{
 \PP( \X_I \leq \x_I | \X_J \in A_J )} \\
 &=&C_{I|J}^{A_{J}}\left(\PP( X_{1} \leq
x_{1} | \X_J \in A_J ),\ldots,\PP( X_{p} \leq x_{p} | \X_J \in A_J )|
\X_J\in A_J\right),
\end{eqnarray*}
for every point $\x_I\in \RR^p$ and every subset
$A_J$ in ${\mathcal A}_J$.
So, it is tempting to replace $\Hc_0$ by
$$\widetilde \Hc_0: C_{I|J}^{A_{J}}(\u_I | \X_J \in A_J) \, \text{does not depend on} \, A_J \in {\mathcal A}_J, \text{for any } \u_I.$$

{\color{black} For any $\x_J$, consider a sequence of boxes $(A^{(n)}_J(\x_J) )$ s.t. $\cap_n A^{(n)}_J(\x_J) = \{\x_J\}$. If the law of $\X$ is sufficiently regular, then 
$\lim_n C_{I|J}^{A_{J}^{(n)}}(\u_I | \X_J \in A_J^{(n)}) =  C_{I|J}(\u_I | \X_J =\x_J)$ for any $\u_I$. Therefore, $\tilde\Hc_0$ implies $\Hc_0$. This is stated formally in the next proposition.

\mds

\begin{prop}
Assume that the function $h:\RR^{d} \rightarrow [0,1]$, 
defined by $h(\y):=\PP(\X_I \leq \y_I | \X_J=\y_J)$ is continuous everywhere.  
Then, for every $\x_J\in \RR^{d-p}$ and any sequence of boxes $(A^{(n)}_J(\x_J) )$ such that $\cap_n A^{(n)}_J(\x_J) = \{\x_J\}$, we have 
$$\lim_n C_{I|J}^{A_{J}^{(n)}(\x_J)}(\u_I | \X_J \in A_J^{(n)}(\x_J)) =  C_{I|J}(\u_I | \X_J =\x_J),$$ 
for every $\u_I\in [0,1]^p$.
\end{prop}

{\it Proof:} 
Consider a particular $\u_I\in [0,1]^p$. 
If one component of $\u_I$ is zero, the result is obviously satisfied.
If one component of $\u_I$ is one, this component does not play any role. 
Therefore, we can restrict ourselves on $\u_I\in (0,1)^p$. 
By continuity, there exists $\x_I\in \RR^p$ s.t. $u_i=F_i(x_i|\x_J)$ for every $i=1,\ldots,p$. 
Let the sequences $(x_i^{(n)})$ such that
$u_i=F_i(x^{(n)}_i|\X_J\in A_J^{(n)})$ for every $n$ and every $i=1,\ldots,p$. 
First, let us show that $x^{(n)}_i \to x_i$ when $n$ tends to the infinity. 
Indeed, by the definition of conditional probabilities (Shiryayev 1984, p.220), we have
    $$ u_i=\PP(X_i \leq x_i^{(n)} | \X_J\in A^{(n)}_J ) = \frac{1}{\PP(\X_J\in A_J^{(n)})} \int_{\{\y_J\in A_J^{(n)}\}} \PP(X_i \leq x_i^{(n)} | \X_J=\y_J)\, d\PP_{\X_J}(\y_J)  ,$$
and
\begin{eqnarray*}
\lefteqn{ u_i = \PP(X_i\leq x_i |\X_J=\x_J)= 
\frac{1}{\PP(\X_J\in A_J^{(n)})} \int_{\{\y_J\in A_J^{(n)}\}} \PP(X_i \leq x_i^{(n)} | \X_J=\x_J)\, d\PP_{\X_J}(\y_J)   }\\
&+& \PP( X_i\leq x_i |\x_J ) - \PP(X_i\leq x_i^{(n)} | \x_J). \hspace{5cm}
\end{eqnarray*}

By substracting the two latter identities, we deduce
\begin{eqnarray}
\lefteqn{    
\frac{1}{\PP(\X_J\in A_J^{(n)})} \int_{\{\y_J\in A_J^{(n)}\}} 
    \left[ \PP(X_i \leq x_i^{(n)} | \X_J=\y_J) - \PP(X_i \leq x_i^{(n)} | \X_J=\x_J)   \right]\, d\PP_{\X_J}(\y_J) \nonumber}\\
    &=& \PP(X_i\leq x_i |\x_J) - \PP(X_i\leq x_i^{(n)} |\x_J).  \label{Dinis} \hspace{5cm}
\end{eqnarray}

But, by assumption, $F_i(t |\y_J)$ tends towards $F_i(t|\x_J)$ when 
$\y_J $ tends to $\x_J$, for any $t$ (pointwise convergence). A straightforward application of Dini's Theorem shows that the latter convergence is uniform on $\RR$: $\|  F_i(\cdot | \y_J) - F_i(\cdot | \x_J) \|_{\infty}$ tends to zero when $\y_J \to \x_J$. 
From~(\ref{Dinis}), we deduce that $\PP(X_i\leq x_i^{(n)} |\x_J)\to \PP(X_i\leq x_i |\x_J)  $. By the continuity of $F_i(\cdot |\x_J)$, we get $x_i^{(n)}\to x_i$, for any $i=1,\ldots,p$.

\mds

Second, let us come back to conditional copulas: setting 
$\x_I^{(n)}:=(x_1^{(n)},\ldots,x_p^{(n)})$, we have
\begin{eqnarray*}
    \lefteqn{ C^{A_J^{(n)}}_{I|J} (\u_I | A_J^{(n)} ) - C_{I|J}(\u_I | \x_J) }\\
    &=&  
    C^{A_J^{(n)}}_{I|J} (F_1(x_1^{(n)}|A_J^{(n)}),\ldots, F_p(x_p^{(n)}|A_J^{(n)}) | A_J^{(n)} ) 
    - C_{I|J}(F_1(x_1|\x_J),\ldots, F_p(x_p|\x_J) | \x_J)\\
    &=& F_{I|J}(\x_I^{(n)} | A_J^{(n)}) - F_{I|J}(\x_I | \x_J) \\
    &=& \frac{1}{\PP(\X_J\in A_J^{(n)})} \int_{\{\y_J\in A_J^{(n)}\}} \left[\PP(\X_I\leq \x_I^{(n)} |\y_J)- \PP(\X_I\leq \x_I |\x_J) \right]  \, d\PP_{\X_J}(\y_J).
\end{eqnarray*}
Since $\x_I^{(n)}$ tends to $\x_I$ when $n\to \infty$ and invoking the continuity of $h$ at $(\x_I,\x_J)$, we get that 
$C^{A_J^{(n)}}_{I|J} (\u_I | A_J^{(n)} ) \to C_{I|J}(\u_I | \x_J)$ when $n\to \infty$. $\Box$
}

\mds

Unfortunately, the opposite is false. Counter-intuitively, $\tilde\Hc_0$ does not lead to a consistent test of the simplifying assumption. Indeed, under $\Hc_0$, we can see that $ C_{I|J}^{A_{J}}(\u_I | \X_J \in A_J)$ {\bf depends on $A_J$} in general, even if $ C_{I|J}(\u_I | \X_J =\x_J)$ {\bf does not depend on $\x_J$}!

\mds

This is due to the nonlinear transform between conditional (univariate and multivariate) distributions and conditional copulas. In other words, for a usual $d$-dimensional cdf $H$, we have
\begin{equation}
    \label{condsubsetdistr}
    H(\x_I | \X_J \in A_J)=\frac{1}{\PP(A_J)} \int_{A_J} H(\x_I | \X_J =\x_J) \, d\PP_{\X_J}(\x_J),
\end{equation}
for every measurable subset $A_J \in {\mathcal A}_J$ and $\x_I\in \RR^p$. At the opposite and in general, for conditional copulas,
\begin{equation}
    \label{condsubsetcopula}
    C_{I|J}^{A_{J}}(\u_I | \X_J \in A_J)\neq
    \frac{1}{\PP(A_J)} \int_{A_J} C_{I|J}(\u_I | \X_J =\x_J) \, d\PP_{\X_J}(\x_J),
\end{equation}
for $\u_I\in [0,1]^p$. And even if we assume $\Hc_0$, we have in general,
\begin{equation}
    \label{condsubsetcopula_0}
    C_{I|J}^{A_{J}}(\u_I | \X_J \in A_J)\neq
    \frac{1}{\PP(A_J)} \int_{A_J} C_{s,I|J}(\u_I) \, d\PP_{\X_J}(\x_J) = C_{s,I|J}(\u_I).
\end{equation}

\mds

As a particular case, taking $A_J = \RR^{d-p}$, this means again that $C_{I}(\u_I) \neq C_{s,I|J}(\u_I)$.

\mds

Let us check this rather surprising feature with the example of Remark~\ref{ex_simple} for another subset $A_J$. Recall that $\Hc_0$ is true and that $C_{s,1,2|3}(u,v)=uv$ for every
$u,v \in [0,1]$. Consider the subset $(X_3 \leq a)$, for any real number $a$. The probability of this event is $\Phi(a)$.
Now, let us verify that
$$ uv \neq H(F_{1|3}^- (u | X_3\leq a),F_{2|3}^- (v | X_3\leq a) | X_3\leq a),$$
for some $u,v$ in $(0,1)$. Clearly, for every real number $x_k$, we have
$$\PP(X_k \leq x_k | X_3\leq a)=
\frac{1}{\Phi(a)}\int_{-\infty}^a  \Phi(x_k - z)\phi(z) \, dz, k=1,2,\; \text{and}$$
$$ \PP(X_1\leq x_1,X_2\leq x_2 | X_3\leq a)= \frac{1}{\Phi(a)}\int_{-\infty}^a \Phi(x_1-z)\Phi(x_2-z)\phi(z)\,dz.$$

In particular, $\PP(X_k \leq 0 | X_3\leq a)= (1+\Phi(-a))/2$. Therefore, set $u^*=v^*=(1+\Phi(-a))/2$ and
we get
\begin{eqnarray*}
\lefteqn{ H(F_{1|3}^- (u^* | X_3\leq a),F_{2|3}^- (v^* | X_3\leq a) | X_3\leq a)=  H(0,0 |X_3\leq a)   }\\
&=&\frac{1}{3}\left(1 + \Phi(-a)+ \Phi^2(-a) \right)\neq u^* v^*.
\end{eqnarray*}
In this example, $ C_{s,1,2|3}(\cdot)\neq C_{1,2|3}^{]-\infty,a]}(\cdot | X_3 \leq a),$
for every $a$, even if $\Hc_0$ is satisfied.

\bigskip

Nonetheless, getting back to the general case, we can easily provide an equivalent of Equation (\ref{condsubsetdistr}) for general conditional copulas, i.e. without assuming $\Hc_0$.

\begin{prop}
    \label{condsubsetcopgeneral}
    For all $\u_I \in [0,1]^p$ and all $ A_J \in {\mathcal A}_J$,
    \begin{align*}
          \nonumber \\
        &C_{I|J}^{A_{J}} (\u_I | \X_J \in A_J) = \frac{1}{\PP(A_J)}
        \int_{A_J} \psi(\u_I , \x_J, A_J) d\PP_{\X_J}(\x_J),\; \text{with}
    \end{align*}

    \begin{align*}
        \psi(&\u_I , \x_J, A_J) \nonumber \\
        &= C_{I|J} \Bigg(
        F_{1|J} \Big( F_{1|J}^- (u_1 | \X_J \in A_J) \big| \X_J=\x_J \Big)
        , \dots,
        F_{p|J} \Big( F_{p|J}^- (u_p | \X_J \in A_J) \big| \X_J=\x_J \Big)
        \Bigg| \X_J = \x_J \Bigg).
    \end{align*}
\end{prop}

{\color{black}
{\it Proof:} From (\ref{condsubsetdistr}), we get :
\begin{align*}
    H&(\x_I | \X_J \in A_J) \\
    &=\frac{1}{\PP(A_J)} \int_{A_J} H(\x_I | \X_J =\x_J) \, d\PP_{\X_J}(\x_J) \\
    &=\frac{1}{\PP(A_J)} \int_{A_J} C_{I|J} \Big(
    F_{1|J}(x_{1} | \X_J = \x_J ) , \ldots , F_{p|J}(x_{p} | \X_J = \x_J )
    \, \big| \, \X_J = \x_J
    \Big) \, d\PP_{\X_J}(\x_J).
\end{align*}
We can conclude by using the following definition of the conditional copula
\begin{align*}
    C_{I|J}^{A_{J}} (\u_I | \X_J \in A_J) =
    H(F_{1|J}^- (u_1 | \X_J \in A_J) , \dots, F_{p|J}^- (u_p | \X_J \in A_J) | \X_J \in A_J).\;\; \Box
\end{align*}
}

\medskip

Now, we understand why (\ref{condsubsetcopula}) (and (\ref{condsubsetcopula_0}) under $\Hc_0$) are not identities:
the conditional copulas, given the subset $A_J$, still depend on the conditional margins of $\X_I$ given $\X_J$ pointwise in general.

\mds

Note that, if $X_i$ is independent of $\X_J$ for every $i=1,\ldots,p$, then, for any such $i$,
\begin{equation*}
    F_{i|J} \Big( F_{i|J}^- (u_i | \X_J \in A_J) \big| \X_J=\x_J \Big)
    = F_{i} \Big( F_{i}^- (u_i) \Big)= u_i,
    \label{indep_XiXJ}
\end{equation*}
{\color{black} and we can revisit the identity of Proposition~\ref{condsubsetcopgeneral}: under $\Hc_0$, we have 
\begin{eqnarray*}
\lefteqn{C_{I|J}^{A_{J}} (\u_I | \X_J \in A_J)  
= \frac{1}{\PP(A_J)} \int_{A_J} C_{I|J}(\u_I | \X_J=\x_J) \, d\PP_{\X_J}(\x_J)  }\\
&=& \frac{1}{\PP(A_J)} \int_{A_J} C_{s,I|J}(\u_I ) \, d\PP_{\X_J}(\x_J) = C_{s,I|J}(\u_I).
\end{eqnarray*}
This means $\Hc_0$ and $\tilde \Hc_0$ are equivalent.} We consider such circumstances as very peculiar and do not have to be confused with a test of $\Hc_0$. Therefore, we advise to lead a preliminary test of independence between $\X_I$ and $\X_J$ (or at least between $X_i$ and $\X_J$ for any $i=1,\ldots,p$) before trying to test $\Hc_0$ itself.

\mds

Now, let us revisit the characterisation of $\Hc_0$ in terms of the independence property, as in Subsection~\ref{IndepProp_SA}.
The latter analysis is confirmed by the equivalent of Proposition \ref{prop_indep_H0} in the case of conditioning subsets $A_J$. Now, the
relevant random vector would be
\begin{equation*}
    \Z_{I|A_J} := \left(F_{1|J}(X_{1}|\X_J \in A_J) , \dots,F_{p|J}(X_{p}|\X_J \in A_J) \right),
\label{def:Z_I_AJ}
\end{equation*}
that has straightforward empirical counterparts.
Then, it is tempting to test
\begin{equation*}
    \widetilde \Hc_0^{*}: \Z_{I|A_J} \text{ and } (X_J \in A_J)\;\text{are independent for every borelian subset}\; A_J\subset \RR^{d-p}.
    \label{atester**}
\end{equation*}
Nonetheless, it can be proved easily that this is not a test of $\Hc_0$, unfortunately.
\begin{prop}
    $\Z_{I|A_J}$ and $ (X_J \in A_J)$ are independent for every measurable subset $A_J\subset \RR^{d-p}$ iff $\X_I$ and $\X_J$ are independent.
\end{prop}

\mds

{\it Proof:}
For any measurable subset $A_J$ and any $\u_I\in [0,1]^p$, under $\widetilde \Hc_0^{*}$, we have
$$ \PP\left( \Z_{I|A_J} \leq \u_I, \X_J\in A_J\right)
= \PP\left( \Z_{I|A_J} \leq \u_I  \right) \PP (  \X_J\in A_J  ) .$$
Consider $\x_I\in \RR^p$. Due to the continuity of the conditional cdfs', there exists $u_k$ s.t.
 $F_k(x_k |\X_J\in A_J)=u_k$, $k=1,\ldots,p$. Then, using the invertibility of $x\mapsto F_k(x | \X_J\in A_J)$, we get
$ \PP\left( \Z_{I|A_J} \leq \u_I ,  \X_J\in A_J \right)
= \PP\left( \X_I \leq \x_I,  \X_J\in A_J \right).$
This implies that $\widetilde \Hc_0^{*}$ is equivalent to the following property: for every $\x_I\in \RR^p$ and $A_J$,
$$ \PP\left( X_I\leq \x_I, \X_J\in A_J\right) =
\PP\left( \X_I\leq \x_I \right)\PP\left( \X_J\in A_J\right).\; \Box$$

\mds

{\color{black} The previous result shows that a test of $\tilde\Hc_0^*$ is a test of independence between $\X_I$ and $\X_J$. When the latter assumption is satisfied, $\tilde \Hc_0$ and then $\Hc_0$ are true too, but the opposite is false.}

\mds

Previously, we have exhibited a simple trivariate model where $\Hc_0$ is satisfied when $\X_I$ and $\X_J$ are not independent.
Then, we see that it is not reasonable to test whether the mapping $A_J \mapsto C_{I|J}^{A_{J}} (\cdot | \X_J \in A_J)$ is constant
over ${\mathcal A}_J$, the set of {\bf all} $A_J$ such that $\PP_{\X_J}(A_J) > 0$, with the idea of testing $\Hc_0$.

\mds

Nonetheless, one can weaken the latter assumption,
and restrict oneself to a {\bf finite} family $\bar {{\mathcal A}}_J$ of subsets with positive probabilities. For such a family, we could test the assumption
$$\bar \Hc_0:
A_J \mapsto C_{I|J}^{A_{J}} (\, \cdot \, | \X_J \in A_J)
\, \text{is constant over} \; \bar{{\mathcal A}}_J .$$

\mds

To fix the ideas and w.l.o.g., we will consider a given family of disjoint subsets $\bar{\Ac}_J = \{A_{1,J},\ldots,A_{m,J}\}$ in $\RR^{d-p}$ hereafter.
Note the following consequence of Proposition~\ref{condsubsetcopgeneral}.
\begin{prop}
    Assume that, for all $A_J \in \bar {{\mathcal A}}_J$ and for all $i \in I$,
    \begin{equation}
        F_{i|J} ( x | \X_J=\x_J ) = F_{i|J} (x | \X_J \in A_J), \;\; \forall \x_J \in A_J, x \in \RR.
        \label{cnd_technik}
    \end{equation}
    Then, $\Hc_0$ implies $\bar\Hc_0$.
\end{prop}

\mds

Obviously, if the family $\bar {{\mathcal A}}_J$ is too big, then~(\ref{cnd_technik}) will be too demanding: $\bar \Hc_0$ will be close to a test of independence between $\X_I$ and $\X_J$, and no longer a test of $\Hc_0$. Moreover, the chosen subsets in the family $\bar{\Ac}_J$ do not need to be disjoint, even if this would be a natural choice.
As a special case, if $\RR^{d-p} \in \bar{{\mathcal A}}_J$, the previous condition is equivalent to the independence between $X_i$ and $\X_J$ for every $ i \in I$.

\mds

{\color{black} Note that~(\ref{cnd_technik}) does not imply that the vector of explanatory variables $\X_J$ should be discretized. Indeed, the full model requires the specification of the underlying conditional copula too, independently of the conditional margins and arbitrarily. For instance, we can choose a Gaussian conditional copula whose parameter is a continuous function of $\X_J$, even if~(\ref{cnd_technik}) is fulfilled. And the law of $\X_I$ given $\X_J$ will depend on the current value of $X_J$.}

\mds

A test of $\bar\Hc_0$ may be relevant in a lot of situations, beside technical arguments as the absence of smoothing.
First, the case of discrete (or discretized) explanatory variables $\X_J$ is frequent.
When $\X_J$ is discrete and takes a value among $\{\x_{1,J},\ldots,\x_{m,J}\}$, set $A_{k,J}=\{\x_{k,J}\}$, $k=1,\ldots,m$.
Then, there is identity between testing $\Hc_0$ and $\bar\Hc_0$, with $\bar\Ac_J=\{ A_{1,J},\ldots,A_{m,J}\}$.
Second, the level of precision and sharpness of a copula model is often lower than the models for (conditional) margins.
To illustrate this idea, a lot of complex and subtle models to explain the dynamics of asset volatilities are available when the dynamics of cross-assets dependencies are often a lot more basic and without clear-cut empirical findings.
Therefore, it makes sense to simplify conditional copula models compared to conditional marginal models. This can be done by considering only a few possible conditional copulas, associated to some events $(\X_J\in A_{k,J})$, $k=1,\ldots,m$.
For example, Jondeau and Rockinger (2006) (the first paper that introduced conditional dependence structures, beside Patton (2006a)) proposed a Gaussian copula parameter that take a finite of values randomly, based on the realizations of some past asset returns.
Third, similar situation occur with most Markov-switching copula models, where a finite set of copulas is managed. In such models,
the (unobservable, in general) underlying state of the economy determines the index of the box: see Cholette et al. (2009), Wang et al. (2013), St\"ober and Czado (2014), Fink et al. (2016), among others.

\mds

Therefore, testing $\bar\Hc_0$ is of interest per se.
Even if this is not equivalent to $\Hc_0$ (i.e. the simplifying assumption) formally, the underlying intuitions are close.
And, particularly when the components of the conditioning variable $\X_J$ are numerous, it can make sense to restrict the information set of the underlying conditional copula to a fixed number of conveniently chosen subsets $ A_J$. And the constancy of the underlying copula when $\X_J$ belongs to such subsets is valuable in a lot of practical situations. Therefore, in the next subsections, we study some specific tests of $\bar\Hc_0$ itself.

\subsection{Non-parametric tests with ``boxes''}
\label{NPApproach_m_SA}

To specify such tests, we need first to estimate the conditional marginal cdfs', for instance by
\begin{equation*}
    \hat F_{k|J}(x|\X_J \in A_J)
    := \frac{ \sum_{i=1}^n
    \1 (X_{i,k} \leq x, \X_{i,J} \in A_J) }
    {\sum_{i=1}^n \1 (\X_{i,J} \in A_J)},
\end{equation*}
for every real $x$ and $k=1,\ldots,p$. Similarly the joint law of $\X_I$ given $(\X_J \in A_J)$ may be estimated by
\begin{equation*}
    \hat F_{I|J}(\x_I|\X_J \in A_J)
    := \frac{ \sum_{i=1}^n
    \1 (X_{i,I} \leq \x_I, \X_{i,J} \in A_J) }
    {\sum_{i=1}^n \1 (\X_{i,J} \in A_J)}\cdot
\end{equation*}
The conditional copula given $(\X_J\in A_J)$ will be estimated by
$$ \hat C_{I|J}^{A_J}(\u_I | \X_J \in A_J)
:= \hat F_{I|J}(
\hat F^-_{1|J}(u_1|\X_J \in A_J), \dots,
\hat F^-_{p|J}(u_p|\X_J \in A_J) |\X_J \in A_J).$$
Therefore, it is easy to imagine tests of $\bar\Hc_0$, for instance
\begin{equation}
    \bar\Tc_{KS,n} := \sup_{\u_I \in [0,1]^d} \sup_{k,l=1,\ldots,m}
    | \hat C_{I|J}^{A_{k,J}}(\u_I | \X_J\in A_{k,J})
    - \hat C_{I|J}^{A_{l,J}}(\u_I | \X_J\in A_{l,J})|,
    \label{bar_Tc_KS}
\end{equation}
\begin{equation}
    \bar\Tc_{CvM,n} := \sum_{k,l=1}^m \int \left(
    \hat C_{I|J}^{A_{k,J}}(\u_I | \X_J\in A_{k,J}) -
    \hat C_{I|J}^{A_{l,J}}(\u_I | \X_J\in A_{l,J}) \right)^2 w(d\u_I) ,
    \label{bar_Tc_CvM}
\end{equation}
for some nonnegative weight functions $w$, or even
\begin{equation}
    \bar\Tc_{dist,n} := \sum_{k,l=1}^m dist\left(
    \hat C_{I|J}^{A_{k,J}}(\cdot | \X_J\in A_{k,J}) ,
    \hat C_{I|J}^{A_{l,J}}(\cdot | \X_J\in A_{l,J})\right) ,
    \label{bar_Tc_dist}
\end{equation}
where $dist(\cdot,\cdot)$ denotes a distance between cdfs' on $[0,1]^p$.
More generally, define the matrix
    \begin{equation*}
        \widehat M(\bar \Ac_J) := \bigg[ \1(k \neq l) \;
        dist \Big(
        \hat C_{I|J}^{A_{k,J}}(\cdot | \X_J\in A_{k,J}) ,
        \hat C_{I|J}^{A_{l,J}}(\cdot | \X_J\in A_{l,J})
        \Big) \bigg]_{1 \leq k,l \leq m},
    \end{equation*}
    and any statistic of the form $||\widehat M (\bar \Ac_J)||$ can be used as a test statistics of $\bar\Hc_0$, when $||\cdot||$ is a norm on the set of $m\times m$-matrices.
Obviously, it is easy to introduce similar statistics based on copula densities instead of cdfs'.

\subsection{Parametric test statistics with ``boxes''}
\label{ParApproach_m_SA}

When we work with subsets $A_J\in \RR^{d-p}$
instead of pointwise conditioning events $(\X_J=\x_J)$, we can adapt all the
previous parametric test statistics of Subsection~\ref{ParApproach_SA}. Nonetheless, the framework will be slightly modified.

\mds

Let us assume that,
for every $A_J \in \bar \Ac_J $, $C_{I|J}^{A_{J}} (\cdot |\X_J \in A_J)$ belongs to the same parametric copula family $ \Cc=\{C_\theta,\theta\in \Theta\}$.
In other words,
$C_{I|J}^{A_{J}} (\cdot | \X_J \in A_J) = C_{\theta(A_J)}(\cdot)$ for every $A_J\in \bar \Ac_J $.
Therefore, we could test the constancy of the mapping $A_J\mapsto \theta(A_J)$, i.e. to test
$$\bar\Hc^c_0: \text{the function } k \in \{ 1, \dots, m \} \mapsto  \theta(A_{k,J}) \text{ is
a constant called} \; \theta_0^b.$$
%
Clearly, for every $A_J \in \bar{\Ac}_J $, we can
estimate $\theta(A_J)$ by
\begin{equation*}
    \hat\theta(A_J):= \arg\max_{\theta\in\Theta} \sum_{i=1}^n \log c_\theta \left(
    \hat F_{1|J}(X_{i,1}|\X_{i,J} \in A_J),\ldots,\hat F_{p|J}(X_{i,p}|\X_{i,J}
    \in A_J )\right)\1( \X_{i,J}\in A_J).
    \label{def_hattheta0}
\end{equation*}
It can be proved that the estimate $\hat\theta(A_J)$ is consistent and asymptotically normal, by revisiting the proof of Theorem 1 in Tsukahara (2005). Here, the single difference w.r.t. the latter paper is induced by the random sample size, modifying the limiting distributions. The proof is left to the reader.

\mds

Under the zero assumption $\bar \Hc_0^c$, the parameter of the copula of
$(F_{1}(X_{1}|\X_{J} \in A_{k,J}),\ldots, F_{p}(X_{p}|\X_{J} \in A_{k,J} ))$ given $(\X_J\in A_{k,J})$ is the same for any $k=1,\ldots,m$.
It will be denoted by $\theta_0^b$, and we can still estimate it by the semi-parametric procedure
\begin{equation*}
    \hat\theta_0^b:= \arg\max_{\theta\in\Theta} \sum_{k=1}^m\sum_{i=1}^n \log c_\theta \left(
    \hat F_{1|J}(X_{i,1}|\X_{i,J} \in A_{k,J}),\ldots, \hat F_{p|J}(X_{i,p}|\X_{i,J}
    \in A_{k,J} )\right)\1( \X_{i,J}\in A_{k,J}).
    \label{def_hatthetaAkJ}
\end{equation*}
Obviously, under some conditions of regularity and under $\bar\Hc^c_0$, it can be proved that
$\hat\theta_0^b$ is consistent and asymptotically normal, by adapting the results of Tsukahara (2005).

\mds

For convenience, let us define the ``box index'' function $k(\x_J) := \sum_{k=1}^m k \1  \{ \x_{J}\in A_{k,J} \} ,$ for any $\x_J\in \RR^{d-p}$.
In other words, $k$ is the index of the box $A_{k,J}$ that contains $\x_{J}$. It equals zero, when no box in $\bar\Ac_J$ contains $\x_J$.
Let us introduce the r.v. $Y_i := k(\X_{i,J}),$
that stores only all the needed information concerning the conditioning with respect to the variables $\X_{i,J}$.
We can then define the empirical pseudo-observations as
\begin{eqnarray*}
    \Z_{i,I|Y} & := & \sum_{k=1}^m \left(
    F_{1|J}(X_{i,1}|\X_{J}\in A_{k,J}), \dots,
    F_{p|J}(X_{i,p}|\X_{J}\in A_{k,J}) \right)
    \1 \{ \X_{i,J}\in A_{k,J} \}  \\
    & = &  \left( F_{1|J}(X_{i,1}|\X_{J}\in A_{k(\X_{i,J}),J}), \dots,
    F_{p|J}(X_{i,p}|\X_{J}\in A_{k(\X_{i,J}),J}) \right) \\
    & = &  \left( F_{1|Y}(X_{i,1}|Y_i), \dots,
    F_{p|Y}(X_{i,p}|Y_i)\right) ,
\end{eqnarray*}
for any $i=1,\ldots,n$.
Since we do not observe the conditional marginal cdfs', we define the observed pseudo-observations that we calculate in practice: for $i=1, \dots, n,$
\begin{equation*}
    \hat \Z_{i,I|Y} := \left(
    \hat F_{1|J} (X_{i,1} | \X_{J} \in A_{Y_i,J}) ,\dots,
    \hat F_{p|J} (X_{i,p} | \X_{J} \in A_{Y_i,J}) \right).
\end{equation*}
Note that we can then rewrite the previous estimators as
\begin{equation*}
    \hat\theta(A_{k,J}) = \arg\max_{\theta\in\Theta} \sum_{i=1}^n \log c_\theta \left(
    \hat \Z_{i,I|Y} \right) \1(Y_i = k),\;
\text{and}\;
    \hat\theta_0^b = \arg\max_{\theta\in\Theta} \sum_{i=1}^n \log c_\theta \left( \hat \Z_{i,I|Y} \right).
\end{equation*}

\mds

Now, let us revisit some of the previously proposed test statistics in the case of ``boxes''.
\begin{itemize}
    \item Tests based on the comparison between $\hat\theta(\cdot)$ and $\hat\theta_0$:
    \begin{equation}
        \bar\Tc_{\infty}^{c} := \sqrt{n} \max_{k=1,\ldots,m} \| \hat\theta(A_{k,J}) - \hat\theta_0 \|,\;
        \bar\Tc_{2}^{c} := n \sum_{k=1}^m \|\hat\theta(A_{k,J}) - \hat\theta_0\|^2 \omega_k,
        \label{bar_Tc_KS_c}
    \end{equation}
    for some weights $\omega_k$.

    \item Tests based on the comparison between $C_{\hat\theta(\cdot)}$ and $C_{\hat\theta_0}$:
    \begin{equation}
        \bar\Tc_{dist}^{c} := \sum_{k=1}^m dist( C_{\hat\theta(A_{k})} , C_{\hat\theta_0} ) \omega_k,
        \label{bar_Tc_dist_c}
    \end{equation}
\end{itemize}
and others.

\subsection{Bootstrap techniques for tests with boxes}
\label{Boot_m_SA}

In the same way as in the previous section, we will need bootstrap schemes to evaluate the limiting laws of the test statistics of $\bar\Hc_0$ or $\bar\Hc^c_0$ under the null.
All the nonparametric resampling schemes of Subsection \ref{resampling_sch} (in particular Efron's usual bootstrap) can be used in this framework, replacing the conditional pseudo-observations $\hat\Z_{i,I|J}$
by $\hat\Z_{i,I|Y}$, $i=1,\ldots,n$.
The parametric resampling schemes of Subsection \ref{resampling_sch} can also be applied to the framework of ``boxes'', replacing $\hat \theta_0$ by $\hat \theta_0^b$ and $\hat \theta(\x_J)$ by $\hat \theta(A_J)$.
In the parametric case, the bootstrapped estimates are denoted by $\hat\theta^*_{0}$ and $\hat\theta^*(A_{J})$. They are
the equivalents of $\hat\theta_{0}^b$ and $\hat\theta_{n}(A_{J})$, replacing $(\hat \Z_{i,I|J},Y_i)$ by $(\Z^*_{i},Y^*_i)$.

\mds

The bootstrapped statistics will also be changed accordingly.
Writing them explicitly is a rather straightforward exercise and we do not provide the details, contrary to Subsection~\ref{Boot_SA}.
For example, the bootstrapped statistics corresponding to~(\ref{bar_Tc_KS_c}) is
\begin{equation*}
    \big(\bar\Tc_{2}^{c} \big)^*:= n \sum_{k=1}^m \|
    \hat\theta^*(A_{k,J}) - \hat\theta(A_{k,J}) - \hat\theta_0^* +  \hat\theta_0^b\|^2 \omega_k,
    \label{T2comp1BootNP}
\end{equation*}
where $\hat\theta_0^*$ is the result of the program $\arg\max_\theta\sum_{i=1}^n \log c_\theta \left( \hat \Z^*_{i,I|Y} \right)$, in the case of Efron's
nonparametric bootstrap.

\mds

As we noticed in Remark~\ref{tricky_schemes}, some changes are required when dealing with the ``parametric independent'' bootstrap.
Indeed, under the alternative,
we observe $\hat\theta^*(A_{k,J}) - \hat\theta_0^* \approx 0$,
because we have precisely generated a bootstrap sample under $\bar\Hc^c_0$.
As a consequence, the law of $\big(\bar\Tc_{2}^{c} \big)^*$ would be close to the law of $\bar\Tc_{2}^{c}$ but under the alternative, providing very small powers.
Therefore, convenient bootstrapped test statistics of $\bar\Hc_0$ under the ``parametric independent'' scheme will be of the type
\begin{equation*}
    \big(\bar\Tc_{2}^{c} \big)^{**}:= n \sum_{k=1}^m \|
    \hat\theta^*(A_{k,J}) - \hat\theta_0^* \|^2 \omega_k.
    \label{T2comp1BootNP2}
\end{equation*}
Such a result is justified theoretically by the following theorem.

\begin{thm}
    \label{thm:ParIndepBoot}
    Assume that $\bar \Hc^c_0$ is satisfied, and that we apply the parametric independent bootstrap. Set
    \begin{equation*}
        \Theta{n,0} := \sqrt{n}
        \big( \hat \theta_0 - \theta_0 \big),
        \Theta_{n,k} := \sqrt{n}
        \big( \hat \theta (A_{k,J}) - \theta_0 \big),k=1, \dots, m,
    \end{equation*}
    \begin{equation*}
        \Theta^*_{n,0} := \sqrt{n}
        \big( \hat \theta^*_{0} - \theta_0 \big),\;\text{and}\;
        \Theta^*_{n,k} := \sqrt{n}
        \big( \hat \theta^*(A_{k,J}) - \theta_0 \big),k=1, \dots, m.
    \end{equation*}
    Then there exists two independent and identically distributed random vectors $\big( \Theta_0, \dots, \Theta_m \big)$
    and $\big( \Theta^{\perp}_0, \dots, \Theta^{\perp}_m \big)$, and a real number $a_0$ such that
    \begin{equation*}
        \Big( \Theta_{n,0}, \dots, \Theta_{n,m}, \Theta^*_{n,0}, \dots, \Theta^*_{n,m} \Big)
        \Longrightarrow
        \Big( \Theta_0, \dots, \Theta_m,
        \Theta^{\perp}_0 + a_0 \Theta_0, \dots,
        \Theta^{\perp}_m + a_0 \Theta_0  \Big).
\end{equation*}

\end{thm}

The proof of this theorem has been postponed in Appendix~\ref{Section_Proofs}.

\mds

As a consequence of the latter result, applying the parametric independent bootstrap procedures for some test statistics based on comparisons between
$\hat \theta_0$ and the $\hat\theta (A_{k,J})$ is valid.
For instance, $\bar\Tc_{2}^{c}$ and $\big( \bar\Tc_{2}^{c} \big)^{**}$ will converge jointly in distribution to a pair of independent and identically distributed variables.
Indeed, we have
\begin{eqnarray*}
    \lefteqn{ \left( \bar\Tc_{2}^{c} \, , \,
    \big( \bar\Tc_{2}^{c} \big)^{**} \right)
    = \left( n \sum_{k=1}^m \|
    \hat \theta_{n,0}^b - \hat \theta_{n} (A_{k,J})
    \|^2 \omega_k\, , \,
    n \sum_{k=1}^m \|
    \hat \theta^*_{n,0} - \hat \theta_{n}^* (A_{k,J})
    \|^2 \omega_k \right)   }\\
    &=& \left( n \sum_{k=1}^m \|
    \hat \theta_{n,0}^b - \theta_0 + \theta_0 - \hat \theta_{n} (A_{k,J})
    \|^2 \omega_k\, , \,
    n \sum_{k=1}^m \|
    \hat \theta^*_{n,0} - \theta_0 + \theta_0 - \hat \theta_{n}^* (A_{k,J})
    \|^2 \omega_k\right)\\
    &\Longrightarrow &
     \left( \sum_{k=1}^m \| \Theta_0 - \Theta_k
    \|^2 \omega_k\, , \,
    \sum_{k=1}^m \|
    \Theta^{\perp}_0 + a_0 \Theta_0 -
    \Theta^{\perp}_k - a_0 \Theta_0 \|^2 \omega_k \right).
\end{eqnarray*}
The same reasoning applies with $\bar\Tc_{\infty}^{c}$ and $\bar\Tc_{dist}^{c}$, for sufficiently regular copula families.

\begin{rem}
    We have to stress that the first-level bootstrap, i.e. resampling among the conditioning variables $\X_{i,J}$, $i=1,\ldots,n$ is surely necessary to obtain the latter result.
    Indeed, it can be seen that the key proposition~\ref{prop:bootstrapped_Donsker_th} is no longer true otherwise, because the limiting covariance functions of the two corresponding processes $\GG_n$ and $\GG_n^*$ will not be the same: see remark~\ref{cov_boot_v2} below.
\end{rem}

\section{Numerical applications}
\label{NumericalApplications}

Now, we would like to evaluate the empirical performances of some of the previous tests by simulation.
Such an exercise has been led by Genest et al. (2009) or Berg (2009) extensively in the case of goodness-of-fit test for unconditional copulas.
Our goal is not to replicate such experiments in the case of conditional copulas and for tests of the simplifying assumption.
Indeed, we have proposed dozens of test statistics and numerous bootstrap schemes.
Moreover, testing the simplifying assumption through $\Hc_0$ or some ``box-type'' problems through $\bar\Hc_0$ doubles the scale of the task.
Finally, in the former case, we depend on smoothing parameters that induce additional degrees of freedom for the fine tuning of the experiments 
(note that Genest et al. (2009) and Berg (2009)
have renounced to consider tests that require
additional smoothing parameters, as the pivotal test statistics proposed in Fermanian (2005)). In our opinion, an exhaustive simulation experiment should be the topic of (at least)
one additional paper. Here, we will restrict ourselves to some partial numerical elements.
They should convince readers that the methods and techniques we have discussed previously provide fairly good results and can be implemented in practice safely.

\mds

Hereafter, we consider bivariate conditional copulas and a single conditioning variable, i.e. $p=2$ and $d=3$.
The sample sizes will be $n=500$, except if it is differently specified.
Concerning the bootstrap, we will resample $N=200$ times to calculate approximated p-values.
Each experiment has been repeated $500$ times to calculate the percentages of rejection. {\color{black} The computations have been made on a standard laptop, and, for the non-parametric bootstrap, they took an average time of $14.1$ seconds for $\Ic_{\chi,n}$ ; $26.9$s for $\Tc^{0,m}_{CvM,n}$, $103$s for $\Ic_{2,n}$, $265$s for $\Tc_2^c$  and $0.922$s for $\bar \Tc_2^c$.}

\mds

In terms of model specification, the margins of $\X = (X_1,X_2,X_3)$ will depend on $X_3$ as
$$  X_1 \sim \Nc(X_3,1),\; X_2 \sim \Nc(X_3,1) \; \text{and} \; X_3 \sim \Nc(0,1).$$
We have studied the following conditional copula families: given $X_3=x$,
\begin{itemize}
\item the Gaussian copula model, with a correlation parameter $\theta(x)$,
\item the Student copula model, with $4$ degrees of freedom and a correlation parameter $\theta(x)$,
\item the Clayton copula model, with a parameter $\theta(x)$,
\item the Gumbel copula model, with a parameter $\theta(x)$,
\item the Frank copula model, with a parameter $\theta(x)$.
\end{itemize}
In every case, we calibrate $\theta(x)$ such that the conditional Kendall's tau $\tau(x)$ satisfies $\tau(x) = \Phi(x)\tau_{\max}$, for some constant $\tau_{\max}\in (0,1)$. 
By default, $\tau_{\max}$ is equal to one. 
In this case, the random Kendall's tau are uniformly distributed on $[0,1]$.

\mds

{\it Test of $\Hc_0$:} we calculate the percentage of rejections of $\Hc_0$, when the sample is drawn under the true law (level analysis) or when it is drawn under the same parametric copula family, but with varying parameters (power analysis). For example, when the true law is a Gaussian copula with a constant parameter $\rho$ corresponding to $\tau = 1/2$, we draw samples under the alternative through a bivariate Gaussian copula whose random parameters are given by $\tau(X_3)= \Phi(X_3)$.
The chosen test statistics are $\Tc_{CvM}^0$, $\tilde \Tc_{CvM}^0$ (nonparametrics test of $\Hc_0$), $\Ic_{\chi,n}$ and $\Ic_{2,n}$  (nonparametric tests of $\Hc_0$ based on the independence property) and $\Tc_2^c$ (a parametric test of $\Hc_0^c$).
To compute these statistics, we use the estimator of the simplified copula defined in Equation (\ref{estimator_meanCI_J}).

\mds

{\it Test of $\bar\Hc_0$:} in the case of the test with boxes, the data-generating process will be 
$$  X_1 \sim \Nc( \gamma(X_3),1),\; X_2 \sim \Nc(\gamma(X_3),1) \; \text{and} \; X_3 \sim \Nc(0,1),$$
where $\gamma(x) = \Phi^{-1} \left(  \lfloor m\Phi(X_3)  \rfloor /m\right)$, so that the boxes are all of equal probability. As $m \to \infty$, we recover the continuous model for which $\gamma(x) = x$.

\mds

In the same way, we calibrate the parameter $\theta(x)$ of the conditional copulas such that the conditional Kendall's tau satisfies
$\tau(X_3) = \lfloor m\Phi(X_3)  \rfloor /m$.

\mds

{\color{black} The choice of ``the best'' boxes $A_{1,J}, \dots, A_{m,J}$ is not an easy task. 
This problem happens frequently in statistics (think of Pearson's chi-square test of independence, for instance), and there is no universal answer.   
Nonetheless, in some applications, intuition can be fuelled by the context.
For example, in finance, it makes sense to test whether past positive returns induce different conditional dependencies between current returns than past negative returns. 
And, as a general ``by default'' rule, we can divide the space of $\X_J$ into several boxes of equal (empirical) probabilities. This trick is particularly relevant when the conditioning variable is univariate.
Therefore, in our example, we have chosen $m=5$ boxes of equal empirical probability for $X_3$, with equal weights.
}

\mds

We have only evaluated $\bar \Tc_{2}^c$ for testing $\bar\Hc_0^c$.
In the following tables, for the parametric tests,
\begin{itemize}
    \item ``bootNP'' means the usual nonparametric bootstrap ;
    \item ``bootPI'' means the parametric independent bootstrap (where $\Z_{I|J}$ is drawn under $C_{\hat \theta_0}$ and $\X_J$ under the usual 
    nonparametric bootstrap);
    \item ``bootPC'' means the parametric conditional bootstrap (nonparametric bootstrap for $\X_J$, and $\X_I$ is 
    sampled from the estimated conditional copula $C_{\hat \theta(\X^*_J)}$);
    \item ``bootPseudoInd'' means the pseudo-independent bootstrap (nonparametric bootstrap for $\X_J$, and 
    draw $\hat \Z^*_{I|J}$ independently, among the pseudo-observations $\hat \Z_{j,I|J}$);
    \item ``bootCond'' means the conditional bootstrap (nonparametric bootstrap for $\X_J$, and $\X_I$ is sampled from the estimated conditional law of $\X_I$ given $\X^*_J$).
\end{itemize}

\mds

\begin{figure}
    \centering
    \includegraphics[width = 16cm]{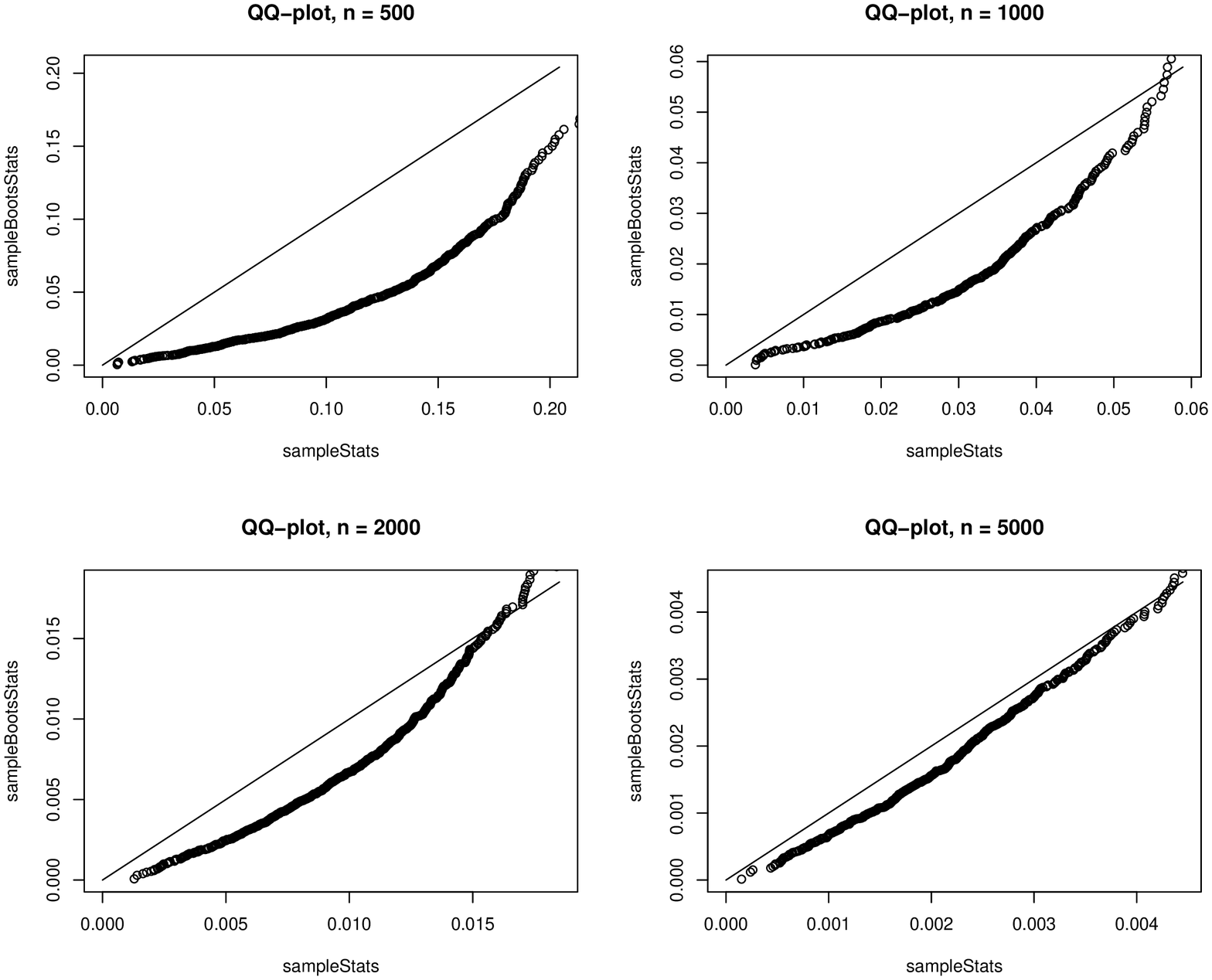}
    \caption{QQ-plot of a sample of the test statistic $\bar \Tc_2^c$ and a sample of the bootstrap test statistic $(\bar \Tc_2^c)^{*}$ using the non-parametric bootstrap for the gaussian copula, with different sample sizes and under $\bar \Hc_0^c$ (conditional Kendall's tau is constant and equal to $0.5$).}
    \label{fig:QQplot_bar_T2c_bootNP}
\end{figure}

\begin{figure}
    \centering
    \includegraphics[width = 12cm]{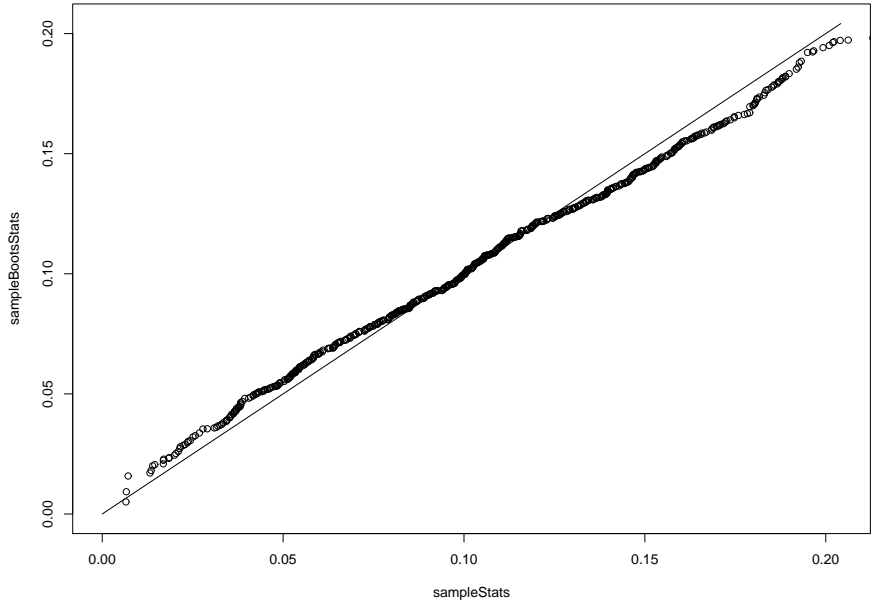}
    \caption{QQ-plot of a sample of the test statistic $\bar \Tc_{2,c}$ and a sample of the bootstrap test statistic $(\bar \Tc_{2,c})^{**}$ using the parametric independent bootstrap for the gaussian copula, with $n = 500$ and under $\bar \Hc_0^c$ (conditional Kendall's tau is constant and equal to $0.5$).}
    \label{fig:QQplot_bar_T2c_bootPI}
\end{figure}

Concerning tests of $\Hc_0$, the results are relatively satisfying. For the nonparametric tests and those based on the independence property (Tables~\ref{NP_tests_level} and~\ref{NP_tests_power}) the rejection rates are large when $\tau_{\max}=1$, and the theoretical levels ($5\%$) are underestimated (a not problematic feature in practice).
This is still the case for tests of the simplifying assumption under a parametric copula model through $\Tc_2^c$: see Tables~\ref{Param_tests_level} and~\ref{Param_tests_power}. The three bootstrap schemes provide similar numerical results. Remind that the bootstrapped statistics is $\left(\Tc_2^c \right)^{**}$ with bootPI (Remark \ref{tricky_schemes}). Tests of $\bar\Hc_0$ under a parametric framework and through $\bar\Tc_2^c$ confirm such observations. 
{\color{black} To evaluate the accuracy of the bootstrap approximations asymptotically, we have compared the empirical distribution of some test statistics and their bootstrap versions under the null hypothesis for two bootstrap schemes (see Figures~\ref{fig:QQplot_bar_T2c_bootNP} and~\ref{fig:QQplot_bar_T2c_bootPI}). For the nonparametric bootstrap, the two distributions begin to match each other at $n = 5 000$ whereas $n = 500$ is enough for the parametric independent bootstrap.}

\mds

\begin{figure}
    \centering
    \includegraphics[width = 11cm]{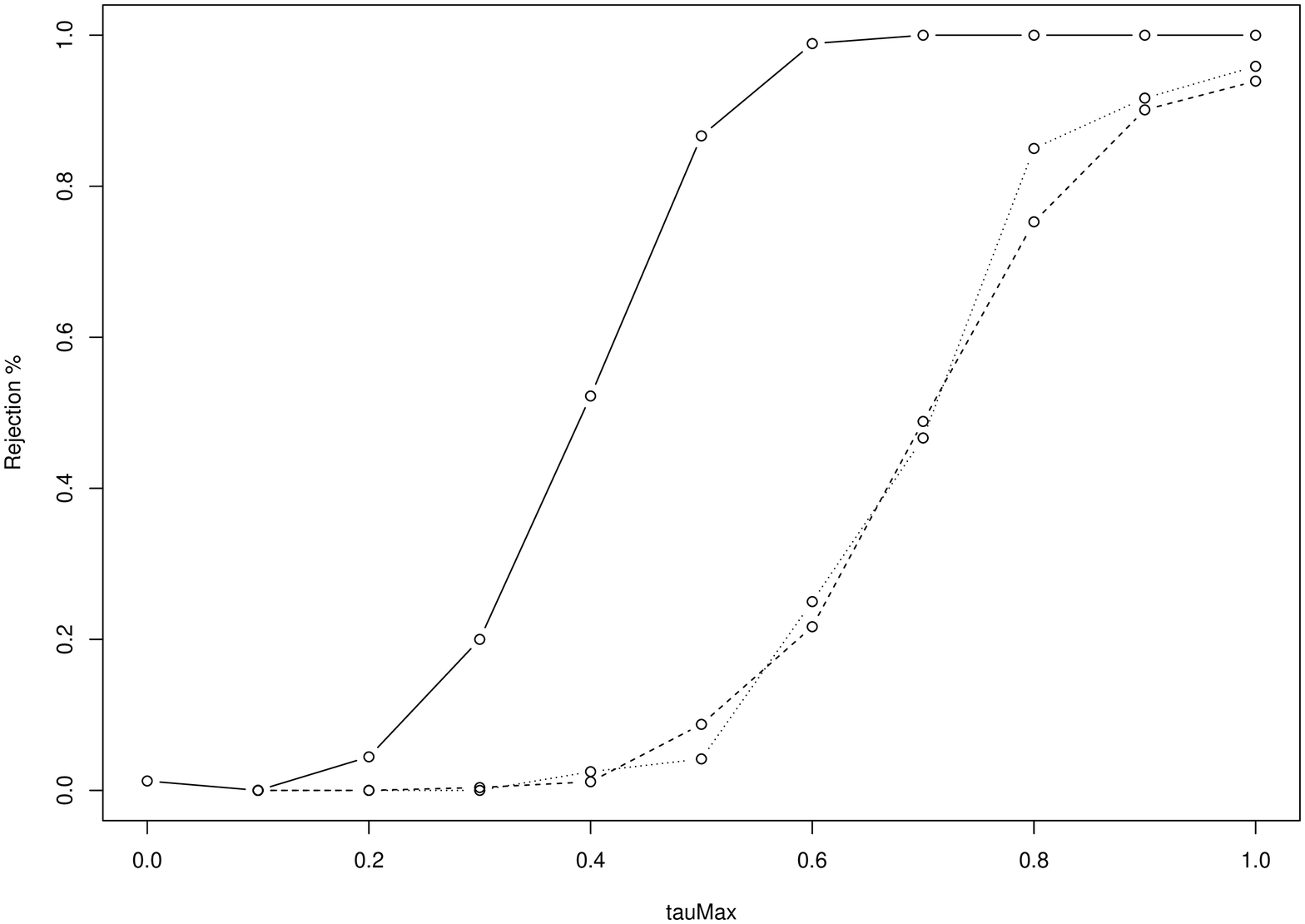}
    \caption{Rejection percentages for the statistics $\Ic_\chi$ (\ref{IcChi}) as a function of $\tau_{max}$: we use the gaussian copula, with a conditional parameter $\theta(x)$ calibrated such that
    the conditional Kendall's tau $\tau(x)$ satisfies $\tau(x) = \tau_{max} \cdot \Phi(x)$. Solid line: bootNP. Dashed line : bootPseudoInd. Dotted line : bootCond.} 
    \label{fig:I_chi_tau_max}
\end{figure}

\begin{figure}
    \centering
    \includegraphics[width = 11cm]{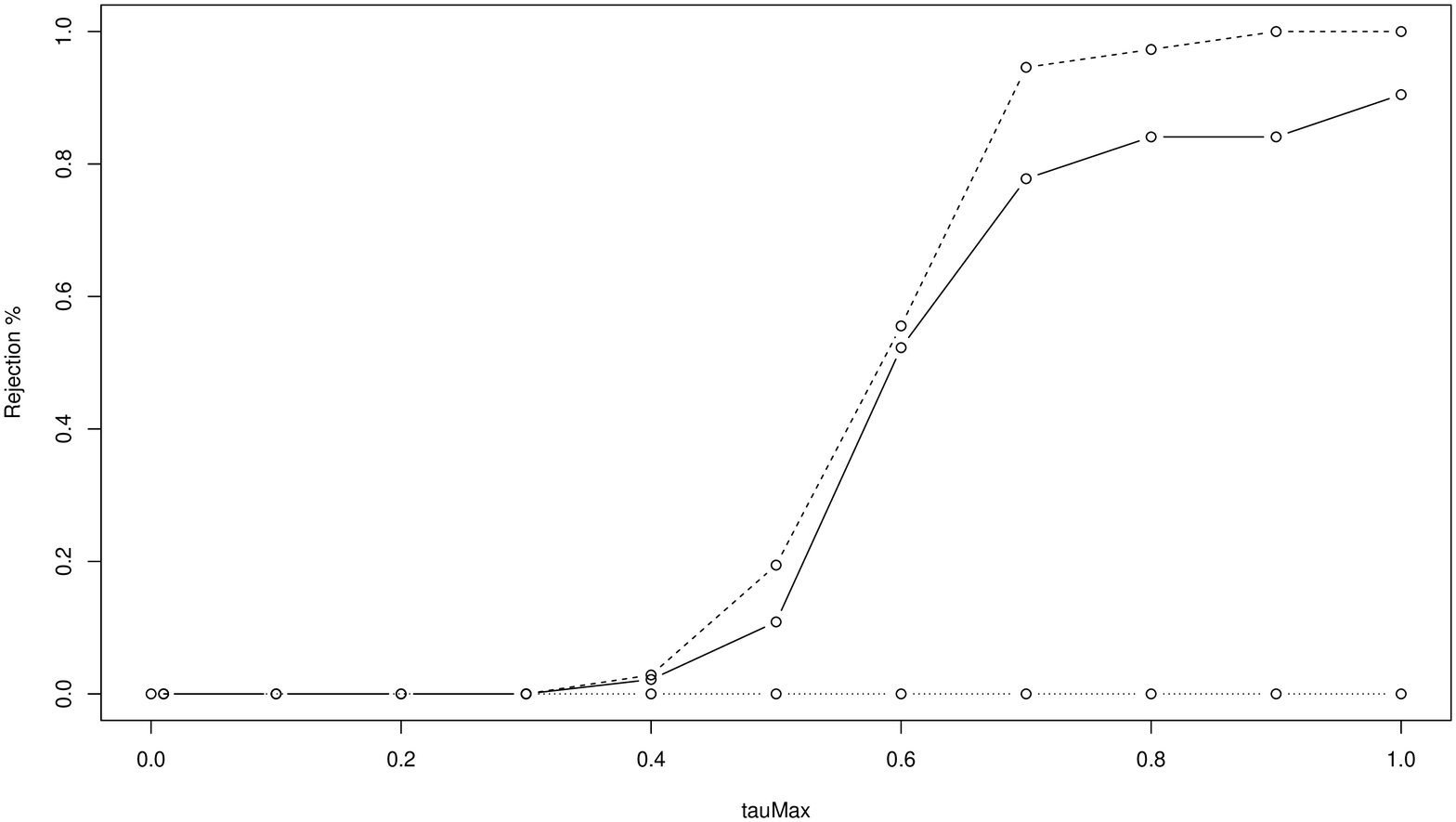}
    \caption{Rejection percentages for the statistics $I_{2,n}$ (\ref{Ic2n}) as a function of $\tau_{max}$: we use the gaussian copula, with a conditional parameter $\theta(x)$ calibrated such that
    the conditional Kendall's tau $\tau(x)$ satisfies $\tau(x) = \tau_{max} \cdot \Phi(x)$. Solid line: bootNP. Dashed line : bootPseudoInd. Dotted line : bootCond.} 
    \label{fig:I2n_tau_max}
\end{figure}

\begin{figure}
    \centering
    \includegraphics[width = 11cm]{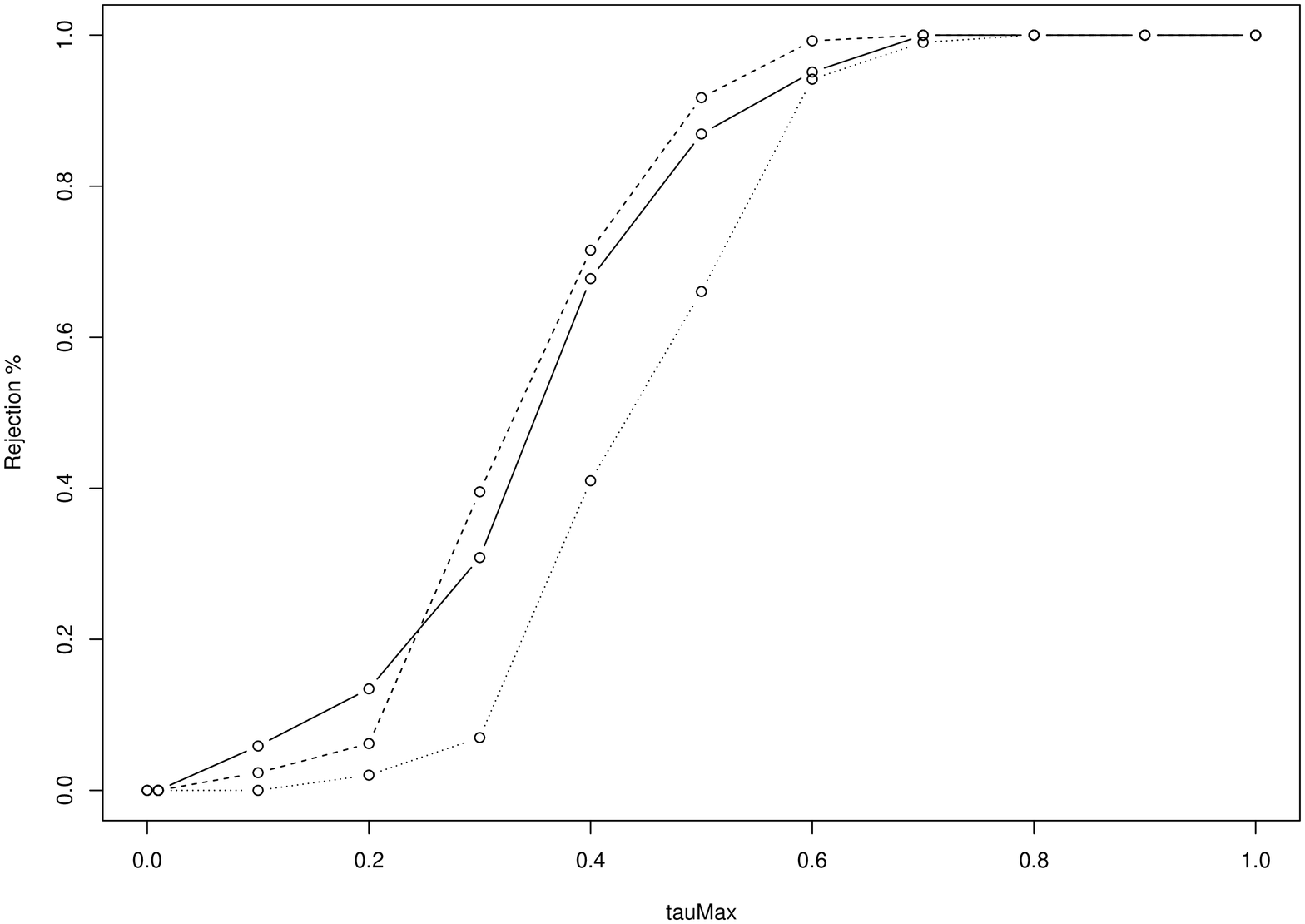}
    \caption{Rejection percentages for the statistics $\Tc_2^c$ (\ref{Tc_infty_C}) as a function of $\tau_{max}$: {\color{black} we use the gaussian copula, with a conditional parameter $\theta(x)$ calibrated} such that
    the conditional Kendall's tau $\tau(x)$ satisfies $\tau(x) = \tau_{max} \cdot \Phi(x)$. Solid line: bootNP. Dashed line : bootPI. Dotted line : bootPC.} 
    \label{fig:T2c_tau_max}
\end{figure}

\begin{figure}
    \centering
    \includegraphics[width = 11cm]{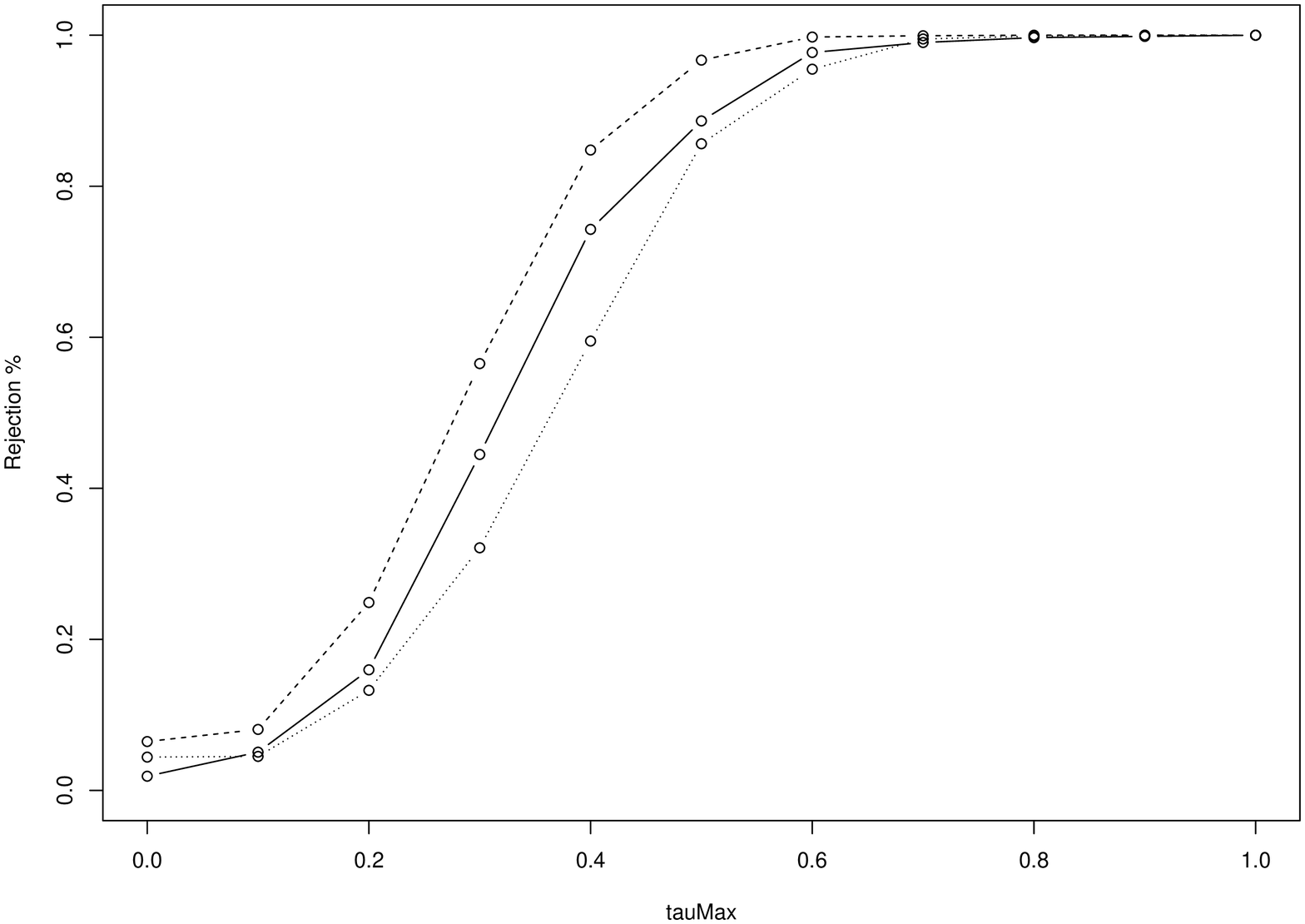}
    \caption{Rejection percentages for the statistics $\bar \Tc_2^c$ (\ref{bar_Tc_KS_c}) as a function of $\tau_{max}$: we use the gaussian copula, with a conditional parameter $\theta(x)$ calibrated such that
    the conditional Kendall's tau $\tau(x)$ satisfies $\tau(x) = \tau_{max} \cdot \lfloor m\Phi(X_3)  \rfloor /m $. Solid line: bootNP. Dashed line : bootPI. Dotted line : bootPC.} 
    \label{fig:bar_T2c_tau_max}
\end{figure}

We have tested the influence of $\tau_{\max}$: the smaller is this parameter, the smaller is the percentage of rejections under the alternative, because the simulated model tends to induce lower dependencies of copula parameters w.r.t. $X_3$: see Figures~\ref{fig:I_chi_tau_max},
\ref{fig:I2n_tau_max}, \ref{fig:T2c_tau_max}, and~\ref{fig:bar_T2c_tau_max}.
{\color{black} Note that, on each of these figures, the point at the left corresponds to a conditional Kendall's tau which is constant, and equal to $0$ (because $\tau_{\max}=0$) whereas the rejection percentages in Tables~\ref{NP_tests_level} and~\ref{Param_tests_level} correspond to a conditional Kendall's tau constant, and equal to $0.5$. As the two data-generating process are not the same, the rejection percentages can differ even if both are under the null hypothesis. Nevertheless, in every case, our empirical sizes converge to $0.05$ as the sample size increases. When $n=5000$, we found that the percentage of rejections are between $4$\% and $6$\%.}

\mds

\begin{figure}
    \centering
    \includegraphics[width = 12cm]{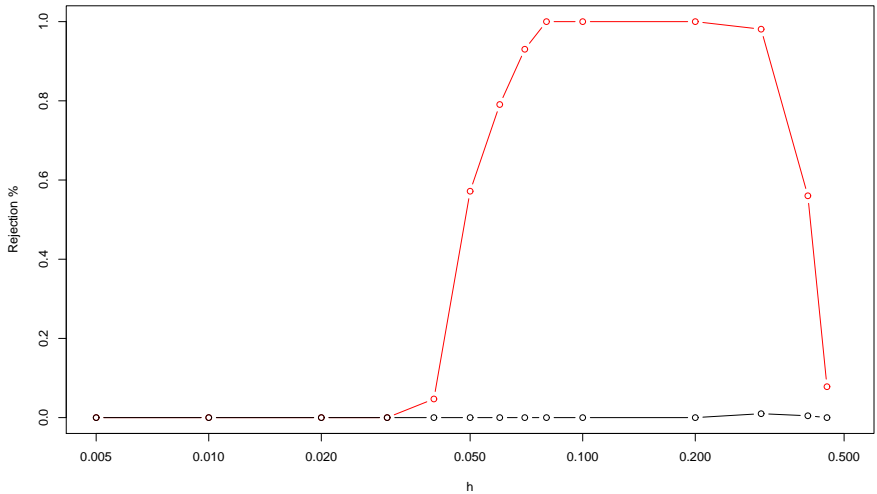}
    \caption{Rejection percentage for the statistic $\Tc^{0,m}_{CvM,n}$ (\ref{TcOCvMm}) with $m=20$ as a function of $h$. The red (resp. black) line corresponds to the alternative (resp. zero) assumption.}
    \label{fig:Rejection_T0_CvM_h}
\end{figure}

We have not tried to exhibit an ``asymptotically optimal'' bandwidth selector for our particular testing problem.
This could be the task for further research. We have preferred a basic ad-hoc procedure.
In our test statistics, we smooth w.r.t. $F_3(X_3)$ (or its estimate, to be specific), whose law is uniform on $(0,1)$.
A reasonable bandwidth $h$ is given by the so-called rule-of-thumb in kernel density estimation, i.e. $h^* = \sigma(F_3(X_3)) / n^{1/5} = 1/ (\sqrt{12} n^{1/5}) = 0.083$.
Such a choice has provided reasonable results. The typical influence of the bandwidth choice on the test results is illustrated in Figure~\ref{fig:Rejection_T0_CvM_h}.
In general, the latter $h^*$ belongs to reasonably wide intervals of ``convenient'' bandwidth values, so that the performances of our 
considered tests are not very sensitive to the bandwidth choice.

\mds

{\color{black} To avoid boundary problems, we have slightly modified the test statistics: we remove the observations $i$ such that $F_3(X_{i,3}) \leq h$ or $F_3(X_{i,3}) \geq 1-h$. This corresponds to changing the integrals (resp. max) on $[0,1]$ to integrals (resp. max) on $[h, 1-h]$.
}

\begin{table}[htb]
    \centering
    \begin{tabular}{c|cccc} 
        Family  & $\Tc_{CvM,n}^{0}$ (\ref{TcOCvMm}) 
        & $\tilde \Tc_{CvM,n}^{0}$ (\ref{TcOCvM_bis})
        & $\Ic_{\chi,n}$ (\ref{IcChi})
        & $\Ic_{2,n}$ (\ref{Ic2n})   \\
        \hline
        Gaussian &    0     &        0    &     0           &   0           \\
        Student  &    0     &        0    &     0           &   0           \\
        Clayton  &    0     &        0    &     1           &   0           \\
        Gumbel   &    1     &        1    &     0           &   1           \\
        Frank    &    0     &        0    &     0           &   0           \\
    \end{tabular}
    \caption{Rejection percentages under the null (nonparametric tests, nonparametric bootstrap bootNP).}
    \label{NP_tests_level}
\end{table}

\begin{table}[htb]
    \centering
    \begin{tabular}{c|cccc}
        Family  & $\Tc_{CvM,n}^{0}$ (\ref{TcOCvMm}) 
        & $\tilde \Tc_{CvM,n}^{0}$ (\ref{TcOCvM_bis})
        & $\Ic_{\chi,n}$ (\ref{IcChi})
        & $\Ic_{2,n}$ (\ref{Ic2n})   \\
        \hline
        Gaussian &    98    &   100       &   100           &   93          \\
        Student  &   100    &    99       &    98           &   90          \\
        Clayton  &    99    &    99       &    99           &   98          \\
        Gumbel   &    99    &    98       &   100           &   95          \\
        Frank    &    98    &   100       &    98           &   50          \\
    \end{tabular}
    \caption{Rejection percentages under the alternative (nonparametric tests, nonparametric bootstrap bootNP).}
    \label{NP_tests_power}
\end{table}

\begin{table}[htb]
    \centering
    \begin{tabular}{c|ccc|ccc}
        \multirow{2}{*}{Family}
        & \multicolumn{3}{c|}{$\Tc_2^c$ (\ref{Tc_infty_C})}
        & \multicolumn{3}{c}{$\bar \Tc_2^c$ (\ref{bar_Tc_KS_c})} \\
                 &  bootPI &  bootPC   &  bootNP &  bootPI &  bootPC & bootNP  \\
        \hline
        Gaussian &     4   &     0     &     0   &      6  &  4 &      1     \\
        Student  &     6   &     0     &     2   &      4  &  5 &      3     \\
        Clayton  &     7   &     0     &     1   &      7  &  1 &      1     \\
        Gumbel   &     3   &     1     &     0   &      9  &  2 &      2     \\
        Frank    &     4   &     0     &     6   &      3  &  5 &      1     \\
    \end{tabular}
    \caption{Rejection percentages under the null (parametric tests).}
    \label{Param_tests_level}
\end{table}

\begin{table}[htb]
    \centering
    \begin{tabular}{c|ccc|ccc}
        \multirow{2}{*}{Family}
        & \multicolumn{3}{c|}{$\Tc_2^c$ (\ref{Tc_infty_C})}
        & \multicolumn{3}{c}{$\bar \Tc_2^c$ (\ref{bar_Tc_KS_c})} \\
                 &  bootPI &  bootPC   &  bootNP &  bootPI &  bootPC & bootNP     \\
        \hline
        Gaussian &   100   &   100     &   100   &  100  &  100   & 100        \\
        Student  &   100   &   100     &   100   &  100  &  100   & 100        \\
        Clayton  &   100   &    62     &    98   &  100  &   98   & 100        \\
        Gumbel   &   100   &   100     &    34   &  100  &   99   &  76        \\
        Frank    &   100   &   100     &   100   &  100  &  100   & 100        \\
    \end{tabular}
    \caption{Rejection percentages under the alternative (parametric tests).}
    \label{Param_tests_power}
\end{table}

\FloatBarrier

\FloatBarrier

\section{Conclusion}
We have provided an overview of the simplifying assumption problem, under a statistical point of view.
In the context of nonparametric or parametric conditional copula models (with unknown conditional marginal distributions), numerous testing procedures have been proposed.
We have developed the theory towards a slightly different but related approach, where ``box-type'' conditioning events replace pointwise ones. This open a new field for research that is interesting per se.
Several new bootstrap procedures have been detailed, to evaluate p-values under the zero assumption in both cases.
In particular, we have proved the validity of one of them (the ``parametric independent'' bootstrap scheme under $\bar\Hc_0$).

\mds

Clearly, there remains a lot of work. We have opened the Pandora box rather than provided definitive answers. Open questions are still numerous:
precise theoretical convergence results of our test statistics (and others!), validity of these new bootstrap schemes, bandwidth choices, empirical performances,...
All these dimensions would require further research. We have made a contribution to the landscape of problems related to the simplifying assumption, and proposed a working program for the whole copula community.

\bigskip

\newpage

\appendix

{\color{black}

\section{Notations}
\label{Section_notations}

\mds

\begin{table}[htb]
    \centering
    \begin{tabular}{l l}
        $\X=(\X_I, \X_J)$ & random vector of size $d$ \\
        $I$, $J$ & $\{1, \dots, p\}$ and $\{p+1, \dots, d\}$ \\
        $\Sc=(\X_{1,1:d}, \dots, \X_{n,1:d})$ & initial sample of $n$ i.i.d. observations \\
        
        $A_J$ & measurable subset in $\RR^{d-p}$ \\
        $\Ac_J$ & collection of all measurable subsets of $\RR^{d-p}$ \\ 
        & \hspace{0.5cm} such that $\X_J$ is in each set with positive probability \\
        $\bar \Ac_J = \{ A_{1,J}, \dots, A_{m,J} \}$ & partition of $\RR^{d-p}$ into $m$ sets\\ 
        & \hspace{0.5cm}  such that $\X_J$ is in each set with positive probability \\
        $Y$ & box index, i.e. $Y$ is the $k$ such that $\X_J \in A_{k,J}$ \\
        $\hat U_{i,k}$ & $i$-th pseudo-observation of the $k$-th variable \\
        $\Z_{I|J}$ & conditional observation of $\X_I$ given $\X_J$ \\
        $\Z_{I|Y}$ & conditional observation of $\X_I$ given the box index $Y$ \\
        \\ \hline \\
        $\Cc$ & copula family indexed by the elements of a set $\Theta$ \\
        $C_\theta$ & copula of the family $\Cc$ with the parameter $\theta \in \Theta$ \\
        $c_\theta$ & density of the copula $C_\theta$ \\
        $\theta_0$ & unconditional parameter of the copula of $\Z_{I|J}$ \\
        $\theta (\x_J)$ & parameter of the conditional copula of $\Z_{I|J}$ given $\X_J = \x_J$ \\
        $\theta_0^b$ & unconditional parameter of the copula of $\Z_{I|Y}$ \\
        $\theta(A_J)$ & conditional parameter of the copula of $\Z_{I|Y}$ given $X_J \in A_J$ \\
        \\ \hline \\
        $F_i (\cdot)$ & marginal cdf of $X_i$, $i=1,\dots,d$ \\
        $F_{i|J}( \cdot | \X_J \in A_J)$ & conditional marginal cdf of $X_i$ given $\X_J \in A_J$, $i=1,\dots,p$ \\
        $F_{i|J}( \cdot | \X_J = \x_J)$ & conditional marginal cdf of $X_i$ given $\X_J = \x_J$, $i=1,\dots,p$ \\
        $F_{I|J}( \cdot | \X_J \in A_J)$ & conditional joint cdf of $\X_I$ given $\X_J \in A_J$ \\
        $F_{I|J}( \cdot | \X_J = \x_J)$ & conditional joint cdf of $\X_I$ given $\X_J = \x_J$ \\
        $G_{I,J} (\cdot)$ & joint cdf of $(\Z_{I|J}, \X_J)$ \\
        $C_{I|J}^{A_J}( \cdot | \X_J \in A_J)$ & conditional copula of $\X_I$ given $\X_J \in A_J$ \\
        $C_{I|J}( \cdot | \X_J = \x_J)$ & conditional copula of $\X_I$ given $\X_J = \x_J$ \\
        $C_{s,I|J} (\cdot)$ & simplified copula of $\X_I$ given $\X_J$ \\
    \end{tabular}
    \caption{Table of notations}
\end{table}

\begin{table}[htb]
    \centering
    \begin{tabular}{l l}
        $\Tc^0_{CvM,n}$ (\ref{TcOCvM})
        & brute-force test statistic of $\Hc_0$, constructed with the $L_2$ distance  \\
        & \hspace{0.5cm} between the conditional and the simplified copula \\
        \hspace{0.5cm} (resp. $\Tc^0_{KS,n}$ (\ref{Tc0KS}))
        & \hspace{0.5cm} (resp. $L_\infty$ distance) \\
        \hspace{0.5cm} (resp. $\Tc^{0,m}_{CvM,n}$ (\ref{TcOCvMm}))
        & \hspace{0.5cm} (resp. $L_2$ distance using a fixed number $m$ of points) \\
        $\tilde \Tc^0_{CvM,n}$ (\ref{TcOCvM_bis})
        & brute-force test statistic of $\Hc_0$, constructed with the $L_2$-distance \\
        \hspace{0.5cm} (resp. $\tilde \Tc^0_{KS,n}$ (\ref{Tc0KS_bis}))
        & \hspace{0.5cm} (resp. $L_\infty$-distance) between all pairs of conditional copulas \\
        \\ & \\
        $\Ic_{\chi,n}$ (\ref{IcChi})
        & chi-square-type test statistic of the independence between $\hat \Z_{I|J}$ and $\X_J$ \\
        $\Ic_{KS,n}$ (\ref{IcKS})
        & test statistic based on the distance between the joint empirical cdf \\
        & \hspace{0.5cm} of $(\hat \Z_{I|J}, \X_J)$ and the product of their empirical cdf, using the $L_\infty$ norm \\
        \hspace{0.5cm} (resp. $\Ic_{2,n} $ (\ref{Ic2n}))
        & \hspace{0.5cm} (resp. using the $L_2$ norm) \\
        \hspace{0.5cm} (resp. $\Ic_{CvM,n}$ (\ref{IcCvM}))
        & \hspace{0.5cm} (resp. using the $L_2$ norm, weighted by the joint empirical cdf as weight) \\
        \\ \hline \\
        $\Tc^c_{\infty}$ (\ref{Tc_infty_C})
        & test statistic based on the $L_{\infty}$ distance between the parameter of the \\
        &\hspace{0.5cm} conditional copula and the constant parameter of the simplified copula \\
        \hspace{0.5cm} (resp. $\Tc^c_2$ (\ref{Tc_infty_C}))
        & \hspace{0.5cm} (resp. $L_2$ distance) \\
        \hspace{0.5cm} (resp. $\Tc^c_{dist}$ (\ref{Tcdfcomp}))
        & \hspace{0.5cm} (resp. using some distance between the estimated copulas) \\
        \hspace{0.5cm} (resp. $\Tc^c_{dens}$ (\ref{Tdensitycomp}))
        & \hspace{0.5cm} (resp. using the $L_2$ distance between the estimated copula densities) \\
        \\ \hline \\
        $\bar\Tc_{dist,n}$ (\ref{bar_Tc_dist})
        & brute-force test statistic of $\bar \Hc_0$ constructed with the distance $dist(\cdot,\cdot)$ \\
        & \hspace{0.5cm} between all pairs of conditional copulas with Borelian subsets\\
        \hspace{0.5cm} (resp. $\bar\Tc_{KS,n}$ (\ref{bar_Tc_KS}))
        & \hspace{0.5cm} (resp. with the $L_\infty$ distance) \\
        \hspace{0.5cm} (resp. $\bar\Tc_{CvM,n}$ (\ref{bar_Tc_CvM}))
        & \hspace{0.5cm} (resp. with the $L_2$ distance) \\
        \\ \hline \\
        $\bar\Tc_{\infty}^{c}$ (\ref{bar_Tc_KS_c})
        & test statistic based on the $L_{\infty}$ distance between the parameters \\
        & \hspace{0.5cm} estimated on each set and the simplified parameter \\
        \hspace{0.5cm} (resp. $\bar\Tc_{2}^{c}$ (\ref{bar_Tc_KS_c}))
        & \hspace{0.5cm} (resp. based on the $L_2$ distance) \\
        \hspace{0.5cm} (resp. $\bar\Tc_{dist}^{c}$ (\ref{bar_Tc_dist_c}))
        & \hspace{0.5cm} (resp. based on some distance between the copulas whose parameters are \\
        & \hspace{0.5cm} estimated on each set and the copula with the simplified parameter) \\
        \\ \hline \\
        $\Tc^*$, $\Tc^{**}$ & bootstrap statistics corresponding to a general test statistic $\Tc$
    \end{tabular}
    \caption{Table of main test statistics}
\end{table}

\FloatBarrier

}

\section{Proof of Theorem~\ref{thm:ParIndepBoot}}
\label{Section_Proofs}

\subsection{Preliminaires}

\mds

Let $(\Z_i)_{i=1, \dots, n}$ be a sequence of i.i.d random vectors in $[0,1]^p$, $\Z_i$ being drawn from the true cdf $C_{\theta_0}$. They have the same law as
 the previously called vectors $\Z_{i,I|A_J}$  or $\Z_{i,I|Y}$ under the zero assumption $\bar\Hc^c_0$.
Let $(\X_{i,J})_{i=1, \dots, n}$ be a sequence of i.i.d random vectors in $\RR^{d-p}$, $\X_{i,J}\sim F_J$.
Let $(\Z_i^*)_{i=1, \dots, n}$ be an independent sequence of i.i.d random vectors in $[0,1]^p$, where $\Z_i^* \sim C_{\theta_0}$ exactly as $\Z_i$.
The three samples $(\Z_i)$, $(\X_{i,J})$ and $(\Z_i^*)$ are mutually independent.
Let $(\X_{i,J}^*)_{i=1, \dots, n}$ be a sequence of i.i.d random vectors in $\RR^{d-p}$, which are drawn from $F_{n,J}$, the empirical cdf of $\X_{1,J}, \dots \X_{n,J}$, and independently of both $(\Z_i)$ and $(\Z_i^*)$.

\mds

In the following, we shall use the notation $f \otimes g := (x,y) \mapsto f(x) g(y)$ when $f$, $g$ are two real functions, possibly from different spaces.
Set $ l(\theta, \cdot) := \log c_\theta (\cdot) $. We will need some conditions of regularity.

\mds

{\it Assumption (R):} $(\theta,\u_I) \mapsto l(\theta,\u_I)$ is three times differentiable  with respect to $\theta$, for every $\u_I\in (0,1)^p$. Moreover, for every $\epsilon>0$,
$$ \EE\left[
\sup_{\theta | \|\theta - \theta_0\|\leq \epsilon} \sup_{\{\z | \| \z - \Z_{i}\|\leq  \| \hat\Z_{i,I|Y} - \Z_{i}\| \} } ||
\dfrac{\partial^3 l}{\partial \theta^3}
\left( \theta , \z \right) ||  \right]  <+\infty.$$

\mds

The latter technical assumption can be weakened through some trimming techniques, as in Fermanian and Lopez (2015). Since this would require to change the definitions of the parametric estimators, we do not try to improve towards this direction.
We will set $\bigdot{c}_{\theta} := \partial c_{\theta} / \partial \theta$
and $\bigdott{c}_{\theta} :=\partial^2 c_{\theta}/\partial \theta^2$.

\mds

We associate to every $\X_{i,J}$ (resp. $\X_{i,J}^*$) its corresponding index $Y_i$ (resp.
$Y_i^* $) s.t. $ \X_{i,J}\in A_{Y_i}$ (resp. $ \X_{i,J}^*\in A_{Y^*_i}$).
For convenience, we assume that $(A_k)_{k=1,\ldots,m}$ is a partition of $\RR^{d-p}$. Otherwise, we have to restrict our sample to the observations for which $X_{i,J}$ belongs to some ``box'' $A_k$, $k=1,\ldots,m$.
Therefore, denote by $C_n$, $C_n^*$, $P_{n,Y}$ and $P_{n,Y}^*$ the empirical laws of $(\Z_i)$, $(\Z_i^*)$, $(Y_i)$ and $(Y_i^*)$ respectively.
The joint law of $(\Z_1,Y_{1})$ (resp. $(\Z_1,\X_{1,J})$) will be denoted by $\bar G := C_{\theta_0} \otimes P_{Y}$ (resp. $\bar G := C_{\theta_0} \otimes F_J$), with
$P_Y(k)=\PP(Y=k)$, $k=1,\ldots,m$.
Denote by $G_n$ (resp. $\bar G_n$) the empirical law of 
$( \Z_{i} , Y_i )_{i=1, \dots, n}$ (resp. $( \Z_{i} , \X_{i,J} )_{i=1, \dots, n}$)
Moreover, $G_n^*$ and $\bar G_n^*$ will be the empirical distributions of
$( \Z_{i}^* , Y_i^* )_{i=1, \dots, n}$ and $( \Z_{i}^* , \X_{i,J}^* )_{i=1, \dots, n}$ respectively.
Let $\Pc_n$ be the joint probability distribution of
$$\big( \Z_{i} , Y_i,\Z^*_{i} , Y_i^* \big)_{i=1, \dots, n}
\in \big([0 , 1]^{p}
\times \{1, \dots, m \} \big)^{\otimes 2n}.$$

\mds

The following proposition is key. It will be proved in Subsection \ref{proof_bootstrapped_Donsker_th}.
\begin{prop}
    \label{prop:bootstrapped_Donsker_th}
    Consider the empirical process defined on $[0,1]^{p}\times \RR^{d-p}$ by
    $$    \bar \GG _n (\z, \x_J) := \sqrt{n} (\bar G _n - \bar G) (\z, \x_J) :=
        \frac{1}{\sqrt{n}} \sum_{i=1}^n \{\1 \big( (\Z_i,\X_{i,J}) \leq (\z, \x_J) \big)
        - C_{\theta_0} (\z) F_J (\x_J)\},$$
    and the corresponding bootstrapped empirical process
    \begin{align*}
         \bar \GG _n^* (\z, \x_J)
        &:=    \frac{1}{\sqrt{n}} \sum_{i=1}^n \1 \big( (\Z_i^*,\X_{i,J}^*) \leq (\z, \x_J) \big)
        - C_{\theta_0} (\z)
        \frac{1}{\sqrt{n}} \sum_{i=1}^n
        \1 (\X_{i,J} \leq \x_J),
    \end{align*}
    or, equivalently, $\bar \GG _n^* = \sqrt{n}
    (\bar G _n^* - C_{\theta_0} \otimes  F_{n,J})$.
    Then there exist two independent and identically distributed Gaussian processes $\bar\AA_G$ and $\bar\AA_G^{\perp}$ such that
    $\big( \bar \GG _n \, , \, \bar \GG _n^* \big)$ converges to $ (\bar\AA_G \, , \, \bar\AA_G^{\perp})$
    weakly in $\Big( \ell^{\infty} \big([0,1]^p \times \RR^{d-p} \big) \Big)^2$.
\end{prop}

As a Corollary, we deduce the same results when the discrete variables $Y_i$ replace the variables $\X_{i,J}$.
\begin{prop}
    \label{prop:bootstrapped_Donsker_th_Y}
    Under the assumptions of Proposition~\ref{prop:bootstrapped_Donsker_th}, let the empirical process defined on $[0,1]^{p}\times \{1,\ldots,m\}$ by
    \begin{align*}
     \GG _n (\z, k)
        &:= \sqrt{n} (G _n - G) (\z, k) \\
        &:=
        \frac{1}{\sqrt{n}} \sum_{i=1}^n \{ \1 \big( \Z_i \leq \z, Y_{i}=k \big)
        - C_{\theta_0} (\z) P_Y (k) \},
     \end{align*}
    and its bootstrapped empirical process
    \begin{align*}
         \GG _n^* (\z, k)
        &:=    \frac{1}{\sqrt{n}} \sum_{i=1}^n \1 \big( \Z_i^* \leq \z, Y_{i}^*=k \big)
        - C_{\theta_0} (\z)
        \frac{1}{\sqrt{n}} \sum_{i=1}^n
        \1 (Y_i =k),
    \end{align*}
    or equivalently $  \GG _n^* = \sqrt{n}
    ( G _n^* - C_{\theta_0} \otimes  P_{n,Y})$, $P_{n,Y}(k)$ being the empirical proportion of $\Sc_n$-observations into $A_k$.
Then, there exist two independent and identically distributed processes $\AA_G$ and $\AA_G^{\perp}$ such that
    $\big( \GG _n \, , \, \GG _n^* \big)$ converges to $ (\AA_G \, , \, \AA_G^{\perp})$
    weakly in $\Big( \ell^{\infty} \big( [0,1]^p    \times \{1,\ldots,m\} \big) \Big)^2$.
\end{prop}

\begin{rem}
\label{CovAAG}
The covariance function of $\AA_G$ (or $\AA^{\perp}$) is given by
\begin{eqnarray*}
\lefteqn{ \EE[ \AA_G (\z,y) \AA_G (\z',y')]= \lim_n \EE[  \GG_n(\z,y)\GG_n(\z',y') ]      }\\
&=&
\1(y=y') \PP(Y=y) C_{\theta_0}(\z \wedge \z') -  \PP(Y=y)\PP(Y=y')C_{\theta_0}(\z) C_{\theta_0}(\z').
\end{eqnarray*}
\end{rem}

\mds
As a ``toolbox'', we will need the following lemma.
\begin{lemma}
    \label{lemma:properties_Gn_G}
    Let $\hat \theta_{0}^b$ and $\hat \theta(A_{k})$ be the estimators based on the pseudo-sample $(\hat\Z_{i,I|Y},Y_i)_{i=1,\dots,n}$ (and then on the sample $(\Z_i,Y_i)_{i=1,\dots,n}$) as
$$    \hat\theta_{0}^b:= \arg\max_{\theta\in\Theta} \sum_{i=1}^n \log c_\theta (\hat\Z_{i,I|Y} ),\ \text{and}$$
$$    \hat\theta(A_k):= \arg\max_{\theta\in\Theta} \sum_{i=1}^n \log c_\theta (\hat\Z_{i,I|Y} ). \1(Y_{i}= k ), \; k=1,\ldots,m.$$
We will assume they lie in the interior of $\Theta$.
Set $     \Theta_{n,0} := \sqrt{n}\big( \hat \theta_{0}^b - \theta_0 \big)$, and, for $k=1, \dots, m$,
$\Theta_{n,k} := \sqrt{n}\big( \hat \theta(A_{k}) - \theta_0 \big)$.
Moreover, for any distribution $H$ on $[0,1]^p\times \{1,\ldots,m\}$, set
    \begin{equation*}
        \psi_{k,1} (H)
        := \int \frac{\partial l}{\partial \theta}
        \left(\theta_0,
        \left(
        \dfrac{\int \1
        \{ z_q^1 \leq z_q^2, y^1 = y^2 \} dH(\z^1, y^1)}
        {\int \1 \{ y^1 = k \} dH(\z^1, y^1) }
        \right)_{q=1, \dots, p}
        \right) \1 \{ y^2 = k \} \, dH(\z^2, y^2),
    \end{equation*}
    \begin{equation*}
        \psi_{k,2} (H)
        := \int \frac{\partial^2 l}{\partial \theta^2}
        \left(\theta_0,
        \left(
        \dfrac{\int \1
        \{ z_q^1 \leq z_q^2, y^1 = y^2 \} dH(\z^1, y^1)}
        {\int \1 \{ y^1 = k \} \, dH(\z^1, y^1) }
        \right)_{q=1, \dots, p}
        \right) \1 \{ y^2 = k \} \, dH(\z^2, y^2).
    \end{equation*}
    \begin{enumerate}
        \item[(i)] For $k=1, \dots, m$,
        \begin{equation*}
            \Theta_{n,k} = -\dfrac
            {\sqrt{n}  \psi_{k,1}(G_n) }{ \psi_{k,2}(G_n)}+ o_P(1).
        \end{equation*}
        \item[(ii)] For every discrete law $P_Y$ with values in $\{1, \dots, m\}$, the corresponding
        distribution $\tilde G := C_{\theta_0} \otimes P_Y$ satisfies $\psi_{k,1} (\tilde G) = 0$.

        \item[(iii)] $\psi_1 := (\psi_{1,1}, \dots, \psi_{m,1})$ is Hadamard-differentiable at every cdf $H$, and its differential is given by
    \end{enumerate}
    \begin{align*}
        &\bigdot{   \psi}_{k,1} (H) (h) =
        \int \frac{\partial l}{\partial \theta}
        \left(\theta_0,
        \left(
        \dfrac{\int \1
        \{ z_q^1 \leq z_q^2, y^1 = y^2 \} dH(\z^1, y^1)}
        {\int \1 \{ y^1 = k \} dH(\z^1, y^1) }
        \right)_{q=1, \dots, p}
        \right) \1 \{ y^2 = k \} \, dh(\z^2, y^2) \\
        &+ \sum_{j=1}^p
        \int \frac{\partial^2 l}{\partial \theta \, \partial z_j}
        \left(\theta_0,
        \left(
        \dfrac{\int \1
        \{ z_q^1 \leq z_q^2, y^1 = y^2 \} dH(\z^1, y^1)}
        {\int \1 \{ y^1 = k \} dH(\z^1, y^1) }
        \right)_{q=1, \dots, p}
        \right) \1 \{ y^2 = k \} \\
        & \cdot \left(
        \dfrac{\int \1
        \{ z_j^1 \leq z_j^2, y^1 = y^2 \} dh(\z^1, y^1)}
        {\int \1 \{ y^1 = k \} dH(\z^1, y^1) }
        - \dfrac{\int \1
        \{ z_j^1 \leq z_j^2, y^1 = y^2 \} dH(\z^1, y^1)
        \int \1 \{ y^1 = k \} dh(\z^1, y^1)}
        {\left(\int \1 \{ y^1 = y^2 \} dH(\z^1, y^1) \right)^2}
        \right)
        \, dH(\z^2, y^2) \\
    \end{align*}

    \saveenum
    \begin{enumerate}\resetenum
        \item[(iv)]
        \begin{equation*}
            \Theta_{n,0} =- \dfrac
            { \sum_{k=1}^m \sqrt{n} \big( \psi_{k,1}(G_n) \big) }
            {\sum_{k=1}^m \psi_{k,2}(G_n) }
            + o_P(1)
            = \dfrac
            { \sum_{k=1}^m
            \psi_{k,2}(G_n)  \Theta_{n,k} }
            { \sum_{k=1}^m \psi_{k,2}(G_n) }
            + o_P(1)
        \end{equation*}
    \end{enumerate}
\end{lemma}

\mds

\noindent
{ \it Proof :}
Note that $\hat\Z_{i,I|J}$ is an explicit measurable function of the sample $(\Z_{i,I|J})_{i=1,\ldots,n}$. Indeed,
for any $i=1, \dots, n$ and $q = 1, \dots, p$,
\begin{align}
    \hat Z_{i,q|Y}
    :=& \hat F_{n,q} ( X_{i,q}
    | \X_{J}\in A_{Y_i , J}) \nonumber \\
    :=& \dfrac{\sum_{j=1}^n \1
    \{ X_{j,q} \leq X_{i,q},    \X_{j,J} \in A_{Y_i, J}
    \}}
    {\sum_{j=1}^n \1
    \{ \X_{j,J} \in A_{k(\X_{i,J}) , J} \} }\nonumber  \\
    =& \dfrac{\sum_{j=1}^n \1
    \{ F_q(X_{j,q} |\X_{J} \in A_{Y_j,J})
    \leq F_q(X_{i,q} |\X_{J} \in A_{Y_j,J}),
    Y_j=Y_i \}}
    {\sum_{j=1}^n \1
    \{ Y_j=Y_i \} } \nonumber \\
    =& \dfrac{\sum_{j=1}^n \1
    \{ F_q(X_{j,q} |\X_{J} \in A_{Y_j,J})
    \leq F_q(X_{i,q} |\X_{J} \in A_{Y_i,J}),Y_j=Y_i\}}
    {\sum_{j=1}^n \1 \{ Y_j=Y_i \} } \nonumber \\
    =& \dfrac{\sum_{j=1}^n \1
    \{ Z_{j,q|Y} \leq Z_{i,q|Y},
    Y_j = Y_i
    \}}
    {\sum_{j=1}^n \1
    \{ Y_j = Y_i \} }.
    \label{hatZZ}
\end{align}
We deduce that $\hat \theta_{0}^b$ and $\hat \theta(A_k)$ are measurable functions of the unobservable random variables $\Z_{i,I|Y}$ and $Y_i$, for $i=1, \dots, n$.

\mds

(i). Let $k \in \{ 1, \dots, m \}$.
Applying successively the first order condition for the estimator $\hat \theta(A_k)$ and some Taylor series expansions, we have
\begin{align*}
    0
    &= \frac{1}{n} \sum_{i=1}^n
    \frac{\partial l}{\partial \theta}
    (\hat\theta(A_k), \hat\Z_{i,I|J}) \1 \{ Y_i = k \} \\
    &= B_n^{1,k} - B_n^{2,k} \big( \hat \theta(A_{k,J}) - \theta_0 \big)
    + o_P\big( \hat \theta(A_{k,J}) - \theta_0 \big),\; \text{with}
\end{align*}
\begin{equation*}
    B_n^{1,k} := \frac{1}{n} \sum_{i=1}^n
    \frac{\partial l}{\partial \theta}
    \big(\theta_0,\hat\Z_{i,I|J} \big)
    \1 \{ Y_i = k \}
    \text{ and }
    B_n^{2,k} := - \frac{1}{n} \sum_{i=1}^n
    \frac{\partial^2 l}{\partial \theta^2}
    \big(\theta_0,  \hat\Z_{i,I|J} \big)
    \1 \{ Y_i = k \},
\end{equation*}
implying
$$\Theta_{n,k}
:= \sqrt{n} \big( \hat \theta(A_{k,J}) - \theta_0 \big)
= \frac{\sqrt{n} B_n^{1,k}} {B_n^{2,k}}
+ o_P \big( \Theta_{n,k} \big).$$
Now, invoking~(\ref{hatZZ}), let us compute the numerator of this expression:
\begin{align*}
    B_n^{1,k} &= \frac{1}{n} \sum_{i=1}^n
    \frac{\partial l}{\partial \theta}
    \left(\theta_0,
    \left(
    \dfrac{\sum_{j=1}^n \1
    \{ Z_{j,q} \leq Z_{i,q},  Y_q = k \}}
    {\sum_{j=1}^n \1
    \{ Y_j = k \} }
    \right)_{q=1, \dots, p}
    \right) \1 \{ Y_i = k \} \\
    &= \int \frac{\partial l}{\partial \theta}
    \left(\theta_0,
    \left(
    \dfrac{\int \1
    \{ z_q^1 \leq z_q^2, y^1 = k \}
    dG_n(\z^1, y^1)}
    {\int \1 \{ y^1 = k \} dG_n(\z^1, y^1) }
    \right)_{q=1, \dots, p}
    \right) \1 \{ y^2 = k \}
    dG_n(\z^2, y^2) \\
    &= \psi_{k,1} (G_n).
\end{align*}

In the same way, the denominator can be rewritten as
\begin{align*}
    B_n^{2,k}
    &= - \int \frac{\partial^2 l}{\partial \theta^2}
    \left(\theta_0,
    \left(
    \dfrac{\int \1
    \{ z_q^1 \leq z_q^2, y^1 = k \}
    dG_n(\z^1, y^1)}
    {\int \1 \{ y^1 = k \} dG_n(\z^1, y^1) }
    \right)_{q=1, \dots, p}
    \right) \1 \{ y^2 = k \}
    dG_n(\z^2, y^2) \\
    &= - \psi_{k,2} (G_n).
    \end{align*}

\mds

(ii). We now prove the second part of the lemma.
Since $\tilde G = C_{\theta_0} \otimes F_Y$, we get
\begingroup
\allowdisplaybreaks
\begin{align*}
    \psi_{k,1} (\tilde G)
    &:= \int
    \frac{\partial l}{\partial \theta}
    \left(\theta_0,
    \left(
    \dfrac{\int \1
    \{ z_q^1 \leq z_q^2, y^1 = k \}
    d \tilde G(\z^1, y^1)}
    {\int \1 \{ y^1 = k \} d \tilde G(\z^1, y^1) }
    \right)_{q=1, \dots, p}
    \right) \1 \{ y^2 = k \}
    d \tilde G(\z^2, y^2) \\
    &= \int
    \frac{\partial l}{\partial \theta}
    \left(\theta_0,
    \left(
    \dfrac{ \PP \{ Z_q^1 \leq z_q^2,
    Y^1 = k \} }
    {\PP \{ Y^1 = k \}  }
    \right)_{q=1, \dots, p}
    \right) \1 \{ y^2 = k \}
    d \tilde G(\z^2, y^2) \\
    &= \int
    \frac{\partial l}{\partial \theta}
    \left(\theta_0,
    \left( \PP \{ Z_q^1 \leq z_q^2 \} \right)_{q=1, \dots, p}
    \right) dC_{\theta_0}(\z^2)
    \int \1 \{ y^2 = k \} dF_Y(y^2) \\
    &= \PP \{ Y = k \}
    \int
    \frac{\partial l}{\partial \theta}
    \left(\theta_0, \z^2 \right) dC_{\theta_0}(\z^2)= 0.
\end{align*}%
\endgroup

\mds

(iii). We remark that the law $G$ appears three times in $\psi_{k,1}$: two times in the log-density $l$ and one time at the end of the main integral.
By separating the effect of a change from $H$ to $H+h$ in the main integral only (first term of the differential) and the effect of a change in $l$, and
using the standard rule of differential calculus ($l$ is differentiable), we obtain the second part of the given result.

\mds

(iv). As in the proof of (i), we apply successively the first order condition for $\hat \theta_{n,0}^b$ and some Taylor series expansion to get
$$    0= \frac{1}{n} \sum_{i=1}^n
    \frac{\partial l}{\partial \theta}
    (\hat \theta_{0}^b,  \hat\Z_{i,I|Y})
    = B_n^1 - \big( \hat \theta_{0}^b - \theta_0 \big) B_n^2
    + o_P\big( \hat \theta_{0}^b - \theta_0 \big),\;\text{with}$$
$$B_n^1 := \frac{1}{n} \sum_{i=1}^n
\frac{\partial l}{\partial \theta}
\big(\theta_0,  \hat\Z_{i,I|Y} \big)
= \sum_{k=1}^m B_n^{1,k}
\text{ and }
B_n^2 := - \frac{1}{n} \sum_{i=1}^n
\frac{\partial^2 l}{\partial \theta^2}
\big(\theta_0,  \hat\Z_{i,I|Y} \big)
= \sum_{k=1}^m B_n^{2,k}.$$
We deduce
\begin{align*}
    \Theta_{n,0}
    &:= \sqrt{n} \big( \hat\theta_{0}^b - \theta_0 \big)
    = \frac{\sqrt{n} B_n^1} {B_n^2}
    + o_P\big( \Theta_{n,0} \big)
    = \dfrac{\sqrt{n} \sum_{k=1}^m B_n^{1,k}}
    {\sum_{k=1}^m B_n^{2,k}}
    + o_P\big( \Theta_{n,0} \big) \\
    &= \dfrac{\sqrt{n} \sum_{k=1}^m \psi_{k,1} (G_n)}
    {\sum_{k=1}^m - \psi_{k,2} (G_n)}
    + o_P\big( \Theta_{n,0} \big) \\
    &= \dfrac{ \sum_{k=1}^m
     \psi_{k,2} (G_n) \Theta_{n,k}}
    {\sum_{k=1}^m  \psi_{k,2} (G_n)}
    + o_P\big( \Theta_{n,0} \big) .\;\;    \Box
\end{align*}

\mds

\begin{lemma}
    \label{decomp_l_n}

    Let $\ell_n$ be defined by
    \begin{equation*}
        \ell_n := \sum_{i=1}^n \log \left(
        \dfrac{c_{\hat \theta_{0}^b}(\Z_i^*)}
        {c_{\theta_{0}}(\Z_i^*)} \right).
    \end{equation*}
    If there exists a random vector $\Theta_0$ such that
    $\Theta_{n,0} \Longrightarrow \Theta_0$ under $\Pc_n$,
    then we have
    \begin{align*}
        \ell_n = \Theta_{0}^T \WW^{\perp}
        - \frac{1}{2} \Theta_{0}^T I_0 \Theta_{0}
        + o_P \left( 1 \right),
    \end{align*}
    where $\WW^{\perp} \sim \Nc ( 0, I_0 )$ is independent of the sample $\big( \Z_{i,I|Y} , Y_i \big)_{i=1, \dots, n}$ and $I_0$ is the Fisher information matrix
    \begin{equation*}
        I_0 := \EE_{C_{\theta_0}} \left[ \dfrac{
        \bigdot{c}_{\theta_0}^T (\Z) \, \bigdot{c}_{\theta_0} (\Z)
        } {c_{\theta_0}^2 (\Z)} \right].
    \end{equation*}
\end{lemma}

\mds

\noindent
{ \it Proof :}
By a Taylor expansion, we obtain
\begin{align*}
    \ell_n
    &= \sum_{i=1}^n  \{ l \left( \hat\theta_{0}^b , \Z_i^* \right)
    - l \left( \theta_0 , \Z_i^* \right) \} \\
    &= \big( \hat\theta_{0}^b - \theta_0 \big)^T
    \sum_{i=1}^n \dfrac{\partial l}{\partial \theta}
    \left( \theta_0 , \Z_i^* \right)
    + \frac{1}{2} \big( \hat\theta_{0}^b - \theta_0 \big)^T
    \sum_{i=1}^n \dfrac{\partial^2 l}{\partial \theta^2}
    \left( \theta_0 , \Z_i^* \right)
    \big( \hat\theta_{0}^b - \theta_0 \big)
    + R_n \\
    &= \Theta_{n,0}^T \left[ \frac{1}{\sqrt{n}}
    \sum_{i=1}^n \dfrac{\partial l}{\partial \theta}
    \left( \theta_0 , Z_i^* \right) \right]
    + \frac{1}{2} \Theta_{n,0}^T \left[ \frac{1}{n}
    \sum_{i=1}^n \dfrac{\partial^2 l}{\partial \theta^2}
    \left( \theta_0 , Z_i^* \right) \right]
    \Theta_{n,0}
    + R_n.
\end{align*}
First, we have
\begin{align*}
    R_n
    &\leq Cst || \hat\theta_{0}^b - \theta_0 ||^3
    \sup_{\theta | \|\theta - \theta_0\|\leq \|\hat\theta_{0}^b - \theta_0\|} || \sum_{i=1}^n \dfrac{\partial^3 l}{\partial \theta^3}
    \left( \theta , Z_i^* \right) || \\
    &\leq Cst || \Theta_{n,0} ||^3
    \cdot \sup_{\theta | \|\theta - \theta_0\|\leq \|\hat\theta_{0}^b - \theta_0\|} ||
     \frac{1}{n}\sum_{i=1}^n \dfrac{\partial^3 l}{\partial \theta^3}
    \left( \theta , Z_i^* \right) || \cdot
    \frac{1}{\sqrt{n}} = O_P \left( \frac{1}{\sqrt{n}}\right),
\end{align*}
by Assumption (R).
By the usual CLT, we know that
   $ \frac{1}{\sqrt{n}}
    \sum_{i=1}^n \partial l/ \partial \theta \left( \theta_0 , Z_i^* \right)
    \longrightarrow \WW^{\perp}.$
$\WW^{\perp}$ is independent of
$\big( \Z_{i,I|Y} , Y_i \big)_{i=1, \dots, n}$
as a limit of a sequence of variables that have the same property.
Using the law of large numbers, we have also
\begin{align*}
    \frac{1}{n} \sum_{i=1}^n
    \dfrac{\partial^2 l}{\partial \theta^2}
    \left( \theta_0 , Z_i^* \right)
    = \frac{1}{n} \sum_{i=1}^n
    \dfrac{\bigdott{c}_{\theta}} {c_{\theta}} ( Z_i^* )
    - \dfrac{\bigdot{c}_{\theta}^T \bigdot{c}_{\theta}} {c_{\theta}^2} ( Z_i^* )
    \Longrightarrow 0 - I_0 \;\;\Box
\end{align*}

\mds

\subsection{Proof of Theorem \ref{thm:ParIndepBoot}}
We first reason under $\Pc_n$ as in Theorem 1 in Genest and Rémillard (2008).
By Proposition \ref{prop:bootstrapped_Donsker_th_Y},
under $\Pc_n$, there exist two independent and identically distributed processes $\AA_G$ and $\AA_G^{\perp}$ such that
\begin{equation*}
    \sqrt{n} \Big(
    G_n - C_{\theta_0} \otimes P_Y \, , \,
    G_n^* - C_{\theta_0} \otimes P_{n,Y}
    \Big) \Longrightarrow (\AA_G \, , \, \AA_G^{\perp}),
\end{equation*}
weakly in $\Big( \ell^{\infty} \big( [0,1]^p    \times \{1,\ldots,m\} \big) \Big)^2$.
By (iii) of Lemma \ref{lemma:properties_Gn_G}, $\psi_1$ is Hadamard-differentiable and so, using the functional Delta-method, we deduce
\begin{equation*}
    \sqrt{n} \Big(
    \psi_1(G_n) - \psi_1(C_{\theta_0} \otimes P_Y) \, , \,
    \psi_1(G_n^*) - \psi_1(C_{\theta_0} \otimes P_{n,Y})
    \Big) \Longrightarrow
    \Big( \bigdot{\psi}_{1} (G) (\AA_G) \, , \,
    \bigdot{\psi}_{1} (G) (\AA_G^{\perp}) \Big).
\end{equation*}
By (ii) of Lemma \ref{lemma:properties_Gn_G},
$\psi_1(C_{\theta_0} \otimes P_Y)
= \psi_1(C_{\theta_0} \otimes P_{n,Y}) = 0$, implying
\begin{align*}
    \sqrt{n} \Big(
    \psi_{1,1}(G_n), & \dots, \psi_{m,1}(G_n) \, , \, \psi_{1,1}(G_n^*), \dots, \psi_{m,1}(G_n^*)
    \Big) \\
    &\Longrightarrow \Big(
    \bigdot{\psi}_{1,1} (G) (\AA_G), \dots,
    \bigdot{\psi}_{m,1} (G) (\AA_G) \, , \,
    \bigdot{\psi}_{1,1} (G) (\AA_G^{\perp}), \dots,
    \bigdot{\psi}_{m,1} (G) (\AA_G^{\perp})
    \Big).
\end{align*}
By Slutsky's theorem, we have
\begin{align*}
    \sqrt{n} \Big(
    \frac{\psi_{1,1}(G_n)}{ \psi_{1,2}(G_n)}, &\dots,
    \frac{\psi_{m,1}(G_n)}{ \psi_{m,2}(G_n)}
    \, , \,
    \frac{\psi_{1,1}(G_n^*)}{ \psi_{1,2}(G_n^*)}, \dots,
    \frac{\psi_{m,1}(G_n^*)}{ \psi_{m,2}(G_n^*)}
    \Big) \\
    &\Longrightarrow \Big(
    \frac{\bigdot{\psi}_{1,1} (G) (\AA_G)}
    { \psi_{1,2}(G)}, \dots,
    \frac{\bigdot{\psi}_{m,1} (G) (\AA_G)}
    { \psi_{m,2}(G)}\, , \,
    \frac{\bigdot{\psi}_{1,1} (G) (\AA_G^{\perp})}
    { \psi_{1,2}(G)} , \dots,
    \frac{\bigdot{\psi}_{m,1} (G) (\AA_G^{\perp})}
    { \psi_{m,2}(G)} \Big).
\end{align*}
By (i) of Lemma \ref{lemma:properties_Gn_G}, the latter convergence result implies
\begin{align*}
    \Big( \Theta_{n,1}&, \dots, \Theta_{n,m} \, , \,
    \Theta_{n,1}^*, \dots, \Theta_{n,m}^* \Big)\\
    &\Longrightarrow \Big(
    \frac{\bigdot{\psi}_{1,1} (G) (\AA_G)}
    {- \psi_{1,2}(G)}, \dots,
    \frac{\bigdot{\psi}_{m,1} (G) (\AA_G)}
    {- \psi_{m,2}(G)}\, , \,
    \frac{\bigdot{\psi}_{1,1} (G) (\AA_G^{\perp})}
    {- \psi_{1,2}(G)} , \dots,
    \frac{\bigdot{\psi}_{m,1} (G) (\AA_G^{\perp})}
    {- \psi_{m,2}(G)} \Big) \\
    &=: \Big( \Theta_1 , \dots , \Theta_m \, , \,
    \Theta_1^{\perp} , \dots , \Theta_m^{\perp} \Big).
\end{align*}
Moreover, $\Big( \Theta_1 , \dots , \Theta_m \Big)$
and $\Big(\Theta_1^{\perp} , \dots , \Theta_m^{\perp} \Big)$
are independent and identically distributed under $\Pc_n$, by construction.
Because of (iv) of Lemma \ref{lemma:properties_Gn_G},
$\Theta_{n,0}$ can asymptotically be seen as a mean of the $\Theta_{n,k}$ and this provides
\begin{align*}
    \Theta_{n,0} = \dfrac
    { \sum_{k=1}^m \psi_{k,2}(G_n) \Theta_{n,k} }
    { \sum_{k=1}^m \psi_{k,2}(G_n) }
    \Longrightarrow \dfrac
    { \sum_{k=1}^m  \psi_{k,2}(G)  \Theta_{k} }
    { \sum_{k=1}^m \psi_{k,2}(G) } =: \Theta_0.
\end{align*}
Therefore, by the continuous mapping theorem, we deduce
\begin{align*}
    \Big( \Theta_{n,0}, &\dots, \Theta_{n,m} \, , \,
    \Theta_{n,0}^*, \dots, \Theta_{n,m}^* \Big)
    \Longrightarrow
    \Big( \Theta_0 , \dots , \Theta_m \, , \,
    \Theta_0^{\perp} , \dots , \Theta_m^{\perp} \Big),
\end{align*}
and we still have that $\Big( \Theta_0 , \dots , \Theta_m \Big)$
and $\Big(\Theta_0^{\perp} , \dots , \Theta_m^{\perp} \Big)$
are independent and identically distributed under $\Pc_n$.

\mds

Now, we will work under $\Pc_n^*$ the probability measure over
$\big([0 , 1]^{p}
\times \{1, \dots, m \} \big)^{\otimes 2n}$ whose density with respect to $\Pc_n$ is
\begin{align*}
    \dfrac{d\Pc_n^*}{d\Pc_n}
    ( \z_1, y_1, \dots, \z_n, y_n,
    \z^*_1, y_1^*, \dots, \z^*_n, y_n^* )
    &= \prod_{i=1}^n \dfrac{c_{\hat\theta_{0}^b}(\z_i^*)}
    {c_{\theta_{0}}(\z_i^*)},
\end{align*}
where $\hat\theta_{0}^b$ is the estimator of $\theta_0$ when applied to the ``sample'' $(\z_1, y_1, \dots, \z_n, y_n)$.
We remark that
$$\dfrac{d\Pc_n^*}{d\Pc_n}(\Z_{1},Y_1,\ldots,\Z_{n},Y_n,\Z^*_{1},Y_1^*,\ldots,\Z^*_{n},Y^*_n) = \exp(\ell_n).$$

Since we have shown that $\Theta_{n,0} \Longrightarrow \Theta_0$ under $\Pc_n$, use Lemma \ref{decomp_l_n} and obtain
\begin{align*}
    \ell_n = \Theta_{0}^T \WW^{\perp}
    - \frac{1}{2} \Theta_{0}^T I_0 \Theta_{0}
    + o_P \left( 1 \right).
\end{align*}

Therefore, under $\Pc_n$, we have
\begin{align*}
    \left( \dfrac{d\Pc_n^*}{d\Pc_n},
    \Theta_{n,0}, \dots, \Theta_{n,m},
    \Theta^*_{n,0}, \dots, \Theta^*_{n,m} \right)
    \Longrightarrow
    \left( \zeta,
    \Theta_{0}, \dots, \Theta_{m},
    \Theta^{\perp}_0, \dots, \Theta^{\perp}_{m} \right),
\end{align*}
where $    \zeta := \exp \left( \Theta_0^T \WW^{\perp}- \Theta_0^T I_0 \Theta_0 / 2 \right)$.
Note that $\EE[\zeta] = \EE[ \EE[\zeta | \Theta_0] ] = 1$
because  $\Theta_0$ and $\WW^{\perp}$ are independent, and $\WW^{\perp} \sim \Nc(0, I_0)$.
This corresponds to condition (iii) of Theorem 3.10.5 of Van der Waart and Wellner (1996), and we deduce $\Pc_n^*$ is contiguous with respect to $\Pc_n$.
We can then apply Le Cam's Third Lemma (Theorem 3.10.7 of Van der Waart and Wellner (1996)), and we get that, under $\Pc_n^*$,
\begin{align*}
    \left(\Theta_{n,0}, \dots, \Theta_{n,m},
    \Theta^*_{n,0}, \dots, \Theta^*_{n,m} \right)
    \Longrightarrow
    \left( \tilde\Theta_{0}, \dots, \tilde\Theta_{m},
    \Theta^*_0, \dots, \Theta^*_{m} \right),
\end{align*}
where $\EE [\chi(\tilde\Theta_{0:m},\Theta^*_{0:m})] =    \EE [\zeta \chi(\Theta_{0:m},\Theta^{\perp}_{0:m})]$ for any simple function $\chi$.
Choose $w_1$ and $w_2 \in \RR^{m+1}$ and set $\Sigma := Var \big[\Theta_{0:m} \big]$.
Then, we have
\begin{align*}
    &\EE[\exp(i w_1^T \tilde\Theta_{0:m}
    + i w_2^T \Theta^*_{0:m} )] = \EE[\zeta \exp(i w_1^T \Theta_{0:m}
    + i w_2^T \Theta^{\perp}_{0:m} )] \\
    &= \EE[\exp(\Theta_0^T \WW^{\perp}
    - \Theta_0^T I_0 \Theta_0 / 2
    + i w_1^T \Theta_{0:m}
    + i w_2^T \Theta^{\perp}_{0:m} )] \\
    &= \EE \Big[\exp( i w_1^T \Theta_{0:m}
    - \Theta_0^T I_0 \Theta_0 / 2) \,
    \EE[ \exp(\Theta_0^T \WW^{\perp}
    + i w_2^T \Theta^{\perp}_{0:m} )
    \, | \, \Theta_{0:m}] \Big] \\
    &= \EE \left[\exp( i w_1^T \Theta_{0:m}
    - \Theta_0^T I_0 \Theta_0 / 2)
    \exp \left(\frac{1}{2} \left(
    - w_2^T \Sigma w_2 + \Theta^T_0 I_0 \Theta_0
    + 2 i w_2 \EE[{\Theta^{\perp}_{0:m}}^T \WW^{\perp}] \Theta_0
    \right) \right) \right] \\
    &= \EE \left[\exp \left(
    i w_1^T \Theta_{0:m} - w_2^T \Sigma w_2 / 2
    +  i w_2 \EE[{\Theta^{\perp}_{0:m}}^T \WW^{\perp} ] \Theta_0
    \right) \right] \\
    &= \EE \left[ \exp \left(
    i w_1^T \Theta_{0:m} + i w_2 \Theta^{\perp}_{0:m}
    +  i w_2 \EE[{\Theta^{\perp}_{0:m}}^T \WW^{\perp} ] \Theta_0
    \right) \right].
\end{align*}

Therefore, we have proven the following equality:
\begin{equation*}
    \left( \tilde\Theta_{0}, \dots, \tilde\Theta_{m},
    \Theta^*_0, \dots, \Theta^*_{m} \right)
    \stackrel{\text{law}}{=}
    \left( \Theta_{0}, \dots, \Theta_{m},
    \Theta^{\perp}_0 + a_0 \Theta_{0}, \dots, \Theta^{\perp}_{m} + a_m \Theta_{0} \right),
\end{equation*}
where $a_k = \EE[{\Theta^{\perp}_k}^T \WW^{\perp}]$.
To finish the proof, it remains to show that $a_k = a_0$ for all $k \in \{ 1, \dots, m \}$, i.e.
\begin{equation*}
\EE[{\Theta^{\perp}_0}^T \WW^{\perp}] = \EE[{\Theta^{\perp}_k}^T \WW^{\perp}].
\label{EndProof}
\end{equation*}

First, we know from the proof of Lemma \ref{lemma:properties_Gn_G} that
$\Theta_{k,n} = -\bigdot{\psi}_{k,1}(G)(\AA_G)/\psi_{k,2}(G) + o_P(1)$, $k=1,\ldots,m$
and $\Theta_{0,n} = -\bigdot{\psi}_{0,1}(G)(\AA_G)/\psi_{0,2}(G) + o_P(1)$, where
\begin{equation*}
    \psi_{0,1} (G)
    := \int \frac{\partial l}{\partial \theta}
    \left(\theta_0,
    \left( \dfrac{\int \1
    \{ z_q^1 \leq z_q^2, y^1 = y^2 \} dG(\z^1, y^1)}
    {\int \1 \{ y^1 = y^2 \} dG(\z^1, y^1) }
    \right)_{q=1, \dots, p} \right) \, dG(\z^2, y^2),
\end{equation*}
and
\begin{equation*}
    \psi_{0,2} (G)
    := \int \frac{\partial^2 l}{\partial \theta^2}
    \left(\theta_0,
    \left( \dfrac{\int \1
    \{ z_q^1 \leq z_q^2, y^1 = y^2 \} dG(\z^1, y^1)}
    {\int \1 \{ y^1 = y^2 \} dG(\z^1, y^1) }
    \right)_{q=1, \dots, p} \right) \, dG(\z^2, y^2).
\end{equation*}
This implies $\Theta_{k} = -\bigdot{\psi}_{k,1}(G)(\AA_G)/\psi_{k,2}(G)$, $k=0,\ldots,m.$

\mds

Actually, the reasoning is exactly the same when dealing with $\Theta_{k,n}^*$ and $\Theta_{k}^{\perp}$, $k=0,\ldots,m$, replacing $\AA_G$ by $\AA_G^{\perp}$.
We get
\begin{equation*}
    \Theta_k^{\perp} = -\dfrac{\bigdot{\psi}_{k,1} (G) (\AA_G^{\perp})}
    {\psi_{k,2} (G) },\;\text{and}\;
    \Theta_0^{\perp} = -\dfrac{\bigdot{\psi}_{0,1} (G) (\AA_G^{\perp})}
    {\psi_{0,2} (G) } \cdot
\end{equation*}

Second, note that, when $k=1,\ldots,m$,
\begin{align*}
    \psi_{k,2} (G)
    &:= \int \frac{\partial^2 l}{\partial \theta^2}
    \left(\theta_0,
    \left(
    \dfrac{\int \1
    \{ z_q^1 \leq z_q^2 , y^1 = k \} dG(\z^1, y^1)}
    {\int \1 \{ y^1 = k \} \, dG(\z^1, y^1) }
    \right)_{q=1, \dots, p}
    \right) \1 \{ y^2 = k \} \, dG(\z^2, y^2) \\
    &= \PP(Y = k)
    \int \frac{\partial^2 l}{\partial \theta^2}
    \left(\theta_0,
    \left( \int \1 \{ z_q^1 \leq z_q^2 \} dC_{\theta_0} (\z^1)
    \right)_{q=1, \dots, p}
    \right) \, dC_{\theta_0} (\z^2) \\
    &= \PP(Y = k)
    \int \frac{\partial^2 l}{\partial \theta^2}
    \left(\theta_0, \z \right) \, dC_{\theta_0} (\z^2) \\
    &= \PP(Y = k) \psi_{0,2} (G).
\end{align*}

Third, let us calculate $\bigdot{\psi}_{k,1} (G) (h)$, $k=0,1,\ldots,m$.
From Lemma \ref{lemma:properties_Gn_G}, we have
\begin{align*}
    &\bigdot{\psi}_{k,1} (G) (h) =
    \int \frac{\partial l}{\partial \theta}
    \left(\theta_0,\z^2
    \right) \1 \{ y^2 = k \} \, dh(\z^2, y^2) + \sum_{j=1}^p
    \int \frac{\partial^2 l}{\partial \theta \, \partial z_j}
    \left(\theta_0,\z^2
    \right) \1 \{ y^2 = k \} \\
    & \cdot \left(
    \dfrac{\int \1
    \{ z_j^1 \leq z_j^2, y^1 = y^2 \} dh(\z^1, y^1)}
    {\int \1 \{ y^1 = y^2 \} dG(\z^1, y^1) }
    - \dfrac{\int \1
    \{ z_j^1 \leq z_j^2, y^1 = y^2 \} dG(\z^1, y^1)
    \int \1 \{ y^1 = y^2 \} dh(\z^1, y^1)}
    {\left(\int \1 \{ y^1 = y^2 \} dG(\z^1, y^1) \right)^2}
    \right)
    \, dG(\z^2, y^2) .
\end{align*}
for $k=1,\ldots,m$.
Since $G = C_{\theta_0} \otimes F_Y$, we can simplify the latter equalities:
\begin{align*}
    &\bigdot{\psi}_{k,1} (G) (h) =
    \int \frac{\partial l}{\partial \theta}
    \left(\theta_0,\z^2
    \right) \1 \{ y^2 = k \} \, dh(\z^2, y^2)
    + \PP(Y=k) \sum_{j=1}^p
    \int \frac{\partial^2 l}{\partial \theta \, \partial z_j}
    \left(\theta_0, \z^2 \right) \\
    & \cdot \left(
    \dfrac{\int [\1
    \{ z_j^1 \leq z_j^2, y^1 = y^2 \}
    - z_j^2 \1 \{ y^1 = y^2 \}] dh(\z^1, y^1)}
    {\int \1 \{ y^1 = y^2 \} dG(\z^1, y^1) }
    \right)
    dC_{\theta_0}(\z^2).
\end{align*}
Since $\bigdot{\psi}_{0,1} (G) (h) = \sum_{k=1}^m \bigdot{\psi}_{k,1} (G) (h)$, we have
\begin{align*}
    &\bigdot{\psi}_{0,1} (G) (h) =
    \int \frac{\partial l}{\partial \theta}
    \left(\theta_0, \z^2 \right) \, dh(\z^2, y^2)
    + \sum_{j=1}^p
    \int \frac{\partial^2 l}{\partial \theta \, \partial z_j}
    \left(\theta_0, \z^2 \right) \\
    & \cdot \left(
    \dfrac{\int [\1
    \{ z_j^1 \leq z_j^2 , y^1 = y^2 \}
    - z_j^2 \1 \{ y^1 = y^2 \} ] dh(\z^1, y^1)}
    {\int \1 \{ y^1 = y^2 \} dG(\z^1, y^1) }
    \right)
    dC_{\theta_0}(\z^2) dF_Y(\y^2).
\end{align*}

Then, we can rewrite $\bigdot{\psi}_{k,1} (G) (h) = M_1(h,k) + \PP(Y=k) M_3(h,k)$
and $\bigdot{\psi}_{0,1} (G) (h) = M_2(h) + \sum_{k'=1}^m \PP(Y=k') M_3(h,k')$, where
\begin{align*}
    &M_1(h,k) := \int \frac{\partial l}{\partial \theta}
    \left(\theta_0, \z^2 \right)
    \1 \{ y^2 = k \} \, dh(\z^2, y^2), \quad
    M_2(h) := \int \frac{\partial l}{\partial \theta}
    \left(\theta_0, \z^2 \right) \, dh(\z^2, y^2), \\
    &M_3(h,k) := \sum_{j=1}^p
    \int \frac{\partial^2 l}{\partial \theta \, \partial z_j}
    \left(\theta_0, \z^2 \right)
    \left(
    \dfrac{\int [ \1
    \{ z_j^1 \leq z_j^2, y^1 = k \}
    - z_j^2 \1 \{ y^1 = k \} ] dh(\z^1, y^1)}
    {\int \1 \{ y^1 = k \} dG(\z^1, y^1) }
    \right)
    dC_{\theta_0}(\z^2).
\end{align*}

Substituting $h$ by $\AA_G^{\perp}$, we get
\begin{align*}
    \Theta_k^{\perp} &= -\dfrac{\bigdot{\psi}_{k,1} (G) (\AA_G^{\perp})}
    {\psi_{k,2} (G) }
    = -\frac{M_1(\AA_G^{\perp},k)} {\PP(Y = k) \psi_{0,2} (G)}
    -  \frac{M_3(\AA_G^{\perp},k)} {\psi_{0,2} (G)},
    \\
    \Theta_0^{\perp} &= -\dfrac{\bigdot{\psi}_{0,1} (G) (\AA_G^{\perp})}
    {\psi_{0,2} (G) }
    = -\frac{M_2(\AA_G^{\perp})} {\psi_{0,2} (G)}
    -  \frac{\sum_{k'=1}^m \PP(Y=k') M_3(\AA_G^{\perp},k')} {\psi_{0,2} (G)}\cdot
\end{align*}

Fourth, since $\WW^{\perp}$ is the weak limit of $ \sum_{i=1}^n \dfrac{\partial l}{\partial \theta} \left( \theta_0 , \Z_i^* \right)/\sqrt{n}$ under $\Pc_n$, this implies
$\WW^{\perp} = \bigdot \psi_{3} (G) (\AA_G^{\perp})$, with
\begin{align*}
    \psi_{3} (G) = \int \frac{\partial l}{\partial \theta}
    \left(\theta_0, \z\right) \, dG(\z, y),\;\text{and}\;
    \bigdot \psi_{3} (G) (h) = \int \frac{\partial l}{\partial \theta}
    \left(\theta_0, \z\right) \, dh(\z, y).
\end{align*}

Finally, by (i) of the following Lemma \ref{lemma:equalityExpectations}, we have
$\EE[M_1(\AA_G^{\perp},k)^T \WW^{\perp}]
= \PP(Y = k) \EE[M_2(\AA_G^{\perp})^T \WW^{\perp}]$,
By (ii) of the latter lemma, we have $\EE[M_3(\AA_G^{\perp},k)^T \WW^{\perp}]
= \EE[M_3(\AA_G^{\perp},k')^T \WW^{\perp}]$ for all $k$ and $k'$.
Finally, we obtain $\EE[{\Theta^{\perp}_0}^T \WW^{\perp}] = \EE[{\Theta^{\perp}_k}^T \WW^{\perp}]$, which finishes the proof. $\Box$

\mds

\begin{lemma}
    \label{lemma:equalityExpectations}
    Assume that $\bar \Hc_0^c$ is satisfied. Then,

    \begin{enumerate}
        \item[(i)] For $k=1, \dots, m$,
    \end{enumerate}
    \begin{align*}
        \EE &\left[ \int \frac{\partial l}{\partial \theta^T}
        \left(\theta_0, \z \right)
        \1 \{ y = k \} \, d\AA_G^{\perp}(\z, y)
        \int \frac{\partial l}{\partial \theta}
        \left(\theta_0, \z'\right) \, d\AA_G^{\perp}(\z', y')
        \right] \nonumber \\
        & \qquad = \PP(Y=k) \;
        \EE \left[ \int \frac{\partial l}{\partial \theta^T}
        \left(\theta_0, \z \right) \, d\AA_G^{\perp}(\z, y)
        \int \frac{\partial l}{\partial \theta}
        \left(\theta_0, \z' \right) \, d\AA_G^{\perp}(\z', y')
        \right].
    \end{align*}

    \saveenum
    \begin{enumerate} \resetenum
        \item[(ii)] The expectations
    \end{enumerate}
    \begin{align*}
        \EE &\left[
        \int \frac{\partial^2 l}{\partial \theta^T \, \partial z_j}
        \left(\theta_0, \z^2 \right)
        \left(
        \dfrac{\int [\1
        \{ z_j^1 \leq z_j^2, y^1 = k \}
        - z_j^2 \1 \{ y^1 = k \} ] d\AA_G^{\perp}(\z^1, y^1)}
        {\int \1 \{ y^1 = k \}  dG(\z^1, y^1) }
        \right) \right.
        dC_{\theta_0}(\z^2)  \\
        & \left. \cdot
        \int \frac{\partial l}{\partial \theta}
        \left(\theta_0, \z^3 \right) \, d\AA_G^{\perp}(\z^3, y^3)
        \right]
    \end{align*}
    do not depend on $k=1,\ldots,m$.
\end{lemma}

\mds

\noindent
{ \it Proof :}
(i) By simple calculations, we obtain
\begin{align*}
    &\EE \left[ \int \frac{\partial l}{\partial \theta^T}
    \left(\theta_0, \z \right)
    \1 \{ y = k \} \, d\AA_G^{\perp}(\z, y)
    \int \frac{\partial l}{\partial \theta}
    \left(\theta_0, \z'\right) \, d\AA_G^{\perp}(\z', y')
    \right] \\
    &= \int \frac{\partial l}{\partial \theta^T}
    \left(\theta_0, \z \right)
    \1 \{ y = k \}
    \frac{\partial l}{\partial \theta}
    \left(\theta_0, \z'\right)
    d_{\z,y,\   z',y'} \left( \EE \left[
    \AA_G^{\perp}(\z, y) \AA_G^{\perp}(\z', y')
    \right] \right) \\
    &= \int \frac{\partial l}{\partial \theta^T}
    \left(\theta_0, \z \right)
    \1 \{ y = k \}
    \frac{\partial l}{\partial \theta}
    \left(\theta_0, \z'\right)
    \left\{  \delta_{y'=y} d\PP(y) [ dC_{\theta_0}(\z) \delta_{\z'=\z} + dC_{\theta_0}(\z') \delta_{\z=\z'}]  \right. \\
    & \left.- dC_{\theta_0}(\z)  dC_{\theta_0}(\z') d\PP(y) d\PP(y')
    \right\} \\
    &= 2\PP(Y=k)\int \frac{\partial l}{\partial \theta^T}
    \left(\theta_0, \z \right)
    \frac{\partial l}{\partial \theta}
    \left(\theta_0, \z\right)
    dC_{\theta_0}(\z)
    -\PP(Y=k) \int \frac{\partial l}{\partial \theta^T}
    \left(\theta_0, \z \right) dC_{\theta_0}(\z)\cdot  \int \frac{\partial l}{\partial \theta}
    \left(\theta_0, \z \right) dC_{\theta_0}(\z).
\end{align*}
By summing up the latter identities w.r.t. $k=1,\ldots,m$, we prove (i).

\mds

(ii) For convenience, let us write $\phi_2(\z) := \partial^2 l / (\partial \theta^T \, \partial z_j)
\left(\theta_0, \z \right)$
and $\phi_3(\z) := \partial l (\theta_0, \z )/ \partial \theta^T$. We get the result if we prove that
\begin{align*}
    A_{1,k} &:= \EE \left[ \int  \phi_2(\z_2)
    \left(
    \dfrac{\int \1 \{ z_j^1 \leq z_j^2\} \1 \{y^1 = k \} d\AA_G^{\perp}(\z^1, y^1)}
    {\int \1 \{ y^1 = k \} dG(\z^1, y^1) }
    \right)
    dC_{\theta_0}(\z^2)
    \int  \phi_3(\z_3) \, d\AA_G^{\perp}(\z^3, y^3)
    \right]\;\text{and} \\
    A_{2,k} &:= \EE \left[ \int  \phi_2(\z_2)
    \left(
    \dfrac{\int z_j^2 \1 \{ y^1 = k \} d\AA_G^{\perp}(\z^1, y^1)}
    {\int \1 \{ y^1 = k \} dG(\z^1, y^1) }
    \right)
    dC_{\theta_0}(\z^2)
    \int  \phi_3(\z_3) \, d\AA_G^{\perp}(\z^3, y^3)
    \right]
\end{align*}
do not depend on $k$. We will do the task for $A_{1,k}$, $k=1,\ldots,m$, and the calculations will be similar for $A_{2,k}$. Note that
\begingroup
\allowdisplaybreaks
\begin{align*}
    A_{1,k} &= \dfrac{1}{\PP(Y=k)} \int \phi_2(\z_2)
    \1 \{ z_j^1 \leq z_j^2\} \1 \{y^1 = k \}
    \phi_3(\z_3) \, dC_{\theta_0}(\z^2)
    d_{\z^1, y^1, \z^3, y^3}
    \EE \left[ \AA_G^{\perp}(\z^1, y^1) \AA_G^{\perp}(\z^3, y^3) \right]  \\
    &= \dfrac{1}{\PP(Y=k)} \int \phi_2(\z_2)
    \1 \{ z_j^1 \leq z_j^2\} \1 \{y^1 = k \}
    \phi_3(\z_3) \, dC_{\theta_0}(\z^2) \\
    &  \left\{  \delta_{y^3=y^1} d\PP(y^1) [ dC_{\theta_0}(\z^1) \delta_{\z^3=\z^1} + dC_{\theta_0}(\z^3) \delta_{\z^1=\z^3}]
    - C_{\theta_0}(\z^1)  dC_{\theta_0}(\z^3) d\PP(y^1) d\PP(y^3) \right\} .
\end{align*}
\endgroup
We deduce that
\begin{align*}
    A_{1,k}
    &= 2\int \phi_2(\z_2)
    \1 \{ z_j^1 \leq z_j^2\} \phi_3(\z^1) \, dC_{\theta_0}(\z^1) dC_{\theta_0}(\z^2) \\
     &- \int \phi_2(\z_2)
    \1 \{ z_j^1 \leq z_j^2\} \phi_3(\z^3) \, dC_{\theta_0}(\z^1) dC_{\theta_0}(\z^2)  dC_{\theta_0}(\z^3),
\end{align*}
that does not depend on $k$. $\Box$

\subsection{Proof of Proposition \ref{prop:bootstrapped_Donsker_th}}
\label{proof_bootstrapped_Donsker_th}

\mds

As usual with the nonparametric bootstrap, we rewrite the bootstrapped empirical process by counting the number of times every observation of the initial sample is drawn:
$$d\bar G_n^* = \frac{1}{n} \sum_{i=1}^n M_{n,i} \delta_{(\Z_i^*,\X_{i,J})},$$
where $M_{n,i}$ denotes the number of times $(\X_{i,J}^*)$ has been redrawn in a $n$-size bootstrap resampling with replacement.
It is well-known that $M_n:=(M_{n,1},\ldots,M_{n,n})$ follows a multinomial distribution $\Mc(n, n^{-1}, \dots, n^{-1})$: its mean is $n$ and the associated probabilities are $1/n,\ldots,1/n$.
In other words, $ \bar\GG_n^*(\z,\x_J) = \frac{1}{\sqrt{n}} \sum_{i=1}^n  M_{n,i}  \{ \1((\Z_i^*,\X_{i,J})\leq (\z,\x_J)) - C_{\theta_0}(\z) F_{n,J}(\x_J) \}.$

\mds
We can remove the dependence between the random components $M_{n,i}$, $i=1,\ldots,n$ by a ``Poissonization'' procedure. We mimic van der Vaart and Wellner (1996), p.346: instead of drawing $n$ times the initial observations, this is done $N_n$ times, where $N_n$ follows a Poisson distribution with mean $n$ and $N_n$ is independent of the initial sample.
Then, the $n$ variables $M_{N_n,1},\ldots,M_{N_n,n}$ are i.i.d. Poisson random variables with mean one. And we can build the new process as
$$ \tilde\GG_n^*(\z,\x_J) := \frac{1}{\sqrt{n}} \sum_{i=1}^n  M_{N_n,i}  \{ \1_{(\Z_i^*,\X_{i,J})\leq (\z,\x_J)} -C_{\theta_0}(\z)  F_{n,J}(\x_J) \}. $$

\mds

Actually, the distance between $\bar \GG_n^*$ and $ \tilde\GG_n^*$ is negligible. Indeed, for every $(\z,\x_J)$,
$$ \Delta_n (\z,\x_J):= (\tilde\GG_n^* - \bar\GG_n^*)(\z,\x_J) = \frac{1}{\sqrt{n}}\sum_{i=1}^n (M_{N_n,i} - M_{n,i}) \{\1_{(\Z_i^*,\X_{i,J})\leq (\z,\x_J)} -C_{\theta_0}(\z)  F_{n,J}(\x_J) \}$$
is centered. Moreover, by independence between the observations and by the resampling scheme, we have
\begin{eqnarray*}
\lefteqn{
\EE[\| \Delta_n \|_{\infty}^2]=
\EE[\sup_{\z,\x_J}\Delta_n^2(\z,\x_J)]
\leq
\frac{1}{n^2}\sum_{i,j=1}^n \EE[|(M_{N_n,i} - M_{n,i})(M_{N_n,j} - M_{n,j})|]  }\\
& \leq & \frac{1}{n^2}\sum_{i,j=1}^n \EE[(M_{N_n,i} - M_{n,i})^2]^{1/2} \EE[(M_{N_n,j} - M_{n,j})^2]^{1/2} \\
&\leq &  \EE[(M_{N_n,1} - M_{n,1})^2],
\end{eqnarray*}
because the sequence $(M_{N_n,i} - M_{n,i})_{i=1,\ldots,n}$ is exchangeable.
Given $N_n=k$, the $i$-th variable $|M_{N_n,i}- M_{n,i}|$ is binomial with the parameters $(|k-n|,1/n)$, i.e.
$$ P( |M_{k,i} - M_{n,i}|=l)=C_{|k-n|}^l \frac{1}{n^l} \left(1-\frac{1}{n}  \right)^{|k-n|-l} ,\; l=0,\ldots,|k-n|.$$
Therefore, we obtain
$$ \EE[(M_{N_n,i} - M_{n,i})^2] = \sum_{k=0}^{\infty} \exp(-n)\frac{n^k}{k!} \left\{ \frac{|k-n|}{n}(1-\frac{1}{n}) + \left( \frac{|k-n|}{n} \right)^2 \right\}.$$
Simple calculations provide
$$ \sum_{k=0}^{\infty} \exp(-n)\frac{n^k}{k!}  \frac{|k-n|}{n} =  \frac{2 n^n}{n!}\exp(-n) \sim \left( \frac{2 }{\pi n}\right)^{1/2},$$
by Stirling's formula, and
$$ \sum_{k=0}^{\infty} \exp(-n)\frac{n^k}{k!}  \left( \frac{k-n}{n} \right)^2 = \frac{\exp(-n)}{n^2} \sum_{k=0}^{\infty} \frac{n^k}{k!}  (k(k-1) + k(1-2n) + n^2)= \frac{1}{n}\cdot$$
We deduce $ \EE[(M_{N_n,i} - M_{n,i})^2]=O(n^{-1/2})$ and $ \PP\left( \|\Delta_n \|_{\infty} > \epsilon\right) \longrightarrow 0,$
when $n$ tends to the infinity, given almost all sequences $\Sc_n:=(\Z_i,\X_{i,J})_{i=1,\ldots,n}$.
This means that we can safely replace $\bar\GG_n^*$ by $\tilde\GG_n^*$, and the theorem follows if we prove the weak convergence of $(\bar\GG_n,\tilde\GG_n^*)$.

\mds

Note that we can rewrite
\begin{eqnarray*}
\lefteqn{\tilde\GG_n^*(\z,\x_J) =
\frac{1}{\sqrt{n}}\sum_{i=1}^n (M_{N_n,i}-1)\left\{\1(\Z_i^*,X_{i,J})\leq (\z,\x_J) - C_{\theta_0}(\z) F_J(\x_J)\right\}  }\\
&+& \frac{1}{\sqrt{n}}\sum_{i=1}^n \left\{\1(\Z_i^*,X_{i,J})\leq (\z,\x_J) - C_{\theta_0}(\z) F_J(\x_J)\right\}
- C_{\theta_0}(\z) \sqrt{n}\left( F_{n,J} - F_J \right)(\x_J) \\
&+& \left(1- \frac{1}{n}\sum_{i=1}^n M_{N_n,i}\right) C_{\theta_0}(\z) \sqrt{n}\left( F_{n,J} - F_J \right)(\x_J) \\
&:=& \tilde\GG_{n,1}^*(\z,\x_J) + \tilde\GG_{n,2}^*(\z,\x_J) - \GG_{n,3}(\z,\x_J) + R_n(\z,\x_J) .
\end{eqnarray*}
Obviously, the last remaining term is $o_P(1)$ uniformly w.r.t. $(\z,\x_J)$, and it can be forgotten.
Moreover, since the variables $(M_{N_n,i}-1)_{i=1,\ldots,n}$ are i.i.d., centered with variance one and independent of the data, we can invoke some
multiplier bootstrap results. Consider we live in the space $\Wc:=[0,1]^{p}\times  [0,1]^{p}\times \RR^{d-p}$ that is related to our observations $W_i:=(\Z_i,\Z_i^*,\X_{i,J})$, $i=1,\ldots,n$. The true distribution of $W_i$ under the null is $P_W$, whose cdf is $C_{\theta_0} \otimes C_{\theta_0} \otimes F_J$.
Applying Corollary 2.9.3. in van der Vaart and Wellner (1996), the sequence of processes
$$ (\WW_n,\WW_n^*):=\left(n^{-1/2}\sum_{i=1}^n (\delta_{W_i} - P_W) , n^{-1/2}\sum_{i=1}^n (M_{N_n,i}-1)(\delta_{W_i} - P_W)    \right)$$
 converges weakly in $\ell^{\infty}(\Fc)\times \ell^{\infty}(\Fc)$ to a vector of independent Gaussian processes,
 where $\Fc$ denotes any Donsker class of measurable functions from $\Wc$ to $\RR$.

\mds

Now, let us consider the class $\Fc$ of functions
$$f_{\z_0,\z_0',\x_{J,0}}:(\z,\z',\x_J)\mapsto \1(\z \leq \z_0,\z' \leq \z_0',\x_J\leq \x_{J,0}),$$
for any triplet $(\z_0,\z_0',\x_{J,0})$ in $[0,1]^{p}\times  [0,1]^{p}\times \bar\RR^{d-p}$.
Note that $\Fc$ is Donsker, that $\tilde\GG_{n,1}^*(\z,\x_J)=\WW_n^* f_{1,\z,\x_J}$, $\tilde\GG_{n,2}^*(\z,\x_J)=\WW_n f_{1,\z,\x_J}$ and
that $\tilde\GG_{n,3}^*(\z,\x_J)=C_{\theta_0}(\z)\WW_n f_{1,1,\x_J}$.
By the permanence of the Donsker property (see Section 2.10 in van der Vaart and Wellner, 1996), and the continuity of $C_{\theta_0}$, the process $\tilde\GG_n^*$ converges in
$\ell^{\infty}([0,1]^p \times \RR^{d-p})$ to a gaussian process $\bar \AA^{\perp}$. Obviously, $\bar\GG_n$ tends in distribution in $\ell^{\infty}([0,1]^p \times \RR^{d-p})$ to a Gaussian process $\bar \AA$, whose
covariance function is given by
$$  \EE\left[ \bar\GG_n(\z,\x_J) \bar\GG_n(\z',\x'_J)   \right]
= C_{\theta_0}(\z\wedge\z') F_J(\x_J\wedge\x_J') - C_{\theta_0}(\z)F_J(\x_J) C_{\theta_0}(\z')F_J(\x'_J),$$
for every $\z,\z',\x_J,\x_J')$.
By some standard calculations, we check that $\EE[\tilde\GG_n^*(\z,\x_J) \tilde\GG_n^*(\z',\x'_J)]=\EE[\bar\GG_n(\z,\x_J) \bar\GG_n(\z',\x'_J)]$ for every couples $(\z,\x_J)$ and $(\z',\x_J')$, implying that $\bar \AA$ and $\bar \AA^{\perp}$ have the same covariance functions.
Moreover, the two limiting processes $\bar \AA$ and $\bar \AA^{\perp}$ are uncorrelated because
$$ \EE[\bar\GG_n (\z,\x_J) \tilde\GG_n^*(\z',\x'_J)]= \EE[\bar\GG_n (\z,\x_J) \EE[ \tilde\GG_n^*(\z',\x'_J) | \Sc_n]]=0,$$
for every couples $(\z,\x_J)$ and $(\z',\x_J')$.
Therefore, the $\bar \AA$ and $\bar \AA^{\perp}$ are two independent versions of the same Gaussian process.

\begin{rem}
\label{cov_boot_v2}
If there were no resampling of the observations $\X_{i,J}$ at the first level, this would no longer be true. Indeed, the corresponding bootstrapped process would be
given by
$$ \GG_n^{**}(\z,\x_J):=\frac{1}{\sqrt{n}}\sum_{i=1}^n \1(\X_{i,J}\leq \x_J) \left\{ \1(\Z_i^* \leq \z) - C_{\theta_0}(\z)   \right\},$$
implying
$$ \EE\left[ \bar\GG^{**}_n(\z,\x_J) \bar\GG^{**}_n(\z',\x'_J)   \right]
= F_J(\x_J\wedge\x_J')[ C_{\theta_0}(\z\wedge\z') - C_{\theta_0}(\z)C_{\theta_0}(\z')],$$
that is different of $\EE\left[ \bar\GG_n(\z,\x_J) \bar\GG_n(\z',\x'_J)   \right]$.
\end{rem}

\mds

To conclude, we apply Corollary 1.4.5. in van der Vaart and Wellner (1996): for every bounded nonnegative Lipschitz function $h$ and $\tilde h$,
 \begin{eqnarray*}
\lefteqn{ \EE[ h(\bar\GG_n) \tilde h (\tilde \GG_n^*)] - \EE[h(\bar \AA)\tilde h(\bar \AA^{\perp})]= \EE[ h(\bar\GG_n) \left( \EE[\tilde h (\tilde \GG_n^*) |\Sc_n] - \EE[\tilde h(\bar \AA^{\perp})]\right)] }\\
&+& \EE[ \left( h(\bar\GG_n) - \EE[h(\bar \AA)] \right)]  \EE[\tilde h(\bar \AA^{\perp})]. \hspace{5cm}
\end{eqnarray*}
The first (resp. second) term tends to zero by the weak convergence of $\tilde\GG_n^*$ (resp. $\bar\GG_n$). This concludes the proof.
$\Box$

\end{document}